\documentclass[11pt]{amsart}

\usepackage{amssymb}

\newtheorem{propo}{Proposition}[section]
\newtheorem{defi}[propo]{Definition}

\newtheorem{lemma}[propo]{Lemma}
\newtheorem{corol}[propo]{Corollary}

\newtheorem{theo}[propo]{Theorem}

\newcommand{\diag}{\mathop{\rm diag}\nolimits}

\newcommand{\Id}{\mathop{\rm Id}\nolimits}

\newcommand{\Irr}{\mathop{\rm Irr}\nolimits}

\newcommand{\NN}{\mathop{\mathbb     N}\nolimits}

\newcommand{\al}{\alpha}

\newcommand{\ep}{\varepsilon}

\newcommand{\lam}{\lambda }
\newcommand{\ld}{,\ldots ,}
\newcommand{\si}{\sigma}
\newcommand{\om}{\omega}
\newcommand{\ra}{ \rightarrow }

\newcommand{\up}{^{-1}}

\newcommand{\lan}{\langle}
\newcommand{\ran}{\rangle}
\newcommand{\el}{\end{lemma}}

\def\d12{{_{12}}}

\def\acf{{algebraically closed field }}

\def\au{{automorphism }}

\def\ei{{eigenvalue }}
\def\eis{{eigenvalues }}

\def\f{{following }}

\newcommand{\ii}{if and only if~\,}

\def\ir{{irreducible }}

\def\irr{{irreducible representation }}

\def\irrs{{irreducible representations }}
\def\itf{{It follows that }}

\def\mult{{multiplicity }}

\def\pe{{permutation }}
\def\po{{polynomial }}
\def\pos{{polynomials }}

\def\rep{{representation }}
\def\reps{{representations }}

\def\syl{{Sylow $p$-subgroup }}

\newcommand{\med}{\medskip}

\begin{document}

\title{Almost cyclic  elements in Weil
representations of finite classical  groups}
\author{Lino Di Martino and A.E. Zalesski}

\maketitle

\medskip

\centerline{ Dedicated to Otto H. Kegel on the occasion of his 80th birthday}

\bigskip

\bigskip
Abstract.  This paper is a significant part of a general project aimed to classify all irreducible \reps
of finite quasi-simple groups over an algebraically closed field, in which  the
 image of at least one element is represented by
an almost cyclic matrix (that is, a square matrix $M$ of size $n$ over a field $F$ with the property that
there exists $\al\in F$ such that $M$ is similar to
 $\diag(\al\cdot \Id_k, M_1)$, where $M_1$ is cyclic and $0\leq k\leq n$).
The paper focuses on the Weil \reps of
finite classical groups, as there is strong evidence that these representations  play a key role in the general picture.

\section{Introduction}

Let $V$ be a vector space of finite dimension $n$ over an arbitrary field $F$, and let $M$ be a square matrix of size $n$ over $F$. Then 
$M$ is said to be cyclic if the characteristic \po and the minimum \po of $M$ coincide. Note that a matrix $M\in {\rm Mat}\,(n,F)$ is cyclic \ii the
$F\langle M \rangle$-module $V$ is cyclic, that is, is generated by
a single element. This is  standard terminology in module theory, and
the source of the term `cyclic matrix'. Matrices with simple spectrum
often arising in applications are cyclic.  We consider a generalization of the notion of cyclic matrix, namely, we define a matrix $M\in {\rm Mat}\,(n,F)$ to be almost cyclic if
there exists $\al\in F$ such that $M$ is similar to
 $\diag(\al\cdot \Id_k, M_1)$, where $M_1$ is cyclic and $0\leq k\leq n$.

Examples of almost cyclic matrices arise naturally in the study
of matrix groups over finite fields. For instance, if an element
$g\in GL(V)$ acts irreducibly on $V/V'$, where $V'$ is some eigenspace of
$g$ on $V$,  then $g$ is almost cyclic. Reflections and transvections 
 are important examples. Other relevant examples  are provided by unipotent matrices with
Jordan  form consisting of a single non-trivial block.

Possibly, the strongest motivation to study groups
containing an almost cyclic matrix  is to contribute to the
recognition of linear groups and finite group representations
by the property of a single matrix.  Our main inspiration is
a paper by Guralnick, Penttila, Praeger
and Saxl \cite{GPPS}, in which the authors
classified irreducible linear groups over finite fields generated by `Dempwolff elements'.
If $V=V(n,q)$ is an $n$-dimensional vector space over a finite field of order $q$,  and  $G=GL(V)=GL(n,q)$, we say that $g \in G$ is a  Dempwolff element if
$|g|=p$ for some prime $p$ with $(p,q)=1$ and $g$ acts irreducibly on $V^g:=(\Id -g)V$.
U. Dempwolff in \cite{De} initiated the study of subgroups of $GL(n,q)$ generated by such elements, obtaining a number of valuable results. The main restriction in \cite{De} is the assumption that $2\dim V^g>\dim V$, and this assumption is held in  \cite{GPPS}.
Clearly, Dempwolff elements are almost cyclic (and are reflections if $p=2$).
We have realized that, if one wishes to drop this restriction, and furthermore  obtain  satisfactory results in full generality, a more conceptual approach is available. Namely, one should deal with finite linear groups over an algebraically closed field. Therefore, our general program can be stated as follows:  determine all irreducible finite linear groups over an algebraically closed field, which are
generated by almost cyclic matrices. In addition, we wish to relax the assumption, held in \cite{GPPS}, that $g$ is of prime order.
However, as current applications seem to
focus on $p$-elements, we will limit ourselves to the study of elements $g\in G$ of any $p$-power order.

Since we will make a systematic use of representation theory and will exploit the classification theorem of finite simple groups, a key part of our project necessarily  focuses on finite quasi-simple groups. The sporadic simple groups and their covering groups have been completely dealt with in  \cite{DMPZ}.  
 In \cite{DMZ10},  we started to deal with finite groups of Lie type, and determined all the  \ir \reps
of a quasi-simple group of Lie type $G$ over an algebraically closed field $F$ of characteristic coprime to the defining characteristic of $G$, in which the image of at least one unipotent
element $g$ is represented by an almost cyclic matrix.  The complementary case, when $g$ is unipotent in $G$, and the characteristic of $F$ is the defining characteristic of $G$, has been settled for classical groups by Suprunenko in \cite{Su13}. 
This leaves open the case when $G$ is an exceptional group of Lie type, as well as the general case when $g$ is a semisimple element of prime-power order of $G$.

The present paper focuses on Weil \reps of
finite classical groups (an overview of these representations is given in Section $5.1$).  The reason to treat this case separately is that there is strong evidence that
most examples of semisimple almost cyclic elements 
occur  in  Weil representations.
Furthermore, the study of Weil \reps requires a lot of analysis and technical background
which justifies the choice of treating them in an independent paper. Besides, Weil \reps
play a very significant role in the \rep theory of classical groups, and several features and properties
of them have been the subject of intensive study in many recent papers. So, the present paper can also be viewed as a contribution to this research area.

Before stating the main result,
a few more words are needed about the existing literature.
Before Dempwolff's work in \cite{De}, important results related to our problem had already appeared
in the literature.

Ch. Hering (\cite{He1},  \cite{He2}) essentially classified  the
finite \ir subgroups $G$ of $GL(n,q)$ containing an irreducible
element of prime order, provided $G$ has a composition factor
isomorphic to a group of Lie type (or an alternating group). Also, the finite \ir
linear groups generated by transvections, reflections and pseudo-reflections
were classified by A. Wagner in \cite{W2,Wag}, Pollatsek \cite{PH},
A.E Zalesski and V.N. Serezhkin in \cite{ZS,ZS80}. I.D. Suprunenko and A.E. Zalesski
in \cite{SZ98, SZ00}
classified the \ir \reps of Chevalley groups in the natural characteristic containing
a matrix with simple spectrum. Furthermore, Di Martino and Zalesski in \cite{DMZ0, DMZ8},
following an earlier paper by Zalesski \cite{Z99}, studied the minimum \pos
of elements of prime power order  in
the cross characteristic \reps of classical groups. This was further extended  by Tiep and Zalesski
in \cite{TZ08}. The latter work also extends  part of the results  of
\cite{Za06,Z08} to
representations over fields of prime characteristic.


More information is available in the case where the ground field $F$
is of characteristic zero. Huffman and Wales in \cite{Hu} classified
the finite \ir linear groups
generated by elements $g$ such that $\dim V^g\leq 2$. As a particular case, this result contains a classification of finite \ir linear groups
over the complex numbers generated by almost cyclic elements of order $3$.  Zalesski in \cite{Za06} determined the \ir linear groups over the complex numbers generated by  elements $g$ of prime order $p>3$ that  have at most  $p-2$ distinct eigenvalues.    In addition, in \cite{Z08}
Zalesski determined  the \ir \reps of quasi-simple groups in which an element of prime order $p$
has at most $p-1$ distinct eigenvalues.  Another relevant work for the characteristic 0 case is \cite{RZ}. See also the surveys \cite{TZ00,Z09}
for further  details.

\medskip

\medskip
Now, we state our main result.

\begin{theo}\label{mt1}

 Let $G$ be one of the following groups:  $G=Sp(2n,q)$, where $n>1$ and $q$ is odd;  $SU(n,q)\subseteq G\subseteq U(n,q)$,
 where  $n>2$;  $SL(n,q)\subseteq G\subseteq GL(n,q)$, where $n>2$. Let  $g\in G$ be a non-scalar $p$-element, where $p$ is a prime not dividing $q$. Let $F$ be an \acf of characteristic $\ell$ not dividing $q$, and
$\tau$
an \ir Weil $F$-representation of $G$. Then the matrix of $\tau (g)$ is  almost
cyclic if and only if one of the following occurs:

\medskip

$(1)$ $G= Sp(2n,q)$, and either

$(a)$ $n$ is a $2$-power and $|g|=(q^n + 1)/2$ is odd, or

$(b)$ $q=3$, $n \neq p$ is an odd prime, and $|g|=(3^n - 1)/2$ is odd, or

$(c)$  $n=2$, $q=3$, and one of the following holds:

$(c_1)$ $p=2$, $\ell\neq 2$ and either    $|g|=2$ 
and $\dim\tau=5$, or  $|g|=4$,  $g^2\notin  Z(G)$   and $\dim\tau=4$, or
 $|g|=8$,  $g^4\in Z(G) $ and $\dim\phi=4$ or $5$;


$(c_2)$ $p=5$ and $\dim\tau \in \{4,5\}$,  where $\dim\tau\neq 5$
 if $\ell= 2$; 

$(c_3)$ $p=\ell=2$ and $\dim\tau=4$. In addition, either $|g|=4$ or $|g|=2$. 

\medskip

$(2)$ $SU(n,q)\subseteq G\subseteq U(n,q)$, and either

$(a)$ $|g|=(q^n+1)/(q+1)$ and $n\neq p$ is an odd prime greater than $3$, or 

$(b)$ $(n,q)=(5,2)$, $|g|=9$, $\ell \neq 3$ and $\dim\tau=10$;

$(c)$ $(n,q)=(4,2)$  and one of the following holds:

$(c_1)$  $|g|=3$ or $9$;



$(c_2)$ $|g|=5$; or

$(d)$ $(n,q)=(3,3)$, and one of the following holds:

$(d_1)$ $|g|=7$;

$(d_2)$  $|g|=8$ and either $\dim\tau=6$, or $\ell\neq 2$ and $\dim\tau=7$; or

$(e)$ $(n,q)=(3,2)$, $|g|=3$ or $9$ and $\dim\tau= 2,3$ for $\ell \neq 2$, $\dim\tau=3$ for $\ell =2$.

\medskip

$(3)$ $SL(n,q)\subseteq G\subseteq GL(n,q)$,  and either

$(a)$ $G=SL(n,2)$, where  $n\neq p$ is an odd prime and $|g|=2^n - 1$ is a Mersenne prime, or

$(b)$ $|g|=(q^n-1)/(q-1)$, where $q>2$, and $n\neq p$ is an odd prime.

\end{theo}

For the sake of completeness we have also examined in this paper  the case where  $SL(2,q)\subseteq G\subseteq GL(2,q)$,  without assuming that $\tau $ is Weil. 
Note that the degree of an irreducible $F$-representation of $G$ in this case belongs to the set $\{1,
q-1,q,q+1,(q-1)/2, (q+1 /2)\}$, where in the last two cases $q$ is
odd. 

The results obtained are collected in Theorem 1.2 below. Additionally, these results (as a consequence of Lemma $\ref{un2}$ and Corollary $\ref{un3}$ in Section $4$) can be carried over to any group $G$ such that  $SU(2,q)\subseteq G\subseteq U(2,q)$.

\begin{theo}\label{mt2}  Let  $SL(2,q)\subseteq G\subseteq GL(2,q)$, $q>3$, and let  $g\in G$ be a non-scalar $p$-element, where $p$ is a prime not dividing $q$. Let $F$ be an \acf of characteristic $\ell$ not dividing $q$, and let $M$ be an irreducible $FG$-module with $\dim M >1$, affording the representation $\tau$. The following holds:

\medskip

$(1)$ Suppose $p>2$. Then $\tau(g)$ is almost cyclic if and only if $ \dim M  \leq |g|+1$. (In this case, $(2,q+1)|g| = q \pm 1$).

\medskip

$(2)$ Suppose $p=2$, and let $h$ denote the projection of $g$ into $G/Z(G)$.

Assume first that  $q\equiv 1\pmod 4$.
 Then  $\tau(g)$ is almost cyclic if and only if  
 one of the following occurs:
 
 $(i)$ $\ell \neq 2$, $g \in SL(2,q) \cdot Z(GL(2,q))$,  $\dim M=(q\pm 1)/2$ and $|h|= (q-1)/2$;

$(ii)$ $\ell \neq 2$,  $g \notin SL(2,q) \cdot Z(GL(2,q))$,  and either $\dim M= q$ or $q -1$ and  $|h|=q-1$, or $G=GL(2,5) \cong \tau(G)$, $\dim M = 4$, $|g|=8$ and $|h|=2$;

$(iii)$ $\ell = 2$ and either $\dim M\leq |h|+1$ or $G=GL(2,5)$, $\tau(G)\cong O^-(4,2)$, $\dim M=4$ and $\tau(g)$ is a transvection.

 Next, assume that $q\equiv -1\pmod 4$.  
 
 Then  $\tau(g)$ is almost cyclic if and only if  
 one of the following occurs:
 
$(i)$ $\ell \neq 2$, $g \in SL(2,q) \cdot Z(GL(2,q))$,  $\dim M=(q\pm 1)/2$ and $|h|= (q+1)/2$;

$(ii)$ $\ell \neq 2$,  $g \notin SL(2,q) \cdot Z(GL(2,q))$,  $\dim M= q$ or $q \pm 1$ and  $|h|=q+1$;

$(iii)$ $\ell=2$, and one of the following holds:

$a)$ $\dim M=q \pm 1$ and $|h|=q+1$ (here the case $\dim M=q + 1$ only occurs for $g \notin SL(2,q) \cdot Z(GL(2,q))$);

$b)$ $\dim M= (q-1)/2$ and $|h|= (q+1)/2$;

$c)$ $q=7$, $\dim M =3$ and $|h| = 2$.

\end{theo}

\medskip

\title{ NOTATION}

\medskip

Throughout the paper, unless  stated otherwise, we denote by  $F$  an \acf
of characteristic  $\ell$.

For any finite group $G$, the representations of $G$ we consider in the paper are all over $F$, unless stated otherwise.
We write  $1_G$ for the trivial $FG$-module  
and  $\rho_G^{reg}$ for the regular $FG$-module (that is the free $FG$-module of rank 1).  

If $G$ is a finite group of Lie type of defining characteristic $r$,
we always assume that $\ell$ is coprime to $r$.  

For the reader's sake, it is also convenient to lay down explicitly some of the notation which is used
throughout the paper for finite classical groups.

Let $V$ be a vector space of finite dimension $m>1$ over a field $K$.

If $K$ is a finite field of
order $q$ (where $q$ is a power of a prime $r$), $K$ will be usually denoted by $\mathbb{F}_{q}$,
and the general linear group $GL(V)$ and the special linear group $SL(V)$ will be denoted
by $GL(m,q)$
and $SL(m,q)$, respectively.

Suppose that the space $V$ is endowed with a non-degenerate orthogonal,
symplectic or unitary form. Then $I(V)$ will denote the group of the isometries of $V$, and we will loosely
use the term 'finite classical group' for a subgroup $G$ of $I((V)$ containing $I(V)^{\prime }$.  In particular:
if $V$ is a symplectic space over  $\mathbb{F}_{q}$,  $I(V)$ will be denoted by $Sp(m,q)$; if $V$ is a unitary
space over the field $\mathbb{F}_{q^{2}}$, $I(V)$ will be denoted by $U(m,q)$; and if $V$ is an orthogonal
space over  $\mathbb{F}_{q}$, $I(V)$ will be denoted by $O(m,q)$. 
 It should be noted that in places the term 'classical group'
will be meant to include also the groups $GL(m,q)$ and $SL(m,q)$ (considering $V$ endowed with the identically
zero bilinear form). Finally, at times we will need to consider, for a given classical group $G$ (defined as above)
the corresponding central quotient (projective image), which will be denoted by $PG$.

Finally, we mention that the notation used in the paper for objects of general group theory is fairly standard. E.g., for a group $G$, $Z(G)$ denotes the centre of $G$; for a subgroup $H$ of $G$, $N_G{(H)}$ and  $C_G{(H)}$ denote the normalizer and the centralizer of $H$ in $G$, respectively. Similarly, for $x \in G$,  $C_G{(x)}$ denotes the centralizer of $x$ in $G$. And so on.  

\section{Preliminaries}

For the  reader's convenience we recall the \f definition:

\begin{defi} Let $M$ be an $(n\times n)$-matrix over an arbitrary field $K$.
We say that $M$ is almost
cyclic  if there exists $\al\in K$ such that $M$ is similar to
 $\diag(\al\cdot \Id_k, M_1)$, where $M_1$ is cyclic and $0\leq k\leq n$.
\end{defi}

{\bf Remark 1}.
In the definition above, it has to be understood that for $k=0$ the matrix $M=M_1$ is cyclic, whereas for $k=n$ the matrix $M$ is scalar.

\medskip

{\bf Remark 2}.
Let $\overline{K}$ denote the algebraic closure of
$K$, and for $\lambda \in \overline{K}$ denote by $\lambda J$ a
Jordan block with eigenvalue $\lambda $. Observe that a matrix
$M_{1}$ is cyclic if and only if $M_{1}$ has Jordan form
$\mathrm{diag}(\lambda _{1}J_{1},....,\lambda _{s}J_{s})$, where the
$\lambda _{j}$'s, $1\leq j\leq s$, are pairwise distinct. In
particular, suppose that $M$ has  order $p^{a}$, where $p$ is a prime,  and set $\ell =\,$char$\, K$. Then $M$
is almost cyclic if and only if the eigenvalues of $M_{1}$ in  $\overline{K}$ are pairwise
distinct when $\ell \neq p$, and if and only if $M_{1}$ consists of a single Jordan
block when $\ell =p$.

\medskip
An elementary  observation, which will be useful throughout the paper,  is the following:
if $M\in GL(V)$ is almost cyclic, and $U$ is
an $M$-stable subspace of $V$, then the induced action of  $M$ on $U$ and
on $V/U$ yield almost cyclic matrices.

\medskip
Let us denote by $\deg (X)$ the degree of the minimum \po of a square matrix $X$ over a field $F$.
Then the \f holds:

\begin{lemma}\label{p74}  Let $A,B$ be  non-scalar
square matrices over an arbitrary field $K$, both diagonalizable over $K$, and
let $k=\,{\rm deg}\,(A),l=\,{\rm deg}\,(B)$.  
 Suppose that
$A\otimes B$ is almost cyclic. Then $A,B$ are cyclic and
 $\deg (A\otimes B)\geq kl-\min\{k,l\}+1$.

\end{lemma}

Proof. The claim about the cyclicity of $A$ and $B$ is obvious. Assume
$k\leq l$ and let $\ep_1\ld \ep_k, \eta_1\ld \eta_l$
be the \eis of $ A,B$, respectively. If $A\otimes B$
is   cyclic, then $\deg (A\otimes B)=kl$. Suppose
$A\otimes B$ is not cyclic. We can assume that
 $\lam=\ep_1\eta_1$ is an \ei of $A\otimes B$
 of \mult greater than 1.  Then all the $\ep_i\eta_j$'s such
that $\ep_i\eta_j\neq \lam$   are distinct. The number of
pairs $(i,j)$ such that $\ep_i\eta_j= \lam$ is at most $k$, and hence  $A\otimes B$ has
at least $1+(kl-k)$ distinct eigenvalues, as required.

\medskip

{\bf Remark}.  A typical application of Lemma \ref{p74}  is the following. Let
$X=X_1\times X_2$ be the direct product of two groups $X_1,X_2$ and $g\in X$
be a $p$-element for some prime $p$. Then $g=g_1g_2$, where $g_1\in X_1,g_2\in X_2$.
Let $\phi\in\Irr_F X$, where $F$ is an \acf of characteristic different from $p$. Then $\phi=\phi_1\otimes \phi_2$, where $\phi_i\in\Irr X_i $
for $i=1,2$. Suppose that  $\phi_i(g_i)$ has order $p^{m_i} > 1$ modulo the scalars, and let
$\phi_i(g_i^{p^{m_i}})=\lam_i\cdot \Id$ for some $\lam_i\in F$. We may assume
$m_1\geq m_2$. Clearly, every eigenvalue of $\phi_i(g_i)$ is a $p^{m_i}$-root
of $\lam_i$. Furthermore, we have $\phi(g)=\phi_1(g_1)\otimes \phi_2(g_2)$.
It follows that the \eis of $\phi(g)$ are $p^{m_1}$-roots of $\lam_1\lam_2^{p^{m_1-m_2}}$,
and the minimum \po of
$\phi(g)$ is of degree at most $p^{m_1}$. Lemma \ref{p74} tells us that if $k,l$
are the degrees of the minimum polynomials of
$A=\phi_1(g_1)$ and
$B=\phi_2(g_2)$, respectively,  then
$\phi(g)$ is not almost cyclic unless (a) $A,B$ are cyclic matrices and (b)
$p^{m_1}\geq kl-\min\{k,l\}+1$.

\begin{lemma}\label{tp1} Let $\lam,\mu$ be two completely reducible \reps of a cyclic 
$p$-group $X=\lan x \ran$ of  order $p^a$ over a field $K$, and let $l$ and $k$, where 
$l\geq k>1$, be
the degrees of the minimum \pos  of $\lam (x),\mu (x)$, respectively.
Suppose that $k+l>p^a>3$. Then  $\lam (x)\otimes \mu (x)$ is not
almost cyclic.\end{lemma}

Proof. 
Suppose the contrary.
Then, by Lemma \ref{p74}, $p^a\geq \deg (\lam (x)\otimes \mu
(x))\geq k(l-1)+1$. As $k+l>p^a$, we have $l\geq \frac{p^a+1}{2}$. So $p^a\geq k(\frac{p^a+1}{2} - 1)+1= k(\frac{p^a -1}{2}) + 1=p^a+(k-2)\frac{p^a-1}{2}$, whence $k=2$. This implies $2+l > p^a\geq 2l-1$, whence $2+l>2l-1$, that is $l=2$. This in turn forces $p^a<4$. A
contradiction, as $p^a>3$ by assumption.

\begin{lemma}\label{p99} Let $K$ be
a  field of arbitrary characteristic $\ell$,
and let $J_m,J_n$ be unipotent Jordan blocks of size $m\geq n>1$ over $K$.
Then $J_m\otimes J_n$  is almost cyclic \ii $m=n=2$ and
 $\ell\neq 2$.

In particular, if $\ell>0$, $P=\lan g\ran$ is a cyclic $\ell$-group,
and $M,N$ are non-trivial indecomposable $KP$-modules, then the matrix of $g$ on $M\otimes N$
is almost cyclic \ii $M$ and  $N$ are of dimension $2$ and $\ell\neq 2$.
\end{lemma}

Proof. Let 
$ V_m$ and $V_n$ be vector spaces over $K$ on which $J_m$ and $J_n$ 
act, respectively.
Clearly, for every $1\leq i\leq m$ there is exactly one subspace $V_i$ of dimension $i$ in $V_m$, which is stable
under $J_m$. Similarly, for every $1\leq j\leq n$ there is a single subspace $V_j$
of $V_n$ stable under $J_n$. Moreover, each $V_i$ is indecomposable under the action of $J_m$.
Similarly for each $V_j$.
Now,  $J_m\otimes J_n$ acts on $V_m\otimes V_n$, and stabilizes 
each subspace $V_i\otimes V_j $.  
In order to prove the first part of the statement, it is enough to prove
that $J_m\otimes J_n$  is almost cyclic if $m=n=2$ and $\ell\neq 2$, whereas  it
is not almost cyclic if $m=3$, $n=2$ or $n=m=\ell=2$.

A direct computation shows that, if $m=n=2$, then $V_2\otimes V_2=W_1\oplus W_2$,
where $W_1,W_2$ are $(J_m\otimes J_n$)-stable subspaces and $\dim W_1=\dim W_2=2$ if $\ell=2$, whereas $\dim W_1=3$ and $\dim W_2=1$ if $\ell \neq 2$.

Next, let
 $m=3,n=2$. In this case we do not need to consider $\ell=2$.
A direct computation shows  that $V_3\otimes V_2=W_1\oplus W_2$,
where $W_1,W_2$ are $(J_m\otimes J_n$)-stable subspaces, and $\dim W_1=\dim W_2=3$ if $\ell=3$,
whereas $\dim W_1=4$ and $\dim W_2=2$ if $\ell \neq 3$.

The additional claim of the lemma is a module-theoretic version
that follows straightforwardly from the first part.

\begin{lemma}\label{higm}  
Let $T=RH$  be  a  finite  group where  $H=\langle h\rangle$ is a
 cyclic $p$-subgroup  and $R$ is a normal $r$-subgroup for some prime
$r\neq p.$ Let   $|H/C_H(R)|=p^k$. Let  $\phi$ be  an $F$-\rep of
$T$ faithful on $R.$  Suppose that $(\ell ,r)=1$ and    $1< \deg \phi
(h)<m(h)$ where $m(h)$ is the order of $h$ modulo $Z(H)$. Then $R$ is non-abelian and $p^a=r^b+1$ for some
$a,b\in \NN .$ Additionally,   $\deg \phi (h)\geq
(p^{a}-1)p^{k-a}$.\end{lemma}

For $\ell =p>0$ the proof can be  found, for  instance,  in
\cite[VII.10.2]{Fe}. Observe  that  a  faithful ${\Bbb C}R$-module
remains faithful under reduction modulo $p,$ and the degree of the
minimum \po cannot increase.  So  Lemma \ref{higm}    is   valid
for characteristic 0. Using reduction modulo $\ell \neq r$, one obtains  the
result for $\ell \neq p$, as the character  of  $H$ coincides on
$\ell '$-elements with the Brauer character.

\begin{lemma}\label{zgm}     {\rm \cite[Ch. IX,
Lemma 2.7]{HB}} Let $p,r$ be primes and $a,b$ positive integers
such that $p^a=r^b+1$. Then either $p=2$, $b=1$, or $r=2, a=1$ or
$p^a=9$.
\end{lemma}

\medskip

Recall that an element $g$ of a group of Lie type $G$ of defining characteristic $\ell$ is said to be semisimple if $g$ has order coprime to $\ell$. Furthermore, we will say that $g$ is regular semisimple if its centralizer in $G$ has order coprime to $\ell$ (this definition, convenient in our context, is well-known to be equivalent to that usually given in the context of algebraic groups).

The series of results that follow will be crucial for our purposes.

\begin{lemma}\label{d102}  {\rm (Gow \cite{Gow})} Let $G$ be
 a quasi-simple group of Lie
type in characteristic $r$ and $g\in G$. Suppose that
$(|C_{G}(g)|,r)=1$, that is, $C_{G}(g)$ contains no element of order
$r$. Then every semisimple element of $G$ can be factorized as
  $ab$, where $a,b\in g^{G}$. \el

\begin{lemma}\label{g12}   Let  $G$ be
    a quasi-simple group of Lie
type.  Then the following holds:

$(1)$ $G$ can be generated by two semisimple elements.

$(2)$ Let $g\in G$ be a regular semisimple element. Then
$G$ can be generated by three elements conjugate to $g$. \el

Proof. 
(1) Obviously, it suffices to prove the statement for $G$ simple. Let
$r$ be the definining characteristic for $G$. If $r=2$, then the result
is available from \cite[Theorem 8.1]{GK}. So, let $r>2$. If $G$ is
classical or of type $E_{6}, {}^{2}E_{6}(q)$, then the result is
contained in the proof of Theorem 3.1 in \cite{MSW}, except for
groups $\Omega ^{\pm }(8,2)$ (which are covered by \cite{GK}) and the group
$PSU(3,3)$. (Note that the case of the groups of type $A_1(q)$ goes back to L.E. Dickson.)  The group $PSU(3,3)$ is easily dealt with: it is
generated by an elements of order $7$ and $4$ (from the class $4A$
in \cite{Atl}). Furthermore, it is shown in \cite{LM} that the
groups $G\in \{F_{4}(q), E_{6}(q), {}^{2}E_{6}(q),E_{7}(q),
E_{8}(q)\}$ are generated by a pair of elements $x,y$ with
$x^{2}=y^{3}=1$ (this is called a $(2,3)$-generation).
 Therefore, if $(6,q)=1$, we are done. If $r=3$, the result for these groups again follows from  \cite{LM}, where a $(2,3)$-generation is provided, with the additional property that $c=xy$ is  a semisimple element
of a suitable kind. Clearly, $G=\lan x,c\ran$, and $x,c$ are
semisimple.  The groups  $G\in \{ G_{2}(q), {}^{2}G_{2}(q),
{}^{3}D_{4}(q)
\}$  for $q$ odd are known to be $(2,3,7)$-generated (that is, they are generated by two elements $x$ and $y$ of order 2 and 3, respectively, such that $xy$ has order 7), with the exception of the groups  ${}^{2}G_{2}(3)$,
$ G_{2}(3)$,   ${}^{3}D_{4}(3^{n})$ (see  \cite{Mal1},
\cite{Mal2}).  So the result follows as above, apart
for the quoted exceptions.  The groups  $G_{2}(3), {}^{3}D_{4}(3^{n})$ are covered in \cite{Mal1} (see the proof of the Corollary, p. 350) and in \cite{Mal2} (see the proof of Proposition 3), respectively. Finally the simple group  ${}^{2}G_{2}(3)'$ is isomorphic to $SL(2,8)$, and hence the result follows.

(2)  By (1), every quasi-simple group of Lie type can be generated by
two semisimple elements. Therefore, if $g \in G$ is regular
semisimple,  then, by Lemma  \ref{d102},  $G$ can be generated by four elements conjugate to
$g$. It was proven in  \cite{GK} and  \cite{Sn} that for any non-trivial $g \in G$, there exists a suitable $h \in G$ such that $G = \lan g,h\ran $, and furthermore
that $h$ can be chosen to be
semisimple. It now follows from Lemma \ref{d102} that $G$ can be
generated by three conjugates of any given regular semisimple element $g$.

\medskip
The following Propositions are due to R. Guralnick and  J. Saxl (\cite{GS}).

\begin{propo}\label{GS1}
Let $G$ be 
simple group of Lie type, and let
$1 \neq  x \in  Aut(G)$.  Denote by $\al (x)$ the minimum number
of $G$-conjugates of $x$ sufficient to generate  $\lan x,G\ran$. Then the \f holds:

\medskip
$(1)$ {\rm (\cite[Theorem 4.2]{GS})} Let $G$ be a simple classical group,
and assume that the natural module for $G$
has  dimension $n \geq  5$. Then
$\al(x) \leq  n$, unless $G = PSp(n,q)$ with $q$ even, $x$ is a transvection and $\al (x) = n+ 1$.
 \medskip

$(2)$ {\rm (\cite[Theorem 5.1]{GS})} Let $G$ be a simple  exceptional group of Lie type,
of untwisted Lie rank $m$. Then $\al(x)\leq  m + 3$, except possibly for the case $ G = F_4(q)$
with $x$ an involution, where $\al(x) \leq  8$.
 \end{propo}

In the same paper, the authors prove analogous results for low-dimensional classical groups. We quote the following, which will be needed in the sequel:

\begin{propo}\label{GS2} Under the same assumptions and notation of Proposition $\ref{GS1}$, the following holds:

\medskip

$(1)$  {\rm (\cite[Theorem 4.1(a)]{GS})} If $G=PSL(3,q)$, and $x$ has prime order, then
 $\al(x)\leq 3$, unless $x$ is an involutory graph field automorphism with $\al(x) \leq 4$.

\medskip
$(2)$ {\rm (\cite[Theorem 4.1(c)]{GS})} If $G=PSL(4,q)$, $q>2$, and $x$ has prime order, then
 $\al(x)\leq 4$, unless $x$ is an involutory graph automorphism with $\al(x) \leq 6$.

\medskip
$(3)$ {\rm (\cite[Theorem 4.1(d)]{GS})} If $G=PSL(4,2)$, and $x$ has prime order, then
 $\al(x)\leq 4$, unless $x$ is a graph automorphism with $\al(x)=7$.

\medskip
$(4)$ {\rm (\cite[Lemma 3.3]{GS})} If $G=PSU(3,q)$, $q>2$, and $x$ has prime order, then
 $\al(x)\leq 3$, unless $q=3$ and $x$ is an inner involution with  $\al(x)=4$.

\medskip
$(5)$ {\rm (\cite[Lemma 3.4]{GS})} If $G=PSU(4,q)$ and $x$ has prime order, then $\al(x)\leq 4$, unless one of the following holds:

{\rm (i)} $x$ is an involutory graph automorphism and $\al(x)\leq 6$;

{\rm (ii)} $q=2$ with $x$ a transvection and $\al(x)\leq 5$.

\medskip
$(6)$ {\rm (\cite[Theorem 4.1(f)]{GS})} If $G=PSp(4,q)$ and $x$ is of prime order, then
$\al(x)\leq 4$, unless $x$ is an involution and $\al(x)\leq 5$, or $q=3$ and $\al(x)\leq 6$.

\end{propo}

The following result, which only requires elementary linear algebra, will often be applied 
in this paper in order to establish a connection between the occurrence of almost cyclic 
matrices in representations of irreducible linear groups and the generation of these groups 
by conjugates.  In particular, it will be usually combined with Lemma \ref{g12} and 
Propositions \ref{GS1} and \ref{GS2}.

\begin{lemma} \label{nd5} 
If $G< GL(n,F)$ is a finite irreducible linear group generated by $m$ almost cyclic elements 
of the same order $d$ modulo $Z(G)$, then $$n\leq m(d-1).$$
\end{lemma}

\medskip

Proof. See {\rm (\cite[Lemma 2.1]{DMPZ})}.

\medskip

Furthermore, we quote the following  result, which will be useful in the next sections.

\begin{lemma}\label{wz1}  {\rm (see \cite{PH})}
 Suppose that $charF= \ell=2$. Let $G$ be an \ir subgroup of $GL(n,F)$ generated by transvections. Then $G$ is isomorphic
either to a
symmetric group $S_{n+1}$ or $S_{n+2}$, where
$n$ is  even, or to one of the groups
$SL(n, q)$, $Sp(n,q)$, $O^\pm(n,q)$,  $SU(n,q)$, where $q$ is a $2$-power.
\end{lemma}

\medskip

We close this section by recording some results which follow from the representation theory of finite groups having cyclic Sylow $p$-subgroups for some prime $p$.

\begin{lemma}\label{ww1}  Let $G$ be a finite group
with a non-trivial cyclic  $\ell$-subgroup $P$ of order $\ell^d$, and let
$M$ be an \ir  $FG$-module faithful on $P$. Then the following holds:

$(1)$ if $M$ is of defect zero, then $M|_P=\frac{\dim M
}{|P|}\rho_P^{reg}$;

$(2)$ if $M$ is of defect $d$, then $M|_P=\frac{\dim M - \dim L
}{|P|}\rho_P^{reg}\oplus L$, where $L$ is the direct sum of isomorphic
indecomposable $FP$-modules of dimension $e<|P|$.
In addition, if $N_G(P)/P$ is abelian, then $L|_P$ is indecomposable.

\end{lemma}

Proof.  Suppose first that $M$ has defect zero. It is well known
that  $M |_P$  is a projective module, and hence a multiple of
$\rho_P^{{\rm reg}}$.  Whence $(1)$.

Next, suppose that $M$ has defect $d$. Set $N=N_G(P)$. By \cite[Lemma VII.1.5]{Fe}, $M|_N=L \oplus
A_1\oplus A_2 $, where $A_1$ is projective, $A_2|_P$ is projective
and $L $ is the Green correspondent of $M$. Therefore, $(A_1\oplus
A_2)|_P$ is projective. As $P=\langle y\rangle $ is cyclic, every
projective $FP$-module is free, so $(A_1\oplus A_2)|_P=\frac{\dim
M -\dim L}{|P|}\cdot \rho_P^{reg}$. Recall that $L$ is
indecomposable as an $FN$-module (\cite[Theorem III.5.6]{Fe}) and
uniserial (\cite[Theorem VII.2.4]{Fe}), that is, the submodule
lattice of $L$  is a chain. Set $x=1-y$ in the group algebra $FN$,
$L_0=L$ and $L_i=x^iL$ for $i=1,\ldots ,e$,  assuming $L_{e}=0$
and $L_{e-1}\neq 0 $. Observe that $L_1$ is an $FN$-module
(indeed, for $n\in N $ we have $nL_1=(1-nyn^{-1})L=(1-y^j)L$ for
some integer $j>0$, and $1-y^j=(1-y)+(1-y)y+\cdots
+(1-y)y^{j-1}$). \itf $L_i$ is an  $FN$-module  for every $i$. As
$P$ acts trivially on every quotient $L_i/L_{i+1}$, the latter module is
completely reducible, and hence irreducible, since $L$ is
uniserial, for every $i$.
By \cite[Theorem VII.2.4]{Fe},
all the composition factors of  $L$ are of the same dimension $c$,
say. By the definition of the $L_i$'s, it follows that  $L |_P$ is a direct sum of $c$ copies of an
indecomposable \rep of $P$ of dimension $e$, and   $e= \dim L /c$.  Note that  $e<|P|$, as
otherwise $M$ would be of defect 0.  
In addition, if $N_G(P)/P$ is abelian, then $c=1$. So the $(2)$  follows.

\medskip
{\bf Remark}.  Observe that, if $G$ is a quasi-simple
group of Lie type with a non-trivial cyclic  $\ell$-subgroup $P$ of order $\ell^d$,
then $P$ is a TI-subgroup (see \cite{BL}, or \cite[Lemma 3.3(ii)]{Z99}).
This implies that every \ir $FG$-module is either of defect 0 or defect $d$,
as the defect group is the intersection of two Sylow $p$-subgroups.

\begin{corol}\label{ddw}  Let $G,P,M$ be as in Lemma $\ref{ww1}$, with $M$ of defect $d$,
and let $1\neq g\in P$. Suppose that  
$g$ is almost cyclic on $M$. Then
either $M=L$ and $\dim L<|P|$, or $ M|_P=\rho_P^{reg}\oplus L$
and $P$ is trivial on $L$. In the latter case, $\dim L =c$, where $c$ is
the dimension of an \irr of $N_G(P)/P$; in particular, if
$N_G(P)/P$ is abelian, then $\dim M=|P|+1$.
\end{corol}

Proof. Obviously, $g$ is almost cyclic on $L$. By Lemma \ref{ww1}, this implies that
either $L|_P$ is indecomposable or $L|_P$ is trivial. In the former case,  $\dim L<|P|$.
Suppose $M\neq L$. Then
$M=\rho_P^{reg}\oplus L$, and $L|_P$ is trivial. Observe
(cfr. the proof of Lemma \ref{ww1}) that $L$ is indecomposable as an
$FN_G(P)$-module; hence, as $L|_P$ is trivial, it is in fact an irreducible $F(N_G(P)/P)$-module. It follows
that $\dim L=c$. 
In particular, if
$N_G(P)/P$ is abelian, then $\dim M=|P|+1$.

\section{Groups with a 
normal subgroup of symplectic type}\label{sec4} 

In this Section, we collect miscellaneous results concerning groups containing normal subgroups of 'symplectic type', with special focus on the occurrence of almost cyclic elements in representations of such groups, which will be essential in the sequel of the paper. In particular, some applications to primitive linear groups containing non-central solvable normal subgroups will be obtained (cfr. Lemma \ref{p0} and Theorem \ref{a66}). 

\medskip

Let $E$ be a finite $r$-group for a prime $r$. We recall that $E$ is
said to be `of symplectic type' if it has no non-cyclic
characteristic abelian subgroups. The structure of such groups is
well understood (e.g. cf. \cite[p.109]{Asch}). [Namely, by an old result
of Philip Hall, if $E$ is a $r$-group of symplectic type, then $E$
is the central product of subgroups $A$ and $R$, where: 1) either
$A$ is extraspecial or $A=1$, and 2) either $R$ is cyclic or $R$ is
dihedral, semidihedral or quaternion, of order $\geq 2^{4}$.]
Certain $r$-groups of symplectic type naturally appear in the
Clifford theory of linear groups as irreducible subgroups of
primitive linear groups $G$ over algebraically closed fields, when $G$
has a non-central solvable normal subgroups. Namely, these
$r$-groups either are extraspecial of order $r^{1+2n}$ (of prime
exponent $r$ if $r$ is odd, and of exponent $4$ if $r=2$) or they
are $2$-groups of order $r^{2+2n}$ and exponent $4$ (with cyclic
centre of order $4$ and derived subgroup of order $2$). Their
structure is fully described, e.g. in \cite[p. 149, Table
4.6.A]{KL}. Their faithful irreducible representations over an
algebraically closed field of characteristic $\ell\neq r$ are
also well-known (e.g. cf. \cite[p.335]{Su} and \cite[pp. 149-150]{KL}).
In particular, they are all of degree $r^{n}$. Moreover, they are
uniquely determined by their restrictions to $Z(E)$, and their
characters vanish outside $Z(E)$.

\medskip
In the sequel of the paper, by `$r$-group of symplectic type'
 we always mean a group of
the above kind (i.e. one of the groups listed in \cite[Table
4.6.A]{KL}).

\smallskip
{\bf Remark}. In \cite{LS} Landazuri and Seitz considered a class of
$r$-groups, called 'groups of extraspecial type', which appear as
unipotent radicals of certain parabolic subgroups of finite groups
of Lie type. These groups are closely related to our 'groups of
symplectic type'. Namely, if $G$ is a group of extraspecial type,
for any subgroup $Z_1$ of index $r$ in $Z(G)$ the quotient group
$G/Z_1$ is a group of symplectic type with centre of order $r$.

\begin{lemma}\label{hb7} 
 Let $G$ be a finite group containing a normal subgroup $E$, where $E$ is an
  $r$-subgroup  of symplectic type and $|E/Z(E)|=r^{2n}$. Suppose that $G= E \cdot S$,
where $S=\langle g  \rangle$ 
and  $Z(E)\subseteq Z(G)$. Let $M $ be an \ir $FG$-module non-trivial on $Z(E)$.
Then $M$ is irreducible as $FE$-module and $\dim _FM=r^n$.
\end{lemma}

Proof. See  \cite[Ch.IX, Lemma 2.5]{HB}, where $E$ is supposed to be extraspecial. However, the proof remains valid without changes for $E$ of symplectic type.

\begin{lemma}\label{ex7} 
 Let $G$ be a finite group containing a normal subgroup $E$,
where $E$ is an $r$-subgroup  of symplectic type and $|E/Z(E)|=r^{2n}$. Suppose that $G= E \cdot S$,
where $S=\langle g  \rangle$, 
  $Z(E)\subseteq Z(G)$,
$C_S(E)=1$.  Suppose that   $|g|=r^n-\ep$,
where $\ep\in\{1,-1\}$. Suppose furthermore  that $E$ contains no $g$-invariant
non-abelian subgroups. Let $M$ be an \ir $FG$-module faithful on
$E$. Then the \f holds:

\medskip
$(1)$ If $\ep =-1$, then $M|_S$ is isomorphic to a submodule of codimension $1$ in $\rho_S^{\rm reg}$ and the matrix of $g$ on $M$ is cyclic.

\medskip
$(2)$ If $\ep=1$, then:

\medskip
$(i)$   $M|_S\cong \rho_S^{{\rm reg}}\oplus L$,
where $L$ is a $1$-dimensional $FS$-module.

\medskip
 $(ii)$ Suppose that $\ell={\rm char}~F>0$, and $P\neq 1$  is
  the Sylow $\ell$-subgroup of $S$. Let $S=PB$, where $B=\lan
b\ran$ is an $\ell'$-subgroup of $S$.
Let $U$ be a sum of some eigenspaces of $b$ on $M$.  Then
the matrix of $g$ on $U$ is cyclic if and only if
$\dim U\equiv 0\,({\rm mod} ~\ell)$.

\medskip
 $(iii)$ Let $1\neq z\in P$. Then the matrix of $g$ is cyclic on $(1-z)M$.

\end{lemma}

Proof. Set $V=E/Z(E)$. Then $V$ is a non-degenerate symplectic space over $F_r$ with respect to the bilinear form on $V$
 induced by the commutator
map $(a,b)\ra [a,b] $ ($a,b\in E$). Let $h$ be the automorphism of $V$ induced by the
conjugation action of $g$. Then $h$ can be viewed
as an element of the symplectic group $Sp(2n,r)$. Note that
$|h|=|g|=r^n-\ep$, as $C_S(E)=1$.

Let $t\neq 1$ be a power of $h$. Then $t$
acts fixed-point freely on $V\setminus \{0\}$. Indeed, let $V^t$ be the fixed point
subspace of $t$ on $V; $ it is well known that $V^t$ is non-degenerate, as $t $
is semisimple. As $hV^t=V^t$, it follows that $h$ is orthogonally decomposable
 on $V$.  However, this is equivalent to saying that $E$ has $g$-invariant non-abelian subgroups, against our assumption.

Suppose first that $\ell=0$ or $(\ell,|S|)=1$. Then we can apply
 \cite[Theorem 9.18]{DH} (the case  $r= 2$ being refined in \cite[Lemma 4.4]{HZ}). Thus  $\rho^{{\rm reg}}_S=M|_S+W$ if $\ep=-1$ and $M|_S=\rho^{{\rm reg}}_S+W$ if $\ep=1$, where  $W$ is a 1-dimensional $FS$-module. So in this case the lemma follows.

Next, suppose $(\ell ,|S|)\neq 1$. We first show that the $b$-eigenspaces on $M$
are all of dimension $|P|$, except one of dimension $|P|+\ep$.  Recall that
 $M$ lifts to characteristic zero
(this is true for every \ir \rep of a finite solvable group, e.g. see \cite[p. 135]{serr}). As
$(|B|,\ell )=1$, the dimensions of the $b$-eigenspaces on $M$ are the
same as in the zero characteristic case. In the latter case
the claim follows from (1) and (2)(i),  already proven  for characteristic zero.

 Let $M=M_1\oplus \cdots \oplus M_k \oplus M_0$, where $M_1\ld M_k$
are the $b$-eigenspaces of dimension $|P|$ and $ M_0$ is the $b$-eigenspace of dimension $|P|+\ep$. Obviously, each of the $M_i$'s $(0\leq i\leq k)$ is $P$-stable.

Let $L$ and $ N=N_G(P)$ be as in Lemma \ref{ww1}.
Observe that $N=N_E(P)S$. Moreover, as $[N_E(P),P]=E\cap P=1$,
we have $N_E(P)=C_E(P)$. We claim that $C_E(P)=Z(E)$. Indeed, by the argument above, every non-identity  element of $P$ acts fixed-point freely on the non-identity elements  of $E/Z(E)$.
\itf $N=Z(E)S$ is abelian, and
therefore, by Lemma \ref{ww1}, $L|_P$ is indecomposable. (Note that $M$ is of non-zero defect as $\dim M$ is coprime to $\ell$.)
In particular, $\dim L<|P|$. Also notice that, since $\rho_P^{{\rm reg}}$ is indecomposable,
 the decomposition of $M|_P$ given in Lemma \ref{ww1} consists of indecomposable summands. It follows,
by the Krull-Schmidt theorem, that
 $M_i|_P\cong \rho_P^{{\rm reg}}$ for $i=1\ld k$, by dimension reasons. In addition, if $\ep =-1$ then $M_0\cong L$, whereas if  $\ep =1$ then $M_0|_P=L\oplus \rho_P^{{\rm reg}}$. We conclude that $M|_S$ is isomorphic to a submodule of of codimension 1 in $  \rho_S^{{\rm reg}}$, and $\dim L=|P|-1$ if $\ep =-1$,
whereas if $\ep=1$ then $M|_S=\rho_S^{{\rm reg}} \oplus L$
 and $\dim L=1$. So we get (1) and item (i) in (2).

\med
Let $U$ be as in (2)(ii). It follows from (2)(i) that the matrix of $g$ on $M$ is cyclic if and only if $U$ does not contain $L$. The latter is equivalent to  assertion (ii).

\med
(iii) Obviously, the matrix of $g$ is cyclic on every quotient module $M/X$ provided $X$ contains $L$. Let $X$ be the kernel of the homomorphism $M\ra (1-z)M$. Then $L\subseteq X$ and $M/X\cong (1-z)M$, as desired. (Note that $L\subset X$ because $z$ is an $\ell$-element, and therefore it acts as the identity on $L$.)

\begin{corol}\label{cr2}  Let $g,M$ be as in items $(1)$ or $(2)$ of
Lemma $\ref{ex7}$. Then the matrix of $g$ on $M$ is almost cyclic.
\end{corol}

\begin{corol}\label{cr3}  Let $g,M$ be as in items $(1)$ or $(2)$ of Lemma $\ref{ex7}$,
and let $h\in \lan g\ran$ be such that $1\neq |h|<|g|$.   Then the
matrix of $h$ on $M$ is not almost cyclic, except for the case where
$r^n=3$, $|g|=4$ and $|h|=2$.
\end{corol}

Proof. Set $T=\lan h\ran$ and let $d=|S:T|=|g|/|h|$. In case (2) of
Lemma \ref{ex7} we have that $M|_T=d\cdot \rho_T^{{\rm reg}}\oplus
L|_T$, so the claim is obvious. In case (1) $M|_T$ is of codimension
1 in $d\cdot \rho_T^{{\rm reg}}$ so we are done unless  $d=2$ and
$|T|=2$. However, if $2=|T|=|g|/2= (r^n+1)/2$
then $3= r^n$, whence $r=3$ and $n=1$.  
If $r=3$ and $n=1$, then $|g|=4$ and $|h|=2$. In this case $h$ is
obviously almost cyclic.

\begin{corol}\label{cr4} Let  $G =E\lan g\ran\subset GL(r^n,F)$, where
 $E$ is a normal subgroup  of symplectic type, $|E/Z(E)|=r^{2n}$, $C_{\lan g\ran}(E)=1$,
$C_E(g)=Z(E)$
  and $g$   is of order coprime to $r$.
 Let $\overline{g}$ be the projection of $g$ into $Sp(2n,r)\subset Aut\, E$.
  Suppose that $\overline{g}$ is orthogonally indecomposable
and  $g$ is almost cyclic. Then $\overline{g}$ is of order $r^n+1$
or $r^n-1$.
 \end{corol}

Proof. Set $V=E/Z(E)$. Then $V$ is a non-degenerate symplectic space
and $\overline{g}$ is completely reducible as an element  of
$Sp(2n,r)$. Moreover, it is well known that either $\overline{g}$ is
\ir or it preserves a totally isotropic subspace of $V$. In fact, in
the second case the assumptions that $(|g|,r)=1$ and $\overline{g}$
is orthogonally indecomposable, imply that $\overline{g}$
preserves a maximal totally isotropic subspace of $V$. Now, suppose
we are in the former case. Then $|\overline{g}|$ divides $r^n+1$.
Let $g_1$ be an element of order $r^n+1$ in $Sp(2n,r)$ such that
$\overline{g}\in \lan g_1 \ran$. Then, by Corollaries \ref{cr2} and
\ref{cr3},
$|g_1|=|\overline{g} |$, and we are done.  

In the latter case $\overline{g}$ has order dividing $r^{n}-1$,  and
the result again follows with same argument, as in the former case,
from Corollaries \ref{cr2} and \ref{cr3}.

  \begin{corol}\label{ng3}  Let $g,M$ be as in items $(1)$ or $(2)$ of Lemma $\ref{ex7}.$
In addition, suppose that $|g|=p^a$   for some integer $a>0$ and
some prime $p$.
   Then one of the following holds:

\medskip
$(1)$ $r=2$, and either $|g|=p$ is a Fermat or Mersenne prime,   or
$|g|=9$;

\medskip
$(2)$ $r$ is odd, and either $n=1$, $|g|=2^a$ for some integer a and
$r$ is a Fermat or Mersenne prime, or $r^n=9$ and $|g|=8$.
  \end{corol}

Proof. As $|g|=p^a=r^n+1$ or $r^n-1$,  it follows from Lemma
\ref{zgm}  that  either $p$ or $r$ equals 2. Moreover,
if $p^a=r^n+1$ then either $p=2,n=1$ or $r=2, a=1$ or $p^a=9$.  If
$p^a=r^n-1$, then $p^a+1=r^n$. So again either $r=2,a=1$ or $p=2,
n=1$ or $r^n=9$. Thus, if $r$ is odd, then either $r^n=9$ or $p=2$
and $r$ is a Fermat or Mersenne prime; if $r=2$ then either $|g|=9$
or $|g|$ is a Fermat or Mersenne prime.

\begin{lemma}\label{a55}  Let $G=E\lan g\ran$, where
 $E$ is a normal  subgroup of $G$ of symplectic type,  $|E/Z(E)|=r^{2n}$ and
$|g|$ is a prime-power
 coprime to $r$. Let $\overline{g}$ be the projection of $g$ into $Sp(2n,r)$.
 Let $\phi\in \Irr_F G$ be faithful with $r\neq  \ell$. Suppose that
   $\phi (g)$ is almost cyclic. Then $\overline{g}$ is orthogonally indecomposable
in $Sp(2n,r)$
    and $|\overline{g}|=r^n+1$ or
$r^n-1$.
 Moreover, $|Spec\,\phi( g)|=r^n$ in the former case and $r^n-1$
in the latter case.
\end{lemma}

Proof. Set $V=E/Z(E)$. Then we may write $V=V_1\oplus \cdots \oplus
V_k$, where the $V_i$'s, for $i=1\ld k$, are non-degenerate,
mutually orthogonal, orthogonally indecomposable subspaces of $V$
invariant under the action of $\overline{g}$. Thus
$\overline{g}=\diag (h_1\ld h_k)$, where $h_i=\overline{g}|_{V_i}$
for $i=1\ld k. $ Let $H_i$ denote the group $Sp(V_i)$, for $i=1\ld
k$, so that $h_i\in H_i$. Set $H=H_1\times \cdots\times H_k$. By Lemma \ref{hb7},
 $\phi |_E$ is irreducible, and hence $\phi$ has degree $r^n$. It is also well
known (e.g. see \cite{Ge}) that $\phi |_E$ extends to a \rep $\tau$,
say, of the semidirect product $EH$ such that the restriction $\tau
|_H$ is the tensor product of the generic Weil \reps $\tau_i$ of the
groups $H_i$, having degree $r^{n_i}$, where $n_i=\dim V_i/2$, and
moreover $\tau (\overline{g})$ differs from $\phi (g)$ by a scalar
multiple. In particular, $\tau (\overline{g})$ is also almost
cyclic. As $\tau (\overline{g})=\tau_1(h_1)\otimes
\cdots\otimes\tau_k(h_k) $, it follows that $\tau_i(h_i)$ is
almost cyclic for every $i$. Now, suppose that
$h_i=\overline{g}|_{V_i}\neq \Id_{V_i}$. As $h_i$ satisfies the
assumptions of Corollary \ref{cr4}, it follows that $h_i$ has order
$r^{n_i}+1$ or $r^{n_i}-1$. As  $|g|$ is a prime-power, using
properties of Zsigmondy primes we readily deduce that $|h_i|=|h_j|$,
unless $|h_i|=1$ or $|h_j|=1$ ($1\leq i,j\leq n$). As $\tau
(\overline{g})$ is almost cyclic, it follows that  $|h_i|\neq 1$
only for one of the $i$'s. We can assume that this $i$ is 1. Assume
that $k>1$. Then $\tau (\overline{g})=\tau_1(h_1)\otimes \Id_m $,
where $m>1$. But then $\tau (\overline{g})$ is not almost cyclic,
and hence also $\phi (g)$ is not almost cyclic, against our
assumptions. Thus $k=1$, and Corollary \ref{cr4} applies. The
additional claim on ${\rm Spec}\, \phi (g)$ follows from Lemma
\ref{ex7}.

\begin{lemma}\label{mm5}  Let $H=\lan g\ran$ be a cyclic
$r$-group ($r$  a prime) and let $\phi :H\rightarrow GL(n,F)$
be a complex representation of $H$ with character $\chi$. Suppose
that $\chi (g^i)=0$ for $(i,|g|)=1$ and $\lambda$ be an eigenvalue
of $g^r$ of multiplicity $d$. Then all the $\mu$'s in $F$ such that
$\mu^r=\lambda$ are eigenvalues of $g$  of
 multiplicity $d/r$.
\end{lemma}

Proof. See \cite[Lemma 2.4]{DMZ10}.

\medskip

We recall here that an irreducible subgroup of $GL(V)$, where $V$ is a vector space over a field $F$, is said to be primitive on $V$ if it does not preserve any direct sum decomposition of $V$ into non-trivial subspaces of equal dimension.

\medskip

The following lemma essentially follows from Clifford theory.

\begin{lemma}\label{p0}  Let $G$ be a finite primitive subgroup of $GL(V)$, where $V$
is a finite-dimensional vector space over 
$F$.
 Let $S(G)$ denote the maximal solvable normal
subgroup of $G$. Suppose that $S(G)\neq Z(G)$. Then the following
holds:

\medskip
$(1)$ $G$ contains a normal $r$-subgroup $E$ of symplectic type for
some prime $r$ (so, the group $E$ has exponent $r$ if $r$ is odd,
whereas it has exponent $4$ if $r=2)$. Furthermore, $r\neq \ell$.

\medskip
  $(2)$ If $G$ is tensor-indecomposable, then $E$ is irreducible and $\dim V=r^n$, where $|E/Z(E)|=r^{2n}$.
\end{lemma}

Proof. (1) Let $E$ be a minimal non-central solvable normal subgroup
of $G$. Then, by Clifford's theorem, $E$ is non-abelian.
Furthermore, $Z(E)$ consists of
scalar matrices and the commutator subgroup $E^{\prime }$ is contained in $%
Z(E)$. So $E/Z(E)$ is abelian. As $E$ is nilpotent, again the
minimality assumption implies that $E$ is a $r$-group for some prime
$r$, and hence is an $r$-group with no non-cyclic characteristic
subgroups. Moreover $O_{\ell }(G)=1$, as $G$ is irreducible,
 whence $r\neq \ell$. Next, suppose that $r$ is odd. Then one easily
sees that $E$ contains a non-central element of order $r$. Let
$\Omega _{1}(E)$ denote the subgroup of $E$ generated by all its
elements of order $r$. Then $\Omega _{1}(E)=E$. As $E/Z(E)$ is
abelian, any two elements of $E$ commute mod $Z(E)
$. It follows that $E/Z(E)$ has exponent $r$, which in turn implies that $%
\left\vert E^{\prime }\right\vert =r$. Indeed, let $x,y\in E$. As
$y^{r}\in Z(E)$, $1=[y^{r},x]=[x,y]^{r}$. As $E^{\prime }$ is
cyclic, $\left\vert
E^{\prime }\right\vert =r$. Now, for any $x,y\in E$, $%
(xy)^{r}=x^{r}y^{r}[x,y]^{r(r-1/2)}=x^{r}y^{r}$. Thus, if $x$ and $y$
have order $r$, $xy$ also has order $r$. As $\Omega _{1}(E)=E$, we
deduce that $E$ has exponent $r$. Finally, suppose that $r=2$. If
$E$ does not contain non-central involutions, then $E$ is the
quaternion group of order $8$. Otherwise, arguing as above one sees
that $E/Z(E)$ has exponent $2$ and $E$ has exponent $4$.

\medskip
(2) It follows from Clifford theory (e.g. see \cite[pp. 139 -
141]{Za93}) that if $G$ is a primitive subgroup of $GL(n,R)$, where
$R$ is an algebraically closed field, then $G$ can be viewed as a
subgroup of the tensor product $G_{1}\otimes \cdot \cdot \cdot
\otimes G_{m}$,
 where $G_{i}$ for $1\leq i\leq m$ is a primitive tensor-indecomposable subgroup of
$GL(n_{i},R) $, $n=n_{1}\cdot \cdot \cdot n_{m}$, and every normal
subgroup of $G_{i}$ is either irreducible or scalar. As $G$ in (2)
is assumed to be tensor-indecomposable, we have that $m=1$, and the
result follows from (1).

\begin{theo}\label{a66}  Let $G$ be a
 primitive subgroup of
$GL(m,F)$ with  non-central maximal solvable normal subgroup
$S(G)$. Suppose that $G=\lan g^G\ran$, where $g$ is almost cyclic
and $g^p  \in Z(G)$ for some prime $p>2$. Then $G$ contains an
irreducible normal $r$-subgroup E of symplectic type, and one of the
following holds:

\medskip

$(1)$ $m=p=r$ and $ G=Z(G)\cdot E\cdot Sp(2,r)$.
\medskip

$(2)$ $m=2^n$ for some natural number $n$, $|E/Z(E)|=2^{2n+1}$ and
$\overline{G}:=G/(Z(G)E)$ is
isomorphic to a subgroup of  $Sp(2n,2)$ generated by a
conjugacy class of elements
 $\overline{g}$ of order $p=2^n-1$ or $2^n+1$.
\end{theo}

Proof.   As, by assumption, the $p'$-part of $g$ is scalar, we may
assume that $|g|$ is a $p$-power
 without loss of generality. By Lemma \ref{p0}, $G$ contains a normal
$r$-subgroup $E$ of symplectic type for some prime $r$, where $r\neq
\ell$ and $Z(E)\subseteq Z(G)$. Let $|E/Z(E)|=r^{2n}$. Set $K:=\lan
E,g\ran$. As $G=\lan g^G\ran$ and $C_G(E)$ is normal in $G$, we have
$[E,g]\neq 1$. Let $V$ be the underlying space of $GL(m,F)$. We
shall show that $E$ acts on $V$ irreducibly.

Assume first that $(|g|,\ell )=1$, where $\ell ={\rm char}F$. Then
$V$ is completely reducible as an $FK$-module. Let $V=V_1\oplus
\cdots \oplus V_t$, where the $V_i$'s are irreducible
$FK$-submodules. Therefore, every $V_i$ is a faithful irreducible
$FE$-module (see  Lemma \ref{hb7}). So $\dim V_i=r^n$ for every
$i=1\ld t$. For each $i$ let $g_i$ be the projection of $g$ to
$V_i$. Then $g_i$ is almost cyclic. Let $\mu$ be an eigenvalue of
$g$. Then the $\mu$-eigenspace of $g$ is the sum of the
$\mu$-eigenspaces of some of the $g_i$'s.

We have two cases: (a) $(|g|,r)=1$; (b) $ |g|$ is an $r$-power.

In case (a),  let $g^p=\lambda \cdot Id$,  where
$\lambda\in F$.
 By Lemma \ref{a55}, $g_i^d$ is scalar in $GL(V_i)$
where $d=r^n\pm 1=p$, so $g_i^d=\lambda \cdot \Id$. Moreover, all the
$p$-roots of $\lambda$, except one of them when $d=p=r^n+1$, occur as
eigenvalues of $g_i$. Therefore, at least two eigenvalues $\nu, \mu$
of $g$ are common on $V_1$ and $V_2\oplus \cdots \oplus V_t$. This
contradicts the assumption that $g$ is almost cyclic, unless $t=1$,
that is $V$ is irreducible as an $FE$-module, or $p=3$, $t=2$. In
the latter case, we have $r=2$, $n\leq 2$. 

In case (b), where $r=p$, by \cite[Lemma 2.5]{DMZ10}, either
$|E|=p^3$ or $p=3$ and $|E|=3^5$. In the latter case $n=2$ and,
again by  \cite[Lemma 2.5]{DMZ10}, $g $ has $r^n=3^2$ distinct
eigenvalues which is false as $g^3$ is scalar. In the former case,
by \cite[Lemma 2.5]{DMZ10}, $g$ is cyclic and hence $m\leq p=r$
unless, possibly,  when $p=3$, which implies  $t=2$.
So either $t=1$, or $t=2$ and $p=3$.

Next, suppose that $g$ is an $\ell$-element, that is $p=\ell \neq
r$, and hence $g^p=1$. Let $V_1\subset V_2\subset \cdots \subset V_t=V$ be a
composition series for $K=\lan E,g\ran$. Then every composition
 factor is a faithful \ir  $FE$-module. Denote by $g|_{V_2/V_1}$ the element
of $GL(V_2/V_1)$ induced by $g$ on $V_2/V_1$. Then, by Lemma
\ref{a55}, $p=|g|=r^n\pm 1$, and so by Lemma \ref{ex7} the Jordan
forms of $g|_{V_1}$ and $g|_{V_2/V_1}$ contain a block of size at
least $|g|-1$.  As $g$ is almost cyclic,  again we must have either
$t=1$, or $p=3$ and $V=V_2$. In the latter case we get
$r=2$ and $n\leq 2$.

Now, let $t=2$, $p=3$. Observe that $G$ is tensor-decomposable for
$t>1$ (see Lemma \ref{p0}(2)). As $t=2$, we have $zg=g_1\otimes
g_2$ for some scalar matrix $z$ where $g_1\in GL(m/2,F)$ and
$g_2\in GL(2,F)$. By Lemma \ref{p74}, both $g_1$ and $g_2$ are
cyclic. Recall that $\Id \otimes g_2$ centralizes $E$ and
$g_1\otimes \Id$ normalizes $E$ and produces the same \au on $E$
as $g$. Therefore, $g_1^3$ and $g_2^3$ are scalar. Therefore $\dim
V_1\leq 3$. Suppose first that $\ell\neq 3$. Assume $\dim V_1= 3.
$ Then both  $ g_1$ and $ g^2_1$ have trace zero. It follows that
the traces of $zg=g_1\otimes g_2$ and $z^2g=g^2_1\otimes g^2_2$
are 0. Hence  the traces of $ g$ and $ g^2$ are also zero. By
Lemma \ref{mm5}, $g $ is not almost cyclic. Therefore, $\dim V_1=2
$ and hence $r=2$. So $G/Z(G)E\subseteq SL(2,2) $, and hence
$G/Z(G)E$ is of order 3 (as $G/Z(G)E $ is generated by the
conjugates of $\bar g$). We conclude that  $G=K$, and the claim
that $E$ is irreducible follows, again by Lemma \ref{hb7}.

Next, let $\ell=3$.  Then  $|g_1|=|g_2|=3$. If $\dim V_1=2$, then $r=2$ and we have  $G=K$
 and $t=1$, as above. Let $\dim V_1=3$. As $g_1$ is cyclic, the Jordan form of $g_1$ is a single block, and hence
the Jordan form of $g_1\otimes g_2$ consists of 2 blocks of size 3, which is false as $zg$ is almost cyclic.

Thus, in view of the above,   $V=V_1$, which means that $m=r^n$.
Suppose $r=p$. Then, as already seen above, $|E|=p^3$ and $n=1$,
which implies (1). Next, let $(r,p)=1$. Let $N$ be the normalizer
of $E$ in $GL(m,F)$. Then $N/EZ(N)\cong Sp(2n,r)$. Let
$\overline{g}$ be the projection of $g$ into $Sp(2n,r)$. By Lemma
\ref{a55}, $|\overline{g}|=r^n\pm 1$. So $r=2$ and we have (2).

\section{Some low-dimensional classical groups}

In this Section, we first consider semisimple elements of prime-power order of a group $G$  such that $SL(2,q)\subseteq G\subseteq GL(2,q)$, and determine the irreducible $F$-representations of $G$ in which such elements are represented by almost cyclic matrices. Next, we obtain results of the same kind for some other small dimensional linear groups (see Lemmas  \ref{ss3}, \ref{333}, \ref{s43}), which will be needed in Section 5 in order to deal with the general case when $SL(n,q)\subseteq G\subseteq GL(n,q)$, for any $n>2$. Finally, in Lemmas \ref{423}, \ref{442},  \ref{452} and  \ref{u62} we examine some low-dimensional symplectic and unitary groups which will also play a role in Section 5. We emphasize that in this Section, we do not restrict ourselves to Weil representations.

\smallskip

\begin{lemma}\label{sq2}  Let  $SL(2,q)\subseteq G\subseteq GL(2,q)$.
Then every \ir $F$-\rep of G lifts to characteristic zero. \end {lemma}

Proof. Denote by $Z$  the set of non-zero scalar matrices in $GL(2,q)$.
Let $\tau$ be an \ir $F$-representation of $G$. Obviously, $\tau$ extends to $Z\cdot G$, and $\tau$ lifts
\ii the extension lifts. Note that $Z\cdot G$ is of index at most 2 in $GL(2,q)$, so either
$Z\cdot G=GL(2,q)$ or $Z\cdot G=Z\cdot SL(2,q)$.

The shapes of the decomposition matrices modulo $\ell$ for $SL(2,q)$ show that the lemma is true for this case,
see for instance \cite[Ch. 9]{Bn}. (The reader should note that this does not hold for $PSL(2,q)$, e.g. see  \cite{Bu}.)
Therefore, it suffices to prove the lemma for   $G=GL(2,q)$, $q$ odd.

Assuming this, set $X=Z\cdot SL(2,q)$ and $\tau_1=\tau|_X$.
 Suppose first that $\tau_1$ is irreducible. Let $\phi_1$ be the lift of $\tau_1$. Looking at the character tables
of $SL(2,q)$ and $GL(2,q)$, one observes that $\tau_1$ extends to $G$.  Let $\phi$ be the extension, and  set $\tau_2=
\phi\pmod \ell$.
A priori, $\tau_2$ may not coincide with $\tau$.
 However, for every  $g\in G$ the conjugation by $\tau(g)$ and $\tau_2(g)$ yields the same automorphism of $\tau(X)$.
 By Schur's lemma, $\tau _2(g)=\tau(g)\lam(g) $ for some $\lam(g)\in F$.
 One readily checks that $g\ra\lam(g)$ is a group homomorphism. Therefore,
 $\tau_2=\tau\otimes \lam $. As $\lam$ is one-dimensional, $\lam$ lifts to
 characteristic zero. Let $\mu$ be the lift of $\lam$. Then $(\mu\up\otimes \phi)
 \pmod \ell=\tau$, as required.

 Next, suppose that $\tau_1$ is  reducible, and hence completely reducible by Clifford's theorem. Then $\tau_1$
 has two \ir constituents $\si_1,\si_2$, say, which are   $G$-conjugate, and hence are of equal dimension, which is at most $(q+1)/2$. As $G/X$ is cyclic, it follows from Clifford's theory that $\si_1,\si_2$ are not equivalent (see for instance \cite[Th. 19.13]{Hp2}). Let $\phi_1,\phi_2$ be lifts of $\si_1,\si_2$, respectively. Then $\phi_1,\phi_2$ are not equivalent, and have equal dimension at most $(q+1)/2$. Moreover, $\phi_1|_Z=\phi_2|_Z$; therefore, $\phi_1|_{SL(2,q)}=\phi_2|_{SL(2,q)}$ are not equivalent. It is well known that $SL(2,q)$ has exactly two non-equivalent complex \reps of equal degree (which is either $(q+1)/2$ or $(q-1)/2$), and they are $G$-conjugate.
It follows from  the character table of  $G=GL(2,q)$ that there exists an \irr $\psi$ of $G$ such that $\psi|_X=\phi_1\oplus \phi_2$. Set $\tau'_2:=\psi\pmod \ell$. Then $\tau'_2$ is irreducible (as neither $\si_1$ nor $\si_2$ is $G$-stable). Then
we claim that $\tau'_2$ is equivalent to $\tau$. Obviously, the Brauer character of $\tau'_2|_X$ coincides with that of $\tau|_X$.
Let $g\in G$, $g\notin  X$. Then $g$ permutes $\si_1,\si_2$, and hence
the matrix of $\tau'_2(g)$ has zero trace. More precisely,
both the Brauer character values of $\tau'_2(g)$ and $\tau(g)$ are 0 (see, for instance, \cite[Proposition 2.14]{DMZ0}).   \itf the Brauer characters of $\tau'_2$ and $\tau$  coincide, and hence
$\tau'_2$ and $\tau$  are equivalent.

\begin{lemma}\label{aq1}  Let  $SL(2,q)\subseteq G\subseteq GL(2,q)$, $q>3,$ $Z=Z(GL(2,q))$, and let
 $T$ be the subgroup consisting of the diagonal matrices in $G$.  Let $\tau$ an \ir $F$-\rep of $G$, such that $\tau|_{Z(G)}=\zeta\cdot \Id$,
where $\zeta\in \Irr Z(G)$.

\med
$(1)$  Suppose that $q$ is odd and $\dim\tau=(q-1)/2$.
Then  $ G\subseteq Z\cdot SL(2,q) $  and
$$\tau|_T=\zeta^T.$$

\med
$(2)$ Suppose that $q$ is odd and $\dim\tau=(q+1)/2$. Then  
$ G\subseteq Z\cdot SL(2,q) $  and
$$\tau|_T=\nu\oplus \zeta^T,$$ where $\nu$ is a  $1$-dimensional \rep of  $T$.

\medskip
 $(3)$ Suppose that  $\dim\tau=q-1$. Then
$$\tau|_T= \zeta^T,$$ unless $ G\subseteq Z\cdot SL(2,q) $,
in which case $\tau|_T=2\cdot \zeta^T$.

\med
$(4)$ Suppose that  $\dim\tau=q+1$. If $ G\subseteq Z\cdot SL(2,q) $ and $q$ is odd,
then $$\tau|_T=
\mu\oplus  2\cdot \zeta^T,$$ where $\mu$ is a $2$-dimensional \rep of  $T$.

If $q$ is even or $ G\not \subseteq Z\cdot SL(2,q) $, then  $$\tau|_T=
\mu\oplus   \zeta^T,$$ where $\mu$ is a  $2$-dimensional \rep of
$T$.

\med
$(5)$ Suppose that  $\dim\tau=q$. Then  $\tau|_T=\nu\oplus c\cdot \zeta^T$, where
$\nu$ is a  $1$-dimensional \rep of  $T$, and $c=1$ if $q$ is even or $ G\not \subseteq Z\cdot SL(2,q) $, otherwise $c=2$.
\end{lemma}

Proof.
By Lemma \ref{sq2}, $\tau$ lifts to characteristic zero. Let $\chi$ be the character of the lift.

 Let $U$ be the abelian subgroup of order $q$ consisting of the upper unitriangular
matrices in $G$. Then $T$ normalizes $U$ and $C_T(u)=Z(G)$ for every $1\neq u\in U$.
 Set $K=\Irr U$. Acting on $U$ by conjugation,  $T$ has
a single orbit on $U\setminus \{1\}$ if $q$ is even or $G\not\subseteq Z\cdot SL(2,q)$, and two orbits of size $(q-1)/2$ if  $q$ is odd and $G\subseteq Z\cdot SL(2,q)$.  Then this is also true for the (dual) action
of $T$ on $K$.

Let $M$ be the module afforded by $\tau$. For $\al\in K$, set $M_\al=\{m\in M: \tau(u)m=\al(u)m$ for all $u\in U\}$.
Then $T$ permutes the (non-zero) $M_\al$'s. It follows that,  for every $T$-orbit   $O$ on the $M_\al$'s,
$T$ stabilizes the subspace $M_O:=\sum _{\al\in O}M_\al$.
Note that $Z(G)\subseteq T $. It is easy to observe that the restriction of  $\tau|_T$ to $M_O$
yields  a \rep of $T$ equivalent to $\zeta^T$, where $\zeta\in\Irr Z(G)$, $xm=\zeta(x)m$ for $x\in Z(G)$ and $m \in M_O$.
As $M$ is irreducible, it is clear that $\zeta $ is the same  for every $T$-orbit $O$.

Suppose first that $\tau(1)=(q\pm 1)/2$. Then $G\subseteq Z\cdot SL(2,q)$.
By the above, applying  Clifford's theorem  to $TU$, it follows that,
if $\dim\tau=(q-1)/2$, then $\tau|_U$ is the sum of the characters of a $T$-orbit of length $(q-1)/2$,
whereas, if $\dim\tau=(q+1)/2$, then  $\tau|_U$ is the sum of $1_U$ (with \mult 2) and the  characters belonging to a $T$-orbit of  size $(q-1)/2$.

Next, suppose that  $\chi(1)\in \{q-1,q,q+1\}$. Then  $\tau|_U=\rho_U^{reg}+a\cdot 1_U$, where
$a=\chi(1)-q$. Therefore, for any $\tau$, the restriction $\tau|_U$ is the sum of
one-dimensional representations of $U$, each of \mult one, except when $ \chi(1)=q+1$,
in which case $1_U$ has \mult 2 and  the other \ir constituents have \mult 1.

This immediately implies all the statements of the lemma.

\begin{lemma}\label{cn2}  Let $SL(2,q)\subseteq G\subseteq GL(2,q)$, 
and let $1\neq g\in G $ be a semisimple
element of $p$-power order, where $p$ is an odd prime.
Let $M$ be an \ir $ FG $-module with $\dim M>1$ and let $\tau$ be the representation afforded by $M$. Then $\tau(g)$ is almost cyclic \ii $\dim M\leq |g|+1$.
Moreover,  in this case $(2,q+1)\cdot |g|$ equals  $q+1$ or $q-1$. \end{lemma}

Proof.
Firstly note that, by our assumptions, $q>3$. Let  $P$ be a \syl of $G$. As $Z\cdot SL(2,q)$  has index at most 2 in $G$, we have $P\subset Z\cdot SL(2,q)$. It follows that
it suffices to prove the result for $G=SL(2,q)$. Indeed, this is trivial if $G\subseteq Z\cdot SL(2,q)$. So, assume otherwise.
Then $Z\cdot G=GL(2,q)$, and hence we may assume $G=GL(2,q)$. The claim is
 obvious if  $\tau(SL(2,q))$
is irreducible. If not,  by Clifford's theorem, $\tau(SL(2,q))$ is a direct sum of two \ir constituents, permuted by any element
of $x\in G$, which is not in $Z\cdot SL(2,q)$. Note that an element $x \in GL(2,q)$ belongs to $Z\cdot SL(2,q)$ \ii
$\det x$ is a non-zero square in $F_q$, and hence
there is   $x\in C_{GL(2,q)}(P)$, which is not in $Z\cdot SL(2,q)$.
Therefore, $y$ permutes the \ir constituents of $\tau|_{ SL(2,q)}$ and commutes with $P$.
This implies that both of them have the same restriction to $P$, and hence no non-scalar element of $\tau|_P$ is almost cyclic.

Thus, we may assume that $G=SL(2,q)$.  By Lemma \ref{sq2}, $\tau$ lifts to characteristic zero. So, if $p\neq \ell$,
 it suffices to verify the lemma for $\ell=0$, which can be easily done examining the character table of $G$.
Therefore, from now on we assume $p=\ell$.

In this case $P$ is cyclic, and we assume that $g\in P$.
As $
p$ is odd and $\dim M\in \{(q\pm 1)/2, q\pm 1, q\}$, it follows that
$p$  divides $\dim M$ \ii so does $|P|$. It is also well known (e.g. see \cite {Bn}) that every $FG$-module is either of defect 0 or of defect $d$, where $|P|=p^d$.

If $g$ is diagonalizable (equivalently, $|g|$ divides $q-1$),
then the statement of the lemma about the almost cyclicity of $\tau(g)$  follows from Lemma \ref{aq1}, except, possibly, for the case where $q$ is even, $\dim M=q+1$, $|P|=|g|=q-1$
and $\mu(g)$ in Lemma \ref{aq1}(4) is scalar.
However, as $N_G(P)/P$ is cyclic, this contradicts the almost cyclicity of $\tau(g)$ by Corollary \ref{ddw}.

Therefore, we may assume that  $g$ is not diagonalizable (and hence $|g|$ divides $q+1$).

If $|P|$ divides $\dim M$,
then $M|_P$ is a projective $FP$-module, and hence $M|_P=m\cdot \rho_P^{reg}$
for some integer $m>0$. 
Then, obviously, the matrix of $\tau(g)$ is almost cyclic \ii
$m=1$ and $|g|=|P|=\dim M$,  
in which case   $\tau(g)$ is cyclic.

Thus, from now on we assume that $|g|$ is  coprime to  $\dim M$.

  By Lemma \ref{ww1} and Corollary \ref{ddw}, we have two options: either

  (i) $\dim M<|P|$
and $M|_P$ is indecomposable; or

(ii)
$M|_P=\rho_P^{reg}\oplus L|_P$, where $L\neq 0$ is an \ir $FN_G(P)$-module trivial on $P$.

Suppose that (i) holds, and let $g=h^{p^b}$, where  $P=\lan h\ran$ and $b\geq 0$.
Then the matrix of $\tau(h)$ is a Jordan block of size $t=\dim M$. So we are done if $b=0$, that is, $|g|=|P|$.
Suppose that $b>0$, that is,  $|g|<|h|$.
By \cite[Lemma 5.4]{DMZ8},
the  Jordan form of $g$ on $M$ contains at least two non-trivial blocks of equal size (and hence the matrix of $\tau(g)$ is not almost cyclic)
unless $t=p^{d-1}+1$ and 
$\tau(g)$  is a transvection. In the latter case  $\tau(G)$ is an \ir subgroup of $GL(M)$
generated by transvections.  The finite \ir subgroups of $GL(M)$ generated by transvections
are well known (see  \cite{W2},    \cite{ZS}).
Since $\ell=p\neq 2$,  these are isomorphic to
$SL(t, F_{\ell^m})$,  $SU(t, F_{\ell^m})$, $Sp(t, F_{\ell^m})$, or $SL(2,5)\subset SL(2,F)$ for $\ell=3$.  Clearly, none of these groups are  isomorphic to $\tau(G)$.
(Note that, as $(\ell,q)=1$ and $\ell\neq 2$, the isomorphisms  $\tau(SL(2,7))\cong SL(3,2)$  and  $\tau(SL(2,5))\cong SL(2,4)$ should be ignored.)

Now, suppose that (ii) holds (that is, $M|_P=\rho_P^{reg}\oplus L|_P$). Clearly, $|g|=|P|$, as the matrix of $\tau(g)$ is almost cyclic.
Furthermore, since $N_G(P)$ has an abelian normal subgroup of index 2, by Clifford's theorem $0<\dim L\leq 2$. Whence $|P|<\dim M\leq |P|+2$.

Suppose first that $q$ is even. As $|P| < \dim M$, we have
 $|P| < q+1$, and hence $|P|\leq (q+1)/3$. Therefore, as $\dim M\geq q-1$, $q-1\leq |P|+2\leq (q+7)/3$, and hence  $3q-3\leq q+7$, that is $q=4$. However, this forces  $|P|=5=q+1$, which is not the case.

So, suppose that $q$ is odd.
Then $|P|\leq (q+1)/2$, which implies $\dim M\leq \frac{q+1}{2} +2$.

If $|P|= (q+1)/2$, then $\dim M>|P|=(q+1)/2$ implies
$\dim M=q-1,q$ or $q+1$,  but the latter option is ruled out, as $|P|$ is coprime to $\dim M$. So $q-1\leq 2+\frac{q+1}{2}$, whence
$q\leq 7$.  However, $q\neq 7$, as $p>2$.  So $q=5$, whence $|g|=|P|=3$ and $\dim M=4,5$. Suppose that $\dim M=5$. Then $M$ is a $PSL(2,5)$-module, and in $PSL(2,5)$ the quotient $N(P)/P$ is abelian. This forces  $ \dim L = 1$, whence $ \dim M = 4$, as $M|_P=\rho_P^{reg}\oplus L$. A contradiction. So $\dim M = 4$, and we are done.

If $|P|<(q+1)/2$, then $|P| \leq (q+1)/4$,
and hence $(q-1)/2 \leq \dim M \leq 2+(q+1)/4$, whence $q\leq11$. The case $q=5$ is ruled out, as  $|g|<(q+1)/2$ implies $|g|<3$. As above, the case $q=7$ is also ruled out, as $p>2$.  Finally, in both the cases $q=9,11$, $M$ is a $PSL(2,q)$-module, and in $PSL(2,q)$  the quotient $N(P)/P$ is abelian, whence $ \dim L = 1$. A contradiction, as $M|_P=\rho_P^{reg}\oplus L$ would then imply $\dim M = 6$ for $q=9$ and $\dim M = 4$ for $q=11$, which is impossible.

As for the last claim in the statement, note that $|g|$ divides $q+\ep$, where $\ep=1$ or $-1$. Let $q$ be odd. If $|g|=(q+\ep)/2$, then the claim is true,
otherwise $|g|\leq (q+\ep)/4$. As $\dim M\geq (q-1)/2$, we have
$(q-1)/2\leq \dim M\leq 1+\frac{q+\ep}{4}$, whence $q\leq 6+\ep$.
But then $|g| \leq 2$, a contradiction. Now, suppose that $q$ is even. Then $q+\ep$ is odd, and hence either $|g|=q+\ep$, as required, or
$|g|\leq (q+\ep)/3$. As $\dim M\geq q-1$, we have $q-1\leq \dim M\leq 1+\frac{q+\ep}{3}$, whence $2q <6+ \ep$, a contradiction, as $q>3$ . 

\medskip

At this stage, we are left to deal with the case where  $SL(2,q)\subseteq G\subseteq GL(2,q)$, and $1\neq g\in G $ is a semisimple
element of $2$-power order.

We begin with an auxiliary Lemma:

\begin{lemma}\label{tt2}  Let $SL(2,q)\subseteq G\subseteq GL(2,q)$, where $q>3$ is odd, and let $g$ be a non-scalar
$2$-element of $G$. Let $\ell=2$, and
let $M$ be an \ir $ FG $-module of dimension $m>1$,   affording the representation $\tau$. Suppose that $\tau(g)$ is a
transvection. Then one of the \f holds:

$(1)$ $G=SL(2,5)$, $m=2$ and $\tau(G)=SL(2,4);$

$(2)$ $G=SL(2,7)$, $m=3$ and $\tau(G)=SL(3,2);$

$(3)$ $G=GL(2,5)$, $m=4$ and $\tau(G)=O^-(4,2).$
 
\end{lemma}

Proof. Let $G_1$ be the subgroup of $G$ generated by the $G$-conjugates of $g$. Clearly, $G_1$ is a (normal) subgroup of $G$ containing $SL(2,q)$. It follows that $\tau(G)$ is irreducible. Indeed, otherwise, by Clifford's theorem $M|_{G_1}=M_1\oplus M_2$,
where $M_1,M
_2$ are $G$-conjugate irreducible constituents (this is because $M|_{SL(2,q)}$ is either \ir or the sum of two \ir constituents). Hence, every element of $G_1$ non-trivial on $M_1$ is also non-trivial on $M_2$. However, 
a transvection stabilising $M_1$ and $M_2$ must be trivial either on $M_1$ or on $M_2$. This is a contradiction.

Thus, $M$ is an  \ir $FG_1$-module. Set $G_2=\tau(G_1)$. By Lemma \ref{wz1},
either $m$ is even and  $G_2\in\{ S_{m+1}$,   $S_{m+2}$,
$SL(m, q_1)$,  $Sp(m,q_1)$,  $O^\pm(m,q_1)$, $SU(m,q_1)\}$, or $m$ is odd and $G_2\in\{SL(m, q_1),SU(m, q_1)\}$, where $q_1$ is even in all the cases. \itf 
 one of the \f holds: 
(i) $G_1=SL(2,5)$ and $G_2= SL(2,4)$, $m=2$; (ii) $G_1=SL(2,7)$ and $G_2= SL(3,2)$, $m=3;$ 
(iii) $G_1=GL(2,5)$ and $G_2= O^-(4,2)$,   $m=4$. (Note that the group $Sp(4,2)$ is not isomorphic to $PGL(2,9)$ (e.g., see \cite[p. 4]{Atl}), so the case $G=GL(2,9)$ does not occur in our list.)

In the cases (i) and (ii) $|G:G_1|\leq 2$, so either $G=G_1$ or $G=GL(2,5)$ and $GL(2,7)$, respectively. The latter options are ruled out, as 
neither $GL(2,5)$ nor $GL(2,7)$ have $2$-modular \irrs of degree $2$ or $3$. In case (iii), we have $G=G_1$.
This completes the proof.  

\medskip

{\bf Remark}. In order to simplify the proof of some of the subsequent lemmas, it is worth observing explicitly at this point that, if one wishes to examine the representations of a group $G$, where $SL(2,q)\subseteq G\subseteq GL(2,q)$, it is enough to consider the
cases $G=GL(2,q)$ and $G=SL(2,q)$. Indeed, let $M$ be an \ir $FG$-module. Set $Z=Z(GL(2,q))$, and $G_1 = G \cdot Z$. Obviously, $M$ extends to an $FG_1$-module. Now, $G_1$ contains $SL(2,q) \cdot Z$. As the latter subgroup has index at most 2 in $GL(2,q)$,
it follows that, without loss of generality, we may assume either $G=SL(2,q)$ or  $G=GL(2,q)$.

Furthermore, note that, for $\ell=2$ it is sufficient to deal with the groups $PSL(2,q)$ and $PGL(2,q)$. This is obvious if $G=SL(2,q)$. If 
$G=GL(2,q)$, let $Z_2$ denote the Sylow 2-subgroup of $Z(G)$. Then $Z_2$ is in the kernel of $M$, so
$M$ can be viewed as an $F(G/Z_2)$-module. Set $\overline{G}=G/Z_2$ and observe that  
$\overline{G}= Z(\overline{G})\times K$, where $K \cong PGL(2,q)$. Whence the claim. See also the proof of Lemma 2.2
in \cite{TZ08}.

\begin{lemma}\label{q14}  Let $SL(2,q)\subseteq G\subseteq GL(2,q)$, where $q\equiv 1\pmod 4$, 
let $ g\in G\setminus Z(G) $ be a  $2$-element, and let $h $ be the projection of $g$ into $G/Z(G)$. Let $M$ be an \ir $ FG $-module of dimension $m>1$,  and let $\tau$ be the representation afforded by $M$. Then $\tau(g)$ is almost cyclic \ii the following holds:

$(1)$ $\ell \neq 2$, $g \in SL(2,q) \cdot Z(GL(2,q))$,  $\dim M=(q\pm 1)/2$ and $|h|= (q-1)/2$;

$(2)$ $\ell \neq 2$,  $g \notin SL(2,q) \cdot Z(GL(2,q))$,  and either $\dim M= q$ or $q -1$ and  $|h|=q-1$, or $q=5$, $\dim M = 4$, $|g|=8$ and $|h|=2$;

$(3)$ $\ell = 2$ and either $\dim M\leq |h|+1$ or $G=GL(2,5)$, $\tau(G)\cong O^-(4,2)$, $m=4$ and $\tau(g)$ is a transvection. 
 \end{lemma}

Proof. Note that, by the remark above, we may assume 
either $G=SL(2,q)$ or $G=GL(2,q)$.

Let us suppose that $\ell = 0$.  
Assume first that $|g|$ does not divide $q-1$. In this case
$G \neq  SL(2,q)$ (otherwise $|g|$ divides $q+1$, but $(q+1)/2$ is odd, so $|g|=2$, and hence $g$ would be scalar).
 So, let $G = GL(2,q)$. Then $g$ is irreducible and $g^2$ is scalar  (as $(q^2-1)/(q-1)=q+1$
and $(q+1)/2$ is odd). Thus $\tau(g)$ has exactly two distinct eigenvalues. As $\tau(g)$ is almost  cyclic, $\tau(g)$ is a pseudo-reflection. It then follows from {\rm (\cite[Lemma 3.1]{GS})} that $G$ can be generated by at most 4 conjugates of $g$, whence, by Lemma  \ref{nd5}, $\dim M \leq 4$. This implies $G = GL(2,5)$ and $|g| = 8$, yielding the exceptional case in (2).  

Thus, we may assume that $|g|$ divides $q-1$. In this case, w.l.o.g. we may assume that  $g\in T$, where $T$ is the subgroup of diagonal matrices in $G$. As  $\tau(G)$/$\tau(Z(G))$ is cyclic, items $(1)$ and $(2)$ of the statement follow from Lemma  \ref{aq1}. 
 
Since, by Lemma \ref{sq2},
every \ir $FG$-module lifts to characteristic zero, the above results hold for any $\ell \neq 2$. So, from now on, we assume that $\ell = 2$, and hence $q$ is odd.

Suppose that $\tau(g)$ is almost cyclic. The case where $\tau(g)$ is a transvection (recorded in item $(3)$ of the statement) follows from Lemma \ref{tt2}. So we can assume that $g^2\notin Z(G)$ (otherwise $\tau(g^2)=\ Id$, and hence $\tau (g)$ would be a transvection).
This implies that $g $ is reducible on the natural module of $G$, and hence $|g|$ divides $q-1$. Indeed, assume that $g$ is irreducible. Then $g$ is contained in a cyclic subgroup $X$ of 
$GL(2,q)$ of order $q^2-1$. As $X$ contains the subgroup of scalar matrices (of order $q-1$),
the order of $h$ divides $q+1$. By assumption, $(q+1)/2$ is odd, whence $g^2\in Z(G)$, a contradiction.

Thus, we can assume that $g\in T$. We wish to use Lemma
\ref{aq1}, which is stated in terms of $\zeta^T$, where $\zeta\in\Irr Z(G)$ is such that $\tau|_{Z(G)}=\zeta\cdot \Id$.

Let $a=|T:SZ(G)|$, where $S$ is the Sylow 2-subgroup of $T$. Set $\overline{G }=G/ (S \cap (Z(G))$ and let $\overline{S}$, $\overline{T}$ be the projections of $S,T$ into $\overline{G}$. Then clearly 
$|\overline{T}:\overline{SZ(G)}|=a$. Let us view $\zeta^T$ as a representation of $\overline{T}$ (clearly, $S\cap Z(G)$ is in the kernel of $\tau|_T$ as well as $\zeta^T$). That is, let us express $\zeta^T$ as $\zeta_1^{\overline{T}}$, where $\zeta_1$ is $\zeta$ viewed as a \rep of $Z(G)/(S\cap Z(G))$.  
Then $\zeta_1^{\overline{T}}|_{\overline{S}}=a\cdot \rho_{\overline{S}}^{reg}$. 

Note that $h\in \overline{S}\subseteq \overline{T}$. If $\tau (g)$ is almost cyclic, then so is $\zeta_1^{\overline{T}}(h)$. By Clifford's theorem, this implies $a=1$ and $|h|=|\overline{S}|$, and therefore $\zeta_1^{\overline{T}}(h)=\zeta^T(g)$ is represented by a Jordan block of size  $|h|$. 

Note that $a=1$ means that $|T:Z(G)|$ is a $2$-power. This implies that $G'=SL(2,q)$ has no $2$-modular \irr of degree $q+1$. Indeed, by \cite[9.2]{Bn}, $G'$ has no nilpotent block, and hence $M$ belongs to the principal 2-block. Then the claim follows from \cite{Bu}. In turn, this implies that case (4) of Lemma \ref{aq1} does not occur. Indeed, suppose that $\dim M=q+1$. Then $M|_{G'} $ is reducible. By Clifford's theorem, $M|_{G'} =M_1\oplus M_2$, where $M_1,M_2$ are irreducible $FG'$-modules of  dimension $(q+1)/2$. However,  $SL(2,q)$ has no \ir 2-modular \rep of such dimension (see \cite{Bu}).
This is a contradiction.

Furthermore, as  $\ell=2$, there are no irreducible $F$-representations of $G$ of degree $q$ (see   \cite{Bu}). Hence, case (5) of Lemma \ref{aq1} does not occur. 

Now, we apply Lemma \ref{aq1}. 
Suppose first that $G\subseteq Z\cdot SL(2,q)$. (Recall that $Z$ is the group of scalar matrices in $GL(2,q)$.) Then (see the argument above) $\dim M \neq q+1,(q+1)/2$. Moreover,   $\tau(g)$ is
 not almost cyclic in case (3),  whereas it is so in case 
(1).  So the lemma is true in this case. 
Next, suppose that $G$ is not contained in $Z\cdot SL(2,q)$. Then, cases (1) and (2) of Lemma \ref{aq1} are ruled out, whereas 
the matrix $\tau(g)$ is cyclic in case (3).

\begin{lemma}\label{jj2} 
 Let $J_n\in GL(n,2)$, $n>2$, be a Jordan block, where $2^{k-1}<n\leq 2^k$, $k>1$. Then the Jordan normal form of
$J_n^{2^{k-1}}$ is $\diag(J_2\ld J_2, \Id_{2^k-n})$.
\end{lemma}

Proof. We argue by induction on $k$. Clearly, the statement is trivially true for $k=2$, that is $n=3,4$.

By  \cite[Lemma 5.4]{DMZ8},
 the Jordan form of $J_n^2$ is $\diag(J_{m},J_{m})$ if $n=2m$ is even, and $\diag(J_{m+1},J_m)$ if $n=2m+1$ is odd.
To apply induction, we need the size
$s$ of each Jordan block of $J_n^2$ to satisfy the inequalities $2^{k-2}<s\leq 2^{k-1}$. If $n$ is even,  $2^{k-1}<n\leq 2^k$ implies  $2^{k-2}<n/2\leq 2^{k-1}$, as required. If $n$ is odd, $n<2^k$, so $\frac{n+1}{2}\leq 2^{k-1}$. Similarly,
$2^{k-2}<\frac{n-1}{2}$, except when $n-1=2^{k-1}$.

Suppose first that $n$ is even. Then, by induction, the Jordan form of  $J_{n/2}^{2^{k-2}}$
has $2^{k-1} -{n/2}$ trivial blocks. Hence
  the Jordan normal form of
$J_n^{2^{k-1}}$ has exactly $2^{k-1}+2^{k-1}-n = 2^{k}-n $ trivial blocks, as required.

Next, suppose that  $n$  is odd.

Suppose first that we are in the exceptional case where  $n=2^{k-1}+1$.
Then the Jordan form of $J_n^2$ is
$\diag(J_{m+1},J_m)=\diag(J_{2^{k-2}+1},J_{2^{k-2}})$. It follows that  $J_m^{2^{k-2}}=\Id$,
whereas $J_{m+1}^{2^{k-2}}$ is a transvection.
Therefore, the Jordan form of $J_n^{2^{k-1}}$ is $\diag(J_2,\Id_{n-2})$,
and hence $n-2=2^{k-1}-1=2^k-2^{k-1}-1=2^k-n$, as required.

In the general case, the Jordan form of $J_n^{2^{k-1}}=(\diag(J_{m+1},J_m))^{2^{k-2}}$ has $2^{k-1}-(m+1)+(2^k-m)=
2^k-n$ trivial blocks, as required.

\medskip

Note: A partial version of the result stated in the following Lemma is contained in a paper by Guralnick and Tiep (``Some bounds for $H^2$'', in preparation). For the reader's convenience, we have written down a comprehensive  proof.

\begin{lemma}\label{nn1}  Let $SL(2,q)\subseteq G\subseteq GL(2,q)$, where $q\equiv -1\pmod 4$. 
Let $ g\in G$ be a  $2$-element such that $g^2 \notin Z(G)$, and let $h$ be the projection of $g$ into $ \overline{G}:=G/Z(G)$. For $\ell=2$,
let $M$ be an \ir $ FG $-module of dimension $q-1$ (respectively, $(q-1)/2$),   and let $\tau$ be the representation afforded by $M$. 
Then $\tau(g)$ is almost cyclic if and only if  $q+1$ is a $2$-power (that is, $q$ is a Mersenne prime), and $|h|=q+1$ (respectively, $(q+1)/2)$.\end{lemma}

Proof. The "only if" part. Clearly, $G'=SL(2,q)$ and $\overline{G'}=PSL(2,q)$.  Let $t\in\lan g\ran$ be such that $t\notin Z(G)$, but $t^2\in Z(G)$, and let $\overline{t}$ be the projection of $t$ into  $\overline{G}$. Observe that $\overline{t}\in\overline{G'}$ (as $Z(GL(2,q))\cdot SL(2,q)$ has index 2 in $GL(2,q)$, $g^2\notin Z(G)$ and the index of $G'$ in $Z(GL(2,q))\cdot SL(2,q)$ is odd).

Set $t=g^{2^m}$, so that $\overline{t}=h^{2^m}$, and let $D$ be the image in $\overline{G}$ of the group of diagonal matrices in $SL(2,q)$. Then $|D|=(q-1)/2$, which is odd. 
Note that $N_{\overline{G}'} (D)$ contains an involution inverting the elements of $D$. Since all the involutions of $\overline{G}'$ are conjugate to each other, we may assume that  this inverting involution is $\overline{t}$. 
In addition, note that $C_D(\overline{t})=1$.
 Let $D=\lan d\ran$.
By Lemma \ref{aq1}, items (1),(3), there are $|d|$ distinct 
eigenspaces $M_\lam$ of $d$ on $M$, where $\lam$ is an \ei of $d$. 
As  $|d|$ is odd,  $\overline{t}(M_\lam)=M_{\lam\up}$ and $\overline{t}(M_\lam)= M_\lam$
implies $\lam=1$. \itf the Jordan form of $\overline{t}$ on $M$ is
$\diag(J_2\ld J_2, H)$, where $H$ is  the Jordan form of the matrix of $\overline{t}$ on $M_1$.

Moreover, by Lemma  \ref{aq1},  $\dim M_1=2$ if $\dim M=q-1$, whereas 
$\dim M_1=1$ if $\dim M = (q-1)/2$. 
In the latter case $H=\Id_1$.  
In the former case   
$H$
is either $J_2$ or $\Id_2$. We  show that  
 $H=\Id_2$ actually holds.

Suppose the contrary. Then the Jordan form of $\overline{t}$ is $\diag(J_2\ld J_2)$. Let $K$ be the Jordan form of $h$ on $M$.  
 It follows from Lemma \ref{jj2} that 
$K=\diag(J_{2^{m+1}}\ld J_{2^{m+1}})$. As $(q-1)/2$ is odd, this can only happen when $m=0$, that is, $g=t$. However, this implies that $t^2=g^2 \in Z(G)$, against our assumptions. 
Thus, $H=\Id_2$.

It also follows from Lemma \ref{jj2} that  
$K=\diag(J_{2^{m+1}}\ld J_{2^{m+1}}, J_{2^{m+1}-2})$ or $K=\diag(J_{2^{m+1}}\ld J_{2^{m+1}}, J_{2^{m+1}-1},J_{2^{m+1}-1})$  (resp., 
$K=\diag(J_{2^{m+1}}\ld J_{2^{m+1}}, J_{2^{m+1}-1})$). As $K$  is supposed to be almost cyclic, we must have $K=J_{2^{m+1}-2}$ (resp.,  $K=J_{2^{m+1}-1}$). Therefore,  $q-1=2^{m+1}-2$ (resp., $(q-1)/2=2^{m+1}-1$). So  $q+1$ is a 2-power, and $|h|=2^{m+1}=q+1$ (resp., $(q+1)/2$), as claimed.

The "if" part. We are now given that $q+1$ is a 2-power and $|h|=q+1$ (resp., $(q+1)/2$). Observe that the possible shapes of $K$ given in the previous paragraph do not depend on the assumption that $K$ is almost cyclic, but only on Lemma \ref{jj2} and the assumption that $(q-1)/2$ is odd.  If $\dim M=q-1$, then $|h|=2^{m+1}=q+1$, and hence the only option is  $K=J_{2^{m+1}-2}$. Otherwise, $|h|=2^{m+1}= (q+1)/2$, and $K=J_{2^{m+1}-1}$.

\begin{lemma}\label{ss2} 
Let $SL(2,q)\subseteq G\subseteq GL(2,q)$, where $q\equiv -1\pmod 4$. Let $g \in G
$ be a non-central $2$-element, and let $h$ be the projection of $g$ into $G/Z(G)$.
Let $M$ be an \ir $ FG $-module with  $\dim M>1$. Then $g$ is almost cyclic on $M$  if and only if:

$(1)$ $\ell \neq 2$, $g \in SL(2,q) \cdot Z(GL(2,q))$,  $\dim M=(q\pm 1)/2$ and $|h|= (q+1)/2$;

$(2)$ $\ell \neq 2$,  $g \notin SL(2,q) \cdot Z(GL(2,q))$,  $\dim M= q$ or $q \pm 1$ and  $|h|=q+1$;

$(3)$ $\ell \neq 2$, $q=7$, $\dim M =3$ and $|h| = 2$.

$(4)$ $\ell=2$, and one of the following holds:

$i)$ $\dim M=q \pm 1$ and $|h|=q+1$ (here the case $\dim M=q + 1$ only occurs for $g \notin SL(2,q) \cdot Z(GL(2,q))$);

$ii)$ $\dim M= (q-1)/2$ and $|h|= (q+1)/2$.

$iii)$ $q=7$, $\dim M =3$ and $|h| = 2$.

Additionally,   in all the above cases $q$ is a Mersenne prime.
\end{lemma}

Proof:  Let $\tau$ be the representation of $G$ afforded by $M$, and suppose that $\tau(g)$ is almost cyclic. 

First of all, observe that, by the Remark preceding the statement of Lemma \ref{q14}, we may assume 
that either $G=SL(2,q)$ or $G=GL(2,q)$.  Recall  that, by Lemma \ref{sq2},
every \ir $FG$-module lifts to characteristic zero.  Hence, if $\ell\neq 2$,
 it is enough to verify the lemma for $\ell=0$, which can be done examining the character 
table of $G$. This yields  items (1), (2) and (3) of the statement. (In (1), for $q=7$, $|g| = 8$.)

So, from now on, we may assume that $\ell=2$.

Suppose first that  $g^2\in Z(G)$. As $\ell=2$, then $\tau(g)$ acts as a transvection on $M$. It follows that case (2) of Lemma \ref{tt2}  holds, and hence $G=SL(2,7)$. In this case, $\dim M = 3$, $|g|=4$ and $|h|=2$, which gives item (4), $iii)$ of the statement.  
 
 Thus, from now on we may assume  that $g^2\notin Z(G)$. 
 Note that this implies that  $g$ is \ir on the natural module of $G$ (otherwise $|g|$ divides $q-1$, and
hence  $g^2=1$, as $(q-1)/2$ is odd by assumption). 

It is well known that the irreducible $2$-modular representations of $PSL(2,q)$ of non-zero defect are of degree $1$, $q-1$ or $(q-1)/2$ (see \cite{Bn}). 
Thus the case $G=SL(2,q)$, $\ell=2$
is dealt with in Lemma \ref{nn1}, provided $M$ has non-zero defect as a $PSL(2,q)$-module. Now, recall that a Sylow $2$-subgroup $P$ of 
$PSL(2,q)$ is dihedral of order dividing $q+1$ (as $q\equiv -1\pmod 4$). So, we are left to examine the case where $\dim M = q+1$. Observe that
$M|_P$ is a projective $FP$-module, and hence
$M|_P=m\cdot \rho_P^{reg}$ for some integer $m>0$. As $\tau(G) \cong PSL(2,q)$, we may assume that $\tau(g)\in P$. It follows that the matrix of $\tau(g)$ is not almost cyclic. Indeed,  the Jordan form of the matrix of $\tau(g)$ consists of
$d:=m\cdot |P|/|\tau(g)|$ blocks of equal size. As $P$ is not
cyclic, $d>1$. 

By the above, we may now assume 
that $G=GL(2,q)$ (and $\ell=2$). Set $G'=SL(2,q)$.  

Note that $\dim M$ is not of
degree $q$,  $(q+1)/2$. Indeed suppose the contrary. Then, by Clifford's theorem,
either $M|_{G'}$ is irreducible (which is not the case, e.g. see  \cite[pp. 90-91]{Bu}), or $M|_{G'}=M_1\oplus M_2$, where
$M_1,M_2$ are \ir $FG'$-modules of the same dimension. But this is impossible, considering the degrees of the  \ir $FG'$-modules. 
Next, suppose  that $\dim M=q+1$. Then $M|_{G'}$ is irreducible, by similar reasons, and we may assume that either $\tau(g) \in  P$ or $\tau(g^2) \in P$. In the former case $\tau(g)$ is not almost cyclic, as seen above.
So, let $ \tau(g)\notin P$. Then $\tau(g^2)\in P$, so
the Jordan form of the matrix of $\tau(g^2)$ consists of
$d:=m\cdot |P|/|\tau(g^2)|$ blocks of equal size. It then follows that the Jordan form of the matrix of $\tau(g)$ consists of
$d:=m\cdot |P|/|\tau(g)|$ blocks of equal size, and hence
is almost cyclic \ii $m=1$ and $|\tau(g)|=q+1$. This is part of item $(3), i)$ in the statement of the Lemma.

Finally, we are left to examine the cases where  $\dim M=q-1$ or $(q-1)/2$. These are dealt with in Lemma \ref{nn1}, and the results are stated in item $(4), i)$ and $ii)$ of the statement.

\medskip

At this point, we are in a position to prove Theorem  $\ref{mt2}$, stated in the Introduction:

\medskip

{\bf Proof of Theorem $\ref{mt2}$}.
If $p>2$, item (1) of the statement follows from Lemma \ref{cn2}. So, suppose that $p=2$. We distinguish two cases, according to $q$ being $\equiv 1\pmod 4$ or $\equiv -1\pmod 4$. If $q\equiv 1\pmod 4$, the claims in item (2) of the statement follow from Lemma \ref{q14}. If  $q\equiv -1\pmod 4$, the claims in item (2) of the statement follow from Lemma \ref{ss2}. Thus, the theorem is proven.

\medskip

The following Lemma and its Corollary show that the results stated in Theorem $\ref{mt2}$ carry over to any group $G$ such that  $SU(2,q)\subseteq G\subseteq U(2,q)$.

\begin{lemma}\label{un2}   Let $G=GL(2,q)$ and $H=U(2,q)$, for $q>3$. Then there exists a group $X$ such that
$X=Z(X)\cdot G=Z(X)\cdot H$. 
\end{lemma}

Proof. Let us consider the groups $G$ and $H$ as naturally embedded subgroups of the algebraic group $\overline{G}=GL(2,\overline{F}_q)$. Recall that $G'\cong H'$. By the general theory of representations of Chevalley groups (see also \cite{Bn}, Chapter 10), $G'$ and $H'$ are conjugate in $\overline{G}$. So, up to taking a suitable conjugate of, say, $G$ within $\overline{G}$, we may assume that $G'=H'$.
Set $X=\lan G,H\ran$, so that $X'=G'=H'$. Let $x\in X$, $g\in G'$. Then $xgx\up \in G'$. Let $T$ be a split torus in $G'$, which is a conjugate in $G'$ of the group of  diagonal matrices in $SL(2,q)$.  Then $xTx\up$ is another split torus, and it is conjugate to $T$ in $G'$ (as split tori are conjugate). So we can assume that $xTx\up=T$. It is then easy to check that $x$ is of shape
$\diag(a,b) $ or $\begin{pmatrix}0&a\\ b&0\end{pmatrix}$, where $a,b\in \overline{F}_q$.
Take $g= \begin{pmatrix}0&1\\ -1&0\end{pmatrix}$. Then $xgx\up$ equals $\pm \begin{pmatrix}0&b\up a\\ - ba\up&0\end{pmatrix}\in G'$ in both cases. Therefore, $b\in aF_q$, so $x\in Z\cdot GL(2,q)$, where $z$ is a scalar matrix.
This shows that $X=Z(X)\cdot G$. Next, we show that  $X=Z(X)\cdot H$. Suppose first that $q$ is even. Then $G=Z(G)G'$, whence, as $Z(G)$ and $Z(H)$ are both contained in $Z(X)$, $X = Z(X)G' = Z(X)H' = Z(X) \cdot H$, as claimed. Next, suppose that $q$ is odd. Then $G \nsubseteq Z(X) \cdot G'$. Indeed, let $g \in G$, with $ \det(g)$ a non-square in $F_q$, and assume that $g = zg'$, where $z =\diag (x,x)$, $x \in \overline F_q$, and $g'  \in G'$. Then $x \in F_q$ and $ \det(zg') = x^2$, a contradiction. It follows that $|G : Z(G)G'|=2$, whence also $|X : Z(X)G'|=2$ (as $X=Z(X)G$). Since $G'=H'$, 
it now suffices to show that $H$ is not in $Z(X)G'$. For this, observe that the latter group has a non-trivial complex representation of degree $(q-1)/2$, whereas $H$ does not, unless $q=3$ (e.g. see  \cite{En}) .

\begin{corol} \label{un3}  Let $\rho$ be an \irr of $H=U(2,q)$. Then $\rho$ extends to $X$. Moreover, if $h\in H$
then there exists $g\in G=GL(2,q)$ such that $\rho(h)=\lam \rho(g)$, where $\lam \in F$. In addition, if $ h^k\in Z(U(2,q))$   
then $g^k\in Z(GL(2,q))$, that is, the order of $g,h$ modulo centres are the same. 
 \end{corol}

The following results will be needed for the proof of Proposition  \ref{35} in Section $5$.

\begin{lemma}\label{ss3}  Let $G=SL(n,q)$, where 
$(n,q)\in \{(3,2),(3,4),(4,2)\}$. Let $g$ be a semisimple element of $G$ of $p$-power order, $p$ a prime,  and let
$M$ be an \ir $FG$-module with $ \dim M>1$, on which the matrix of
$g$ is almost cyclic. Then one of the \f holds: 

\medskip
$(1)$  $G=SL(3,2)$, and   either  $|g|=7$ and $M$ is arbitrary, or
$|g|=3$ and $\dim M=3$.

\medskip
 
$(2)$  $G=SL(4,2)$, $|g|\in\{3,5,7\}$ 
and  $\dim M=7$. (Here, if $|g|=3$, then
$g$ belongs to the class $3A$, in the Atlas notation).
\end{lemma}

Proof.  If $(|g|, \ell)=1$, then the result follows by inspection
of the Brauer characters of $G$ (see  \cite{MAtl}). Therefore, we may  assume that $\ell$
divides $|g|$.

\medskip
(1) Let $G=SL(3,2)$. If $|g|=\ell=7$ then, as   $SL(3,2)\cong
PSL(2,7)$, $M$ is realized as a $PSL(2,7)$-module, and the result
follows from the well known fact that a unipotent element $g\neq 1$
of $SL(2,\ell)$ in every \irr of this group in characteristic
$\ell$ is represented by a single Jordan block, and hence the matrix of
$g$ is cyclic. So, let $|g|=3$. Then $\dim M\in \{3,6,7\}$. Thus, by Lemma \ref{ww1} and
Corollary \ref{ddw}, $\dim M= 3$.

  \medskip
(2) Let $G=SL(4,2)$. In this case $|g|\in \{3,5,7\}$, and, for $p>3$, the Sylow $p$-subgroups $S$ of $G$ are cyclic of order $p$. 

First, let
$\ell =p=7$. Then the minimum dimension of $M$ equals $7$ (see \cite{MAtl}), in which case $M$ has defect zero. It follows from Lemma \ref{ww1} that  the matrix of $g$ on $M$ is a single Jordan block, and hence cyclic. Otherwise, $M$ has defect 1. In this case, as $N_G(S)/S$ is abelian, it follows from Corollary \ref{ddw} that $\dim M\leq 8$. However, for $\ell=7$ the only \ir $F$-\rep of dimension at most
$8$  is that of dimension 7, a contradiction.

 Next, let $\ell =p=5$. Then $g$
is regular; hence, by Lemma \ref{g12},  $G$ can be generated by three suitable conjugates
of  $g$. By Lemma \ref{nd5}, $\dim M\leq 12$. As $G$ has no irreducible $F$-representation
of degree $d$ for  $7<d<13$ (see \cite{MAtl}), it follows that $\dim M=7.$ The same
is true for $p=3$. Indeed,
 by Proposition \ref{GS2}, $G$ can be generated by at most 4 conjugates
of $g$; this implies $\dim M\leq 8$, by Lemma \ref{nd5}. It follows that
$\dim M=7$ (see \cite{MAtl}).

 Now, for both $p = 3$ and $p = 5$, there is a unique $FG$-module of
dimension 7. \itf  $M$ is isomorphic to the
only non-trivial constituent of the $8$-dimensional permutation
module for the alternating group $A_8\cong SL(4,2)$. For $p=5$, this implies that $g$ is almost cyclic on $M$. So, let $p=3$. There are exactly two conjugacy classes of elements of order three in $G$, labelled 3A and 3B in \cite{Atl}. The class $3A$  is
represented by a permutation fixing 5 out of 8 points.   \itf
$g\in 3A$ is almost cyclic on $M$. Next, suppose that $g\in 3B$. Then $g$ is
contained in a subgroup $X\cong A_7$, and $M|_{X}$ contains a non-trivial
$6$-dimensional constituent $N$, say, which is also a constituent
of the $7$-dimensional permutation module for $A_7$. We claim that $g$ is not almost cyclic on $M$, and for this it suffices to
show that $g$ is not almost cyclic on $N$. Suppose the contrary. Observe that $X$ can be
generated by two suitable elements from $3B$  (direct computation using GAP). It then follows from Lemma \ref{nd5} that  an irreducible constituent of $N$ must have dimension $ \leq 4$. However, the minimal dimension of a non-trivial $3$-modular representation of $A_7$ equals 6. This yields a contradiction.

\medskip
(3) 
 Let $G=SL(3,4)$. Here $p\in \{3,5,7\}$, and the Sylow 5- and
7-subgroups $S$ of $G$ are cyclic of prime order. Moreover, $N_G(S)/S$ is
abelian. Let first $p=\ell=5$. If $M$ has defect zero, then $\dim M  \geq 15$, and hence $g$ is not almost cyclic on $M$ by Lemma \ref{nd5}. Otherwise, by Corollary \ref{ddw}, $g$ almost cyclic implies $\dim M\leq 6$. However, $G$ has no non-trivial $F$-\reps of such degrees. Now, let $p=\ell=7$. If $M$ has defect zero, then $\dim M  \geq 21$, and again $g$ is not almost cyclic on $M$ by Lemma \ref{nd5}. Otherwise, by Corollary \ref{ddw}, $g$ almost cyclic implies $\dim M\leq 8$. Again, $G$ has no non-trivial $F$-\reps of such degrees. Finally, let $p=\ell=3$. In this case, it suffices to deal with the group $H=PSL(3,4)$. Let $h$ be the projection of $g$ into $H$. Then $|h|=3$, and by Proposition \ref {GS2} $H$
is generated by three suitable conjugates of $h$. Hence $\dim M\leq 6$,  by Lemma \ref{nd5}. But again,  there is no $3$-modular \irr of $H$ of this degree. (Note, however, that  the universal covering of $G$ has \ir  5-modular \reps of degree 6, as well as \ir  7-modular \reps of degree 6 and 8, and
$g$ is almost cyclic on these modules, by Lemma \ref{ww1}).

\medskip

Remark: Observe that the representation of $SL(3,2)$ afforded by the $FG$-module $M$, where $\dim M = 3$, is $ not$ Weil, according to our definitions (see above).

\begin{lemma}\label{sz33}  Let $SL(n,q)\subseteq G\subseteq GL(n,q)$,  where $n>2$ and $(n,q)\neq(3,3),(4,3)$,
and let $g\in G$ be a non-scalar semisimple element of $p$-power order, $p$ a prime. 
Suppose that $g$ stabilizes a $1$-dimensional subspace on the natural module of $G$.
Let $M$
be an \ir $FG$-module with $\dim M>1$, affording the representation $\phi$. Then $\phi (g)$ is
 not almost cyclic, unless one of the \f holds:

\medskip
$(1)$  $G=SL(3,2)$,    
$|g|=3$ and $\dim M=3$.

\medskip
$(2)$  $G=SL(4,2)$, $|g|\in\{3,7\}$, where $g\in 3A$ if $|g|=3,$
and  $\dim M=7$.
\end{lemma}

 Proof.  Suppose that $\phi(g)$ is almost cyclic, and assume that $g\in P$, where $P$ is the stabilizer of a $1$-dimensional subspace.
Let $U$ be the unipotent radical of $P$, and let  $\tau$ be
an  \ir constituent of $\phi |_P$ non-trivial on $U$. Note that $\tau$ is faithful on $U$: indeed, any subgroup $W$ of $U$ on which $\tau_W = Id$, would be normalized by $P$; however, $P$ acts transitively on $U\setminus \{1\}$ by conjugation. Let $T$
be the $FP $-module afforded by $ \tau $. Then $T|_{ U}=\bigoplus
T_{\kappa }$, where $\kappa $ runs over the group $K$ of
$F$-characters of $U$, and $T_{\kappa }=\left\{ t\in T\mid
ut=\kappa (u)t,\forall u\in U\right\} $. Moreover, the action of $P$ on $U$
by conjugation is dual to the action of $P$ on $K$.
Let $h=g^s\notin Z(G)$, where $s$ is such
that $h^p\in Z(G)$ (so $g^{ps}\in Z(G)$). Let $\phi (g^{ps})=\lam
\cdot \Id$. It is straightforward to check that $[h,U]\neq 1$. 
As $U$ acts scalarly on every $T_\kappa$, there is
$\kappa$ such that  $T_\kappa\neq
0$ and $hT_\kappa \neq T_\kappa$ (otherwise $\tau ([h,U])=1$, and hence $[h,U]=1$ as $\tau(U)\cong U$). It
follows that the $g$-orbit containing this $\kappa$ is of size
$ps$. Set $d:=\dim T_\kappa$ and $R=\oplus _{\nu\in \{g^i\kappa\}}T_\nu$.  If $p=\ell$, then the matrix of $g$ on $R$ is similar to the sum of $d$ Jordan blocks $J_{ps}$.  If $p\neq \ell$, then  all the $ps$-roots of $\lam$ are eigenvalues of
$\tau (g)$, each  with multiplicity at least
 $\dim T_\kappa$. Therefore,  $d=1$, since $\phi(g)$ is assumed to be almost cyclic. Furthermore, observe that, if $\tau'$ is another \ir constituent of $\phi |_P$ non-trivial on $U$, then we reach the same conclusion.
As $\phi(g)$ is assumed to be almost cyclic, we conclude that $\tau$ is the only \ir constituent of $\phi |_P$ non-trivial on $U$. It follows that $T':=\bigoplus_{\kappa\in K\setminus \{1_U\}}
T_{\kappa }$ must be an \ir $FP$-module of dimension at most $|K|-1=|U|-1=q^{n-1}-1.$

Now, let $T_1$ denote the subspace $T_\kappa$ with $\kappa=1_U$. Clearly, $T_1$
can be viewed as an $F(P/U)$-module, and by the above $M|_{P}=T'\oplus T_1$, where $T'$ is irreducible.
Observe that $L:=P/U$ is isomorphic to a subgroup of the group $X:=GL(n-1,q)\times GL(1,q)$ containing $SL(n-1,q)$. As  $X/Z(X)\cong PGL(n-1,q)$,
it follows that every normal subgroup  
of $L$ either is contained in $Z(L)$, or it contains $L'\cong SL(n-1,q)$, unless $(n-1,q)=(2,2),(2,3)$. These exceptions, however, do not occur here, as $(n,q)\neq(3,2)$ by Lemma \ref{ss3} and 
$(n,q)\neq(3,3)$ by assumption.

 It was shown in  \cite[p.237]{SZ}, that $P\cap SL(n,q)$ has an irreducible
constituent on $T_1$ of dimension greater than 1, unless $(n,q) \in \left\{(4,2),(3,2),(4,3),(3,4)\right\}$. Observe that the exceptional cases where $(n,q) \in \left\{(4,2),(3,2),(3,4)\right\}$ were dealt with in Lemma \ref{ss3}, yielding items $(1)$ and $(2)$ of the statement, whereas the case $(n,q)=(4,3)$ is ruled out by assumption. Therefore, from now on, we may suppose that $P\cap SL(n,q)$ has an irreducible
constituent on $T_1$ of dimension greater than 1. 
As $U\subset SL(n,q)$, it follows that $P$, and hence $L$, has an irreducible
constituent on $T_1$ of dimension greater than 1.

Observe that, since $g$ is almost cyclic on $M$, $g$ must act scalarly on $T_1$. 
[Otherwise, $g$ would have either a  non-trivial Jordan block on $T_1$ (if $\ell=p$), or
at least 2 distinct eigenvalues on $T_1$, which are also $ps$-roots
of $\lam$ (if $\ell \neq p$). But this would contradict the almost cyclicity of $g$ on $M$, in view of the action of $g$ on $T'$, as described above.] So, we may suppose that $g$ acts on $T_1$ scalarly, and hence that $\rho(g)$ is scalar. Let $N=\{a\in L: \rho(a)$ is scalar$\}$. Clearly, $N$ is a normal subgroup of $L$. So, either
$N\subseteq Z(L)$ or $N$ contains $L'$. The latter cannot happen, as $L/L'$ is abelian, and hence $\rho$ would be one-dimensional, which is false.
So $N\subseteq Z(L)$, and hence $g{\mod U}$ $\in Z(L)$. Let us consider the action of $P$, and hence of $L$, on $U$ by conjugation. Then, viewing $U$ as a vector space over $F_q$, $Z(L)$ acts on $U$ scalarly, and the kernel of the action of $L$ is $Z(G)$.
It readily follows that all the $g$-orbits on $U$, but one, have the same size $ps > 1$, and the number of non-trivial $g$-orbits is at least $(q^{n-1}-1)/(q-1)>1$.
Clearly, this remains true for the action of $g$ on $K$. However, as shown above, $g$ must have only one non-trivial orbit on $K$, which gives a contradiction.

\medskip
In the next two Lemmas we deal with the cases where $(n,q) \in \left\{(3,3),(4,3)\right\}$, which were left open in Lemma  \ref{sz33}.

\medskip

\begin{lemma}\label{333}   Let $SL(3,3)\subseteq G\subseteq GL(3,3)$, and let
 $g\in G$ be a non-scalar semisimple element of $p$-power order, for some prime $p$. Let
$M$ be an \ir $FG$-module with $\dim M >1$, affording a representation $\phi$. Then  the matrix of
$g$ on $M$ is not almost cyclic, unless one of the following holds:

$(1)$ $\ell \neq 2, 13$, $|g|=13$ and $\dim M = 12$ or $13$ (in which cases $g$ is cyclic on $M)$;

$(2)$ $\ell = 2$, $|g|=13$ and $\dim M = 12$ (in which case $g$ is cyclic on $M)$;

$(3)$ $\ell = p = 13$, $|g|=13$ and $\dim M=13$ (in which case $g$ is cyclic on $M)$;

$(4)$ $\ell = p = 13$, $|g|=13$ and $\dim M=11$ (in which case $g$ is cyclic on $M)$.
\end{lemma}
Proof. Observe that, since $GL(3,3)=SL(3,3)\times \{\pm \Id\}$, we may assume that $G=SL(3,3)$. 
Here $p \in \left\{2,13)\right\}$. Suppose first that $p=13$. If $\ell \neq 2,13$ or $\ell=2$, then items $(1)$ and $(2)$ of the statement follow by direct inspection of the character table of $G$ and \cite{MAtl}, respectively. So, suppose that $\ell =p= 13$. If $M$ has defect zero, then $\dim M \in \left\{13,26,39\right\}$. Hence $g$ is almost cyclic on $M$ precisely when $\dim M = 13$, by Lemma \ref{ww1}. If $M$ has positive defect, then $\dim M \in \left\{11,16\right\}$. As $N_G(\langle g\rangle )/\langle g\rangle $ is abelian, Corollary  \ref{ddw} rules out the case $\dim M = 16$, while direct computation using MAGMA shows that $g$ is cyclic on $M$ when $\dim M = 11$. This gives items $(3)$ and $(4)$ of the statement.

 Next, suppose that $p=2$. If $\ell \neq 2$, then the statement follows by inspection of the character table of $G$ and  \cite{MAtl}. So, let $\ell =2$. If $g$ is an involution, then the claim follows from Lemma \ref{wz1} (as $\phi(G)$ is not generated by transvections). Suppose that $g^2\neq 1$. Then one observes that
$C_G(g)$ contains no element of order 3, that is,  $g$ is regular. By Lemma \ref{g12}, 
$G$ is generated by three conjugates of $g$. Then $\dim M\leq
3(|g|-1)$ by Lemma \ref{nd5}. As the minimum dimension of a
non-trivial $F$-\rep of $G$ is 12, it follows that $|g|=8$, and
$\dim M\leq 21$.
So   $\dim M\in \{12,16\}$ (see \cite{MAtl}). 
As the order of a Sylow 2-subgroup $S$ of $G$ is 16, the \reps of degree 16 are of defect zero, and hence $\phi|_S=\rho_S^{reg}$, by Lemma \ref{ww1}.
\itf $\phi(g)$ is not almost cyclic.  If $M$ has dimension 12, then the claim follows by direct computation, using MAGMA.

\begin{lemma}\label{s43}    Let $SL(4,3)\subseteq G\subseteq GL(4,3)$, and let
 $g\in G$ be a non-scalar semisimple element of $p$-power order, for some prime $p$. Let
$M$ be an \ir $FG$-module with $\dim M >1$.
Then  the matrix of
$g$ on $M$ is not almost cyclic.

\end{lemma}

Proof. By way of contradiction, assume that the matrix of $g$ is almost cyclic on $M$. Recall that the minimum dimension of a (projective)
 irreducible $F$-representation of $G$ is 26. Suppose first that $p>2$. Then
$p\in \{5,13\}$. As $|G/G'|\leq 2$, it suffices to verify the lemma for $G=SL(4,3)$. (Indeed, we may assume that $g \in G'$. Moreover, as $g$ is almost cyclic on $M$, $g$ must be almost cyclic on any constituent of $M_{|G'}$.) Observe that $g$ is regular. For $p=5$, it follows from Lemma \ref{g12} and Lemma \ref{nd5} that $\dim M \leq 12$, which is a contradiction. So, let $p=13$. Then, by  Lemma \ref{g12} and Lemma \ref{nd5}, $\dim M = 26$. If $\ell \neq 13$, a direct inspection of the character table and the Brauer characters of $G$ in \cite{Atl} and  \cite{MAtl} shows that $g$ is not almost cyclic on $M$. If $\ell=13$, then $M$ has defect zero, and hence $g$ is not almost cyclic on $M$ by Lemma \ref{ww1}.

Next, let $p=2$, and let $V$ be the natural module for $G$. Note that $g^8\in Z(G)$. (Indeed, if $g$ is \ir,  then $|g|\leq16$, and hence $g^8=\pm \Id$. On the other hand, if $g$ is reducible, then $|g|\leq 8$.)
If $g$ is regular (that is, $C_G(g)$ contains no unipotent element), then $g$ is generated by three
conjugates of $g$, by Lemma \ref{g12}. But  then $\dim M\leq 21$ by Lemma \ref{nd5}, a contradiction.
So, suppose that $g$ is not regular. Then $g$ is reducible (by Schur's Lemma), and hence $|g|\leq 8$. If $g^4\in Z(G)$, then, by Proposition \ref{GS2}(1), $
G$ is generated by four suitable conjugates of $g$. 
Hence $\dim M\leq 12$ by Lemma \ref{nd5}, again a contradiction.
 So, we may assume that $|g|= 8$.  Then
 $V=V_1\oplus V_2$   (a direct sum decomposition), where $\dim V_i=2$ and $gV_i=V_i$
for $i=1,2$. Set $g_i=g|_{V_i}$.
 If both $g_1,g_2$ are of order 8, then $g^4\in Z(G)$, which case has been already ruled out.
So we may assume that $|g_1|=8$ and $|g_2|\leq 4$. Then $g_1$ is \ir on $V_1$, and hence both the \eis of $g_1$
on $V_1\otimes \overline{F}_q$ are primitive 8-roots of unity. As $g$ is not regular, it follows that the \eis of $g_2$
on $V_2\otimes \overline{F}_q$ are not distinct, whence $g_2=\pm \Id$.
So $g$ stabilizes a direct sum decomposition of $V$, say $V=W \oplus U$, where $\dim W = 1$.  Let $H$ denote the stabilizer in $G$ of both $W$ and $U$, so that $g \in H$. If $G=SL(4,3)$, then $H \cong GL(3,3)$,  whereas if
$G=GL(4,3)$, then $H\cong GL(3,3)\times Y$, where $Y= \left\{\pm 1 \right\}$. As $g\in H$,
 the result follows from Lemma \ref{333}.

\begin{lemma}\label{423} 
Let  $G =Sp(4,3)$, and $g\in G$ be a non-scalar semisimple element of $p$-power order, $p$ a prime.
Let $\phi$ be an \ir $F$-representation of $G$. Then the matrix $\phi(g)$
is almost cyclic if and only if one of  the \f holds:

$(1)$ $p=2$, $\ell\neq 2$ and

$(i)$    $|g|=2$ 
and $\dim\phi=5$;

$(ii)$ $|g|=4$,  $g^2\notin  Z(G)$   and $\dim\phi=4$;

$(iii)$ $|g|=8$,  $g^4\in Z(G) $ and $\dim\phi=4$ or $5$.

\noindent Furthermore,   the matrix of   $\phi(g)$ is cyclic only if $|g|=8$ and $\dim\phi=4$.

$(2)$ $p=5$ and $\dim\phi\in \{4,5,6\}$,  where $\dim\phi\neq 6$ if   $\ell= 3$ and $\dim\phi\neq 5$
 if $\ell= 2$. Furthermore, $\phi$ is faithful if and only if  $\dim\phi=4$ and $\ell\neq 2$.

$(3)$ $p=\ell=2$ and $\dim\phi=4$. In addition, either $|g|=4$, or $|g|=2$ and $\phi(g)$ is a transvection in
$SU(4,2)$.
\end{lemma}

Proof.  
First, let $p>2$. Note that a Sylow 5-subgroup $S$ of $G$ is of order 5, and we may assume $g\in S$. If $\ell\neq 5 $, then the claim in $(2)$ follows
 from a direct inspection of the Brauer character tables of $G$ in \cite{MAtl}.


So, let $\ell=5$. Observe that $C_G(g)$ has order 10, and hence $g$ is regular. It follows that $\dim\phi\leq 12$, by Lemmas  \ref{g12} and \ref{nd5}. Thus $\dim\phi \in\{4,5,6,10\}$. If $\dim\phi=5$ or $10$, then $\phi$ is of 5-defect zero, and hence by Lemma \ref{g12}
$\phi|_{S}= \rho_S^{reg}$ or $2\rho_S^{reg}$, respectively. Therefore, the matrix of $g$ is almost cyclic (in fact cyclic, represented by a single Jordan block $J_5$) only when $\dim\phi=5$. If $\dim\phi \in\{4,6\}$, then direct computation using MAGMA shows that $\phi(g)$ is cyclic, yielding $(2)$.

Next, let $p=2$ and $\ell\neq 2$. In this case the claim in $(1)$ follows by direct computation from the data in  \cite{Atl} and \cite{MAtl}.

If $p=\ell=2$, then $\phi$ can be viewed as a \rep of $SU(4,2)$.  It is easy to check, using MAGMA, that in both the classes $4A$ and $4B$ (Atlas notation) can be found two suitable elements generating $SU(4,2)$, and hence, for $g$ in these classes, we only need to examine $\phi(g)$ for $\dim\phi \leq 6$. It turns out that  $\phi(g)$ is almost cyclic only when $\dim\phi =4$ (almost cyclic in case 4A, cyclic in case 4B). Finally, if $|g|=2$, almost cyclicity only occurs when $g$ is a transvection, in which case five conjugates of $g$ are enough to generate the group. This gives $(3)$.

\med
{\bf Remark}. Recall that $O^-(6,2)=SO^-(6,2)\cong SU(4,2)\cdot C_2$. The group $O^-(6,2)$ is generated by transvections, and has an \irr of degree 6 over the complex numbers, in which there exists an element of order  2 represented by an almost cyclic matrix (it belongs to the class $2C$ in the notation of \cite{Atl}). (Of course, there are no transvections  in the commutator subgroup of $O^-(6,2))$. In addition,
 $O^-(6,2)\cong CSp(4,3)$, the conformal symplectic group (see \cite[p. 26]{Atl}), and $|{\rm Aut}\, G:G|=2$.
 
\begin{lemma}\label{442} 
Let  $G =SU(4,2)$. Let $g\in G$ be a non-scalar semisimple element of $p$-power order, $p$ a prime.
Let $\phi$ be an \ir $F$-representation of $G$ such that the matrix $\phi(g)$
is almost cyclic. Then one of the \f holds (we use the Atlas notation for conjugacy classes):

$(1)$ $p=3, \ell\neq 3$, $|g|=3$, $g\in 3D$ and $\dim\phi=5,$  or  $g\in 3C$ and $\dim\phi=6$;

$(2)$ $p=3, \ell\neq 3$, $|g|=9$, $g\in 9A,9B$ and $\dim\phi=5,6$;

$(3)$ $p=3, \ell= 3$,  $|g|=3$,  $g \in 3C,3D$ and $\dim\phi=5$;

$(4)$  $p=3, \ell= 3$,  $|g|=9$, $g\in 9A,9B$ and $\dim\phi=5;$

$(5)$ $p=5$ and $\dim\phi =5,6$.\end{lemma}

Proof. Note that $G\cong PSp(4,3)$. If $p =3$ and $\ell \neq3$, then $(1)$ follows from \cite[Lemma 4.2]{DMZ10}, where the reader can find more details. 

Now, let $p= \ell =3$. 
All the $3$-modular irreducible representations of $PSp(4,3)$ are available on the Atlas on line.  Easy routines using the MAGMA package yield the results listed in $(3)$ and $(4)$.

Finally, let $p=5$. Then the claim in $(5)$ follows from Lemma \ref{423}.

\begin{lemma}\label{452}
Let  $G =SU(5,2)$. Let $g\in G$ be a non-scalar semisimple element of $p$-power order, $p$ a prime.
Let $\phi$ be an \ir $F$-representation of $G$.  Then the matrix $\phi(g)$
is almost cyclic if and only if one of the following occurs:

$(1)$ $p=3, \ell\neq 3$, $|g|=9$, $g\in 9C,9D$ and $\dim\phi=10$;

$(2)$ $p=11$ and $|g|=11$ and $\dim\phi=10$ or $11$.

\noindent (Note that the representations occurring in $(1)$ and $(2)$ are Weil $F$-representations of $G)$.

 \end{lemma}
 
 Proof.  First, observe that, by Proposition  \ref{GS1},(1),  $G$ can be generated by at most five conjugates of $g$. By Lemma  \ref{nd5}, this implies that, whenever $|g|= 3, 5, 9$,  we only  need to examine the $F$-representations $\phi$ of $G$ of degree $10$ and $11$ (since any other irreducible $F$-representation of $G$ has degree $\geq 43$). On the other hand, the same holds when $|g|=11$; indeed, using the MAGMA package, it turns out that, for $|g|=11$, two suitable conjugates of $g$ are enough to  generate $G$.   Then the statement follows by direct computation using the Atlas and the MAGMA package. (Note that in item $(1)$, for $\ell=3$, $\phi(g)$ has Jordan form $\diag(J_8,J_2)$).

\medskip

We close this Section with the following result, which will be needed in the sequel  (see the proof of  Lemma \ref{nd5}).

\begin{lemma}\label{u62}  Let $G=U(6,2)$ and let $g\in G$ be an element of order $9$.
Let $\tau$ be an irreducible $F$-representation of G. Then $\tau(g)$ is not almost cyclic. 
\end{lemma}

Proof. 
Let $V$ be the natural module for $G$. Suppose first that $g$ is not contained in any proper parabolic subgroup of $G$; so, in particular, $g$
is regular. Observe that $g$ stabilizes  a $3$-dimensional subspace, obviously non-degenerate.
One readily observes that there exists an orthogonal basis of $V$ with respect to which $g$ has one of the following  shapes (where $\ep$ is a non-trivial cubic root of $1$):

\medskip
$g_1=\diag\big(\begin{pmatrix}0&1&0\\ 0&0&1\\ \ep&0&0\end{pmatrix},\begin{pmatrix}0&1&0\\ 0&0&1\\ \ep^2&0&0\end{pmatrix}\big)$ or 
$g_2=\diag\big(\begin{pmatrix}0&1&0\\ 0&0&1\\ \ep^i&0&0\end{pmatrix},\begin{pmatrix}1&0&0\\ 0&\ep&0\\0&0& \ep^2\end{pmatrix} \big)$ 

\medskip\noindent
for $i=1,2$ (but the matrices with $i=1,2$ only differ by a scalar, so we may assume $i=1$). 
Note that $g_1\in G'=SU(6,2)$, whereas $\det g_2\neq 1$; however, $g_2 \in G'$ up to a scalar. 

Suppose first that $g=g_1$. By Lemma \ref{d102}, given a semisimple element $x\in G'$ there are two conjugates  of $g$ whose product is equal to $x$. So, if we fix some element $x$ of order 11 in $G'$, we can assume that $x=gg'$, where $g'$ is conjugate to $g$. We claim that $g$ and $g'$ generate $G'$. Indeed, suppose the contrary. It suffices to show that $g$ and $g'$ are not contained in any maximal subgroup of $G'$.  Inspecting  the list of maximal subgroups $M$ of $G'/Z(G')$ (see  \cite{Atl}), we observe that 
11 is coprime to the  order of any such $M$, except for the case where $M \cong SU(5,2)$. Let  $M_1\cong U(5,2)$ be the preimage of $M$ in $G'$. Then $M_1$ (up to conjugacy) is the unique maximal  proper subgroup of $G'$ of order divisible by 11. So we may assume that  that  $x,g,g'\in M_1$. Now,  $M_1$ fixes a $1$-dimensional subspace of $V$, whereas $g = g_1$ does not fix any such subspace, since it has no eigenvalues on $V$. This is a contradiction. 
Thus, $G=\langle g,g'\rangle$. By Lemma \ref{nd5}, $\dim \tau\leq 16$. However, the minimum dimension of a non-trivial irreducible $F$-representation of $G'$ equals $21$. This completes the analysis of this case. 

Next, suppose that
$g=g_2$.  Since $|G:G'|=3$, $G=\langle g,G'\rangle$. Using the MAGMA package, one sees that there is a conjugate  $g'$ of $g$ such that $\langle g,g'\rangle=\langle g,G'\rangle =G$.  As above, $\dim \tau\leq 16$
by Lemma \ref{nd5}. So we have again a contradiction, as in the previous paragraph. 


Now, suppose that $g$ is not regular. Then $g$ is conjugate to an element $g_3$ of shape $$g_3=\diag\big(\begin{pmatrix}0&1&0\\ 0&0&1\\ \ep^i&0&0\end{pmatrix},\begin{pmatrix}\ep_1&0&0\\ 0&\ep_2&0\\0&0& \ep_3\end{pmatrix} \big),$$
where $\ep_1,\ep_2,\ep_3$ are $3$-roots of unity, not all distinct (corresponding to the $20$ classes of non-regular elements of order $9$ contained in $G$). \itf $g$ is contained in a parabolic subgroup  $P$, say, which stabilizes an isotropic $1$-dimensional subspace.
 Let $U$ be the unipotent radical of $P$. Then $U/Z(U)$ is an elementary abelian 2-group of order $4^4=2^8$.
Set $X=\langle g,U\rangle$, and let $\phi$ be 
an irreducible constituent of $\tau|_X$ non-trivial on $Z(U)$. Set $E=\phi(U)$. Then $E$ is a group of symplectic type, and $E/Z(E) \cong U/Z(U)$ has order $2^8$ (see \cite{DMZ0}). However, Lemma \ref{a55} implies that $|g|=2^4 \pm 1$, which is false, as $|g|=9$. 




\section{Almost cyclic elements in Weil \reps}

As mentioned in the Introduction, most non-trivial examples of almost cyclic matrices seem to arise
in Weil \reps of finite classical groups. In this Section we fully analyze such representations, with the aim of providing an exhaustive picture of the occurrence of almost cyclic matrices. 

The use of induction is an essential part of our machinery.
As we deal with  classical groups, the starting point of
 induction will be the study of elements 
that are orthogonally indecomposable on the underlying vector space.
This means that $g$ is an element of  a finite classical group $G$ which does not
stabilize any non-trivial
non-degenerate subspace of $V$, where $V$ is the natural module
of $G$  (in the case of $G=GL(n,q)$ or $SL(n,q)$ the word
`non-degenerate' must be dropped). This implies that one of two
situations holds: either $g$ is irreducible,  or $G\neq
GL(n,q),SL(n,q)$ and $V=V_1\oplus V_2$, where $V_1,V_2$ are
$g$-stable totally singular subspaces of $V$ (e.g. see
\cite[Satz 1 and 2]{Hu}). The orthogonally indecomposable case will be dealt with in Subsection 5.2.
Next, the case where the element $g$ is orthogonally decomposable must be treated. This will be done in  Subsection 5.3.

\subsection{Weil \reps}\label{sbs61}

We recall the notion and the basic properties of  Weil representations.

Let $E$ be an extraspecial $r$-group.  If $r$ is odd,  assume
$E$ to be of exponent $r$.
As always in this paper,  $F$ is an \acf of characteristic $\ell\neq r$, and $F_q$
is a field of order $q$, where $q$ is an $r$-power.
It is well known that $E$ has
faithful \ir $F$-representations,
all of them  of  degree $r^m$, where $|E/Z(E)|=r^{2m}$. Let us single out one of these \reps and identify $E$ with its image. Thus, $E$ is now an \ir subgroup of $GL(r^m,F)$.
 The commutator map $(a,b)\ra [a,b]$
yields a symplectic space structure on $V:=E/Z(E)$. Let $N$ be the normalizer of $E$ in $GL(r^m,F)$.
 Then the conjugation action of $N$ on $E$ preserves
the commutation in $E$, and hence yields  a homomorphism $\eta:N\ra Sp(V)$, which is known to be surjective.
This means that $E/Z(E)$ is isomorphic to the natural module for $Sp(V)$.
Now let $G$ be a non-trivial group. Suppose that there is an injective homomorphism  $j: G\ra N$ such that $j(G)\cap E=1$.
  Then $j$  yields  a \rep $G\ra GL(r^m,F)$, which is called
a {\it  generic Weil representation}, and whose \ir constituents
are called {\it Weil \reps} of $G$. In
practice, it is not reasonable to use this definition for an
arbitrary group $G$; so we assume that $\eta (j(G))$ stabilizes no
non-zero subspace of $V$.
Thus the groups $ Sp(2m,r) ~(r~ {\rm odd}),~
SU(m,r)$, $U(m,r)$, and $SL(m,r) $, $GL(m,r)$ are examples of the group $G$ in question (e.g.,  see \cite{Ge}).

In principle, a Weil \rep of a group $G$ as defined above depends on a faithful \rep of $E$ and on the embedding $j$.
 However,  if $G\in \{SL(m,r)$, $SU(m,r)\}$,  there is in fact only one (up to equivalence) generic Weil representation, whereas if  $G=Sp(2m,r)$, exactly two non-equivalent generic Weil representations can be obtained in this way.
If $G\in \{GL(m,r)$, $U(m,r)\}$,  one obtains several  generic Weil representations, but all of them differ from each other by
tensoring with a one-dimensional \rep of $G$ (e.g., see
 \cite{Ge}). In fact, this is immaterial for
our purposes.  

Additionally, we emphasize that every generic Weil \rep
of $G=GL(m,r)$ is the tensor product of the permutation
$F$-\rep of $G$, associated with the action of $G$ on the vectors of the
standard $F_rG$-module, with a 1-dimensional module.

Now, let $m=nk$ and set $q=r^k$, $k\geq1$. It is well known that there are embeddings $GL(n,q)\ra GL(m,r)$, $Sp(2n,q)\ra
Sp(2m,r)$ and $U(n,q)\ra Sp(2m,r)$ obtained by viewing $F_q$ or $F_{q^2}$ as  vector spaces over $F_r$. We call them  standard embeddings. Composing  each of these embeddings  with a \rep $j$ as defined above, one obtains generic Weil \reps of these groups, and again, the above comments remain valid by replacing $r$ by $q$. Namely, in this way one obtains 
exactly one generic Weil \rep for $SL(n,q)$ and $SU(n,q)$ (up to equivalence), and exactly two generic Weil \reps for $Sp(2n,q)$ (up to equivalence). Likewise, those for  $GL(n,q)$ and $U(n,q)$ can be
 obtained from each other by tensoring with a one-dimensional one.

\medskip
In this section we also use the term   Weil character
referring to the character (Brauer character)
of the $FG$-module afforded by a Weil \rep of $G$. 
It follows from the construction of a generic Weil   \rep   that its Brauer character
(when the characteristic $\ell$ of the ground field is prime) coincides with the
restriction to $\ell'$-elements of the Weil character in characteristic 0.
One can refer to
\cite{Ge} and \cite{War} for more details on the basic properties of Weil  representations.

\medskip

Each of the two ordinary (i.e.  $\ell=0$) generic Weil \reps $\phi$ of  $Sp(2n,q)$, $q$ odd, has two \ir constituents,
 $\phi^+,\phi^-$, say,   of dimension $(q^n+ 1)/2, (q^n- 1)/2$, respectively. They remain \ir under reduction to any
   characteristic $\ell >2$ coprime to $q$. For $\ell =2$ this is only true for $\phi^-$, while the reduction of $\phi^+\,$mod 2
   has two composition factors, one of them one-dimensional  (and in fact trivial
unless $(n,q)=(1,3)$), the other one  equivalent to
 $\phi^-({\rm mod}~ 2)$ (see  \cite{War}.)

As mentioned above, up to tensoring by a one-dimensional representation,
there is a unique ordinary generic Weil \rep
of $U(n,q)$; if $n>2$, it consists of $q+1$ composition factors,
not equivalent to each other. If $n$ is odd, then the dimensions of the
\ir constituents are $-1+\frac{q^n+1}{q+1}=\frac{q^n-q}{q+1}$ or  $\frac{q^n+1}{q+1}$. If $n$ is even, then the dimensions
of the \ir constituents are $\frac{q^n-1}{q+1}$ or
$1+\frac{q^n-1}{q+1}=\frac{q^n+q}{q+1}$.  These \ir constituents
remain \ir and pairwise  non-equivalent under restriction to
$SU(n,q)$. The representations of lower degree remain \ir under
reduction modulo any prime $\ell$ coprime to $q$. The other
representations remain \ir provided $(\ell ,q+1)=1$.  More
precisely, if $(\ell,q+1)\neq 1$, the following holds (see \cite[Proposition 9]{HM}.) Assume first that
$n$ is odd, and $\psi$, say, is an ordinary \ir Weil \rep of degree $\frac{q^n-q}{q+1}$.
 Then 
the reduction modulo $\ell $  of a representation of degree $\frac{q^n+1}{q+1}$ either remains \ir or it has two composition
factors, one of them 1-dimensional, the other one
 equivalent to $\psi\, ({\rm mod}~ \ell)$ tensored by a 1-dimensional one.
Next, suppose that  $n$ is even. Then an ordinary  Weil \rep of degree $\frac{q^n+q}{q+1}$ is
reducible modulo any prime $\ell $ dividing $q+1$; its reduction modulo $\ell$
has two \ir constituents, one  of dimension 1, the other one
of dimension $\frac{q^n-1}{q+1}$. In fact, it is known that every $\ell$-modular
\ir Weil \rep lifts to characteristic 0. (This follows from
results in \cite{DT, HM}, but it is not stated there  explicitly. For $n$ even,
see the last paragraph of the proof of \cite[Theorem 7.2]{DT}; for
 $n$ odd, see the proof of \cite[Proposition 9]{HM}.)

In the case where $G=GL(n,q)$, as noticed above, a generic Weil \rep  coincides with the permutation $F$-representation
of $G$ associated with the action of $G$ on the vectors of the natural $G$-module, up to tensoring
with a one-dimensional representation. It follows that the dimensions of the \ir
constituents of
a generic Weil \rep of $GL(n,q)$ are the same as those of the permutation \rep in question. These
are known to be $\frac{q^n-1}{q-1}$, $\frac{q^n-1}{q-1}-1,$ $\frac{q^n-1}{q-1}-2$ or $1$;
for details see for instance  \cite[Theorem 9.1.4]{GPPS}.

\medskip

Finally, in the following Lemma we state  a crucial property of Weil representations, concerning their restrictions to 'standard' subgroups:

\begin{lemma}\label{wr3}  Let $G\in\{GL(n,q), U(n,q)$,
$n>2$, $Sp(2n,q)$, $n>1$ and $q$ odd\} and let $V$ be the natural module for $G$.
Let $V=W\oplus W'$ be a  decomposition of V as a direct sum of subspaces,
where $W$ is non-degenerate if $G\neq GL(n,q)$, and set 
$S=\{g\in G | gW=W~{\rm and} ~gw'=w'$ for all $w'\in W'\}$. 
Let 
$\om$ be a generic, respectively  \ir Weil $F$-\rep of $G$.
Then $\om|_S$ is a direct sum of generic, respectively \ir Weil $F$-\reps of $S$.
\end{lemma}

Proof. The statement follows for arbitrary $\ell$ (coprime to $q$) if it holds for
$\ell=0$, by the very definition
of $\ell$-modular Weil representations. So let $\ell=0$.
It is known that  the restriction of $\om$ to $S$ is the sum of  generic Weil \reps of $S$.
(The proof is available in \cite{Z86}, and can be easily deduced
from properties of extraspecial $r$-groups and their representations. See also \cite[Proposition 2.2]{TZ2}.)
This immediately implies the claim for \ir Weil representations.

\medskip
At this point, it is worth to recall that every abelian subgroup $A$ of a finite classical group $G$ consisting of semisimple elements and orthogonally indecomposable,
is cyclic. If $A$ is irreducible and of maximal order,
then $A$ is called a Singer subgroup and  its generators are called
Singer cycles. If $n$ is even, $U(n,q)$ and $SU(n,q)$ do not have Singer cycles. Likewise,
$O^+(2n,q)$ and $O(2n+1,q)$ do not have Singer cycles.  If $G\in \{GL(n,q)$; $SL(n,q)$; $U(n,q)$, $n$ odd; $SU(n,q)$, $n$ odd; $Sp(2n,q)$; $O^-(2n,q) \}$, then the order of a Singer cycle is  known to be $q^n-1,(q^n-1)/(q-1)$,
$q^n+1$, $(q^n+1)/(q+1)$, $q^n+1$, $q^n+1$,  respectively.

Now, suppose that  $A$ is reducible. Clearly, by Maschke's theorem, such an $A$ cannot occur in the groups $GL(n,q)$ and
$ SL(n,q)$. So we assume that $G$ is not one of these two groups. It is well known (for
details, see \cite{Hup}) that $V$ is a direct sum of
 two maximal totally singular $A$-stable subspaces $V_1,V_2$ of equal dimension. So $V$ is of even
dimension and of Witt index $\dim V/2$ in the case of unitary and orthogonal groups. Furthermore,
$A$ acts irreducibly on both $V_1$ and $V_2$, and the actions of $A$ on these subspaces are dual  to each other.
 In particular,  if $G\in \{Sp(2n,q)$, $O^+(2n,q)$, $U(n,q)$, $n$ even$\}$,  then
$|A|$ divides $q^n-1$. Moreover, if $A$ is reducible and of maximal order, any generator of $A$
will be called a Singer-type cycle of $G$.

 An additional, simple but useful  observation  is that if an element  $g\in G$ 
is  semisimple and orthogonally indecomposable,
then it is a power of a Singer cycle or of  a Singer-type cycle. (This is a well known fact.
For detailed arguments see \cite[Lemmas 7.1 and 8.1]{EZ}.)

\subsection{Orthogonally indecomposable elements}

In this subsection we deal with the case when $g$ is a
semisimple and orthogonally indecomposable element of $G$.

As always, let $F$ be an \acf of characteristic $\ell\neq r$. Recall that
 $1_G$ denotes the trivial $FG$-module, and
$\rho^{{\rm reg}}_G$  the regular $FG$-module. 

\medskip
We first consider the  generic Weil \reps of $G$.

\begin{lemma}\label{p33} 
$(1)$ Let $G=Sp(2n,q)$,  where $n >1$ and $q$  is odd, or
  $G=U(n,q)$, where $n>1$ is odd.
Let $ S=\lan g \ran$, where $g$ is  a Singer cycle  in $G$,  and let $\phi$ be a  generic Weil $F$-\rep of
$G$. Then $\phi (g)$ is a cyclic matrix.  

\medskip
$(2)$ Let $G\in \{Sp(2n,q)$, where $n>1$ and $q$ is odd; $U(n,q)$,
 where  $n>2$ is even;  $GL(n,q)$, $n>2 \}$, and let $\phi$ be
a generic Weil $F$-representation of $G$.
Let $g$  be either a Singer cycle for $GL(n,q)$, or
a Singer-type cycle for $G\neq GL(n,q)$ (in each case
the order of $g$ equals $q^n-1$). Then $\phi (g)$ is an almost
cyclic matrix and $\deg \phi (g)=|g|$.

\end{lemma}

Proof. The cyclicity (respectively, almost cyclicity) of $\phi (g)$ in
case (1) (respectively,  case (2)) follows  from the definition of a Weil
representation and Lemma  \ref{ex7}, items (1) and  (2)(i), respectively.
Lemma \ref{ex7}(2)(i) also implies the claim on $\deg\phi(g)$ in (2). We
only have to observe that if $g\in U(n,q)$ is orthogonally
indecomposable (resp., $g\in GL(n,q)$ is irreducible), then $g$ is
orthogonally indecomposable in its action on $E/Z(E)$ when it is
viewed as a symplectic space. (Recall that the natural module for
$G=U(n,q)$ can be embedded into a symplectic space of dimension
$2n$ over $F_q$,
 preserving orthogonality. This is well known (e.g. see  \cite[4.3, p.117]{KL}).
In the case of $G=GL(n,q)$, the natural module can be embedded into
 $E/Z(E)$ as a maximal totally isotropic subspace.

\begin{corol}\label{bc44} 
 Suppose that $h\in \lan g\ran$, where $g$ is as in $(1)$ or $(2)$
 of Lemma $\ref{p33}$, and $g\in  \lan h, Z(G)\ran$.
Let $\tau$ be an \ir Weil \rep of $G$, with  $\dim\tau>1$. Then $\tau (h)$ is almost
cyclic. Furthermore, $\deg\tau(g)=\deg\tau(h)={\rm min}\, \{|g|/|Z(G)|,
\dim\tau\}$. In particular, $\deg\tau(g)\geq
\frac{|g|}{|Z(G)|}-1$.\end{corol}

Proof. Let $M$ be the $FG$-module afforded by a generic Weil
 representation $\phi$. It follows from
 Lemma \ref{p33}, (1) and (2), that $g$ yields an almost cyclic
 matrix in its action
on  every composition factor $M'$ of $M$.  Set $Z(G)= \lan z\ran$.
As $Z(G)$ acts scalarly on $M'$, the matrix of $gz^i$ on $M'$  is
almost cyclic too. As $g=z^jh^k$ for some $j,k$, the matrix of
$h^k$ on $M'$ is almost cyclic. This implies a similar claim for
$h$. Whence the first statement of the Corollary.

For the second one, let $M'$ be the module affording $\tau$.
If item (1) of Lemma \ref{p33} holds, then $g$ is cyclic on $M$, and hence
 $g$ is cyclic on $M'$. So  $\deg\tau (g)=\dim M' $. A case-by-case inspection
(as
$|g|/|Z(G)|=(q^n+1)/|Z(G)|$ and $\dim\tau$ are known), shows that
$|g|/|Z(G)|\geq \dim\tau$, whence the result.

Next, suppose that item (2) of Lemma \ref{p33} holds.

If $G=GL(n,q)$,  
there is a one-dimensional $FG$-module $L_1$, say, such that
$M\otimes L_1$ is isomorphic to the permutation $G$-module $L$ associated with
the $G$-action
on the vectors of the natural $F_qG$-module $V$ (see comments at the end of Section \ref{sbs61}).
So it suffices to assume that $M=L$. Let  $N$ be
the submodule generated by the zero vector in
$V$, so that $M/N$ is  isomorphic to the permutation $G$-module associated with the $G$-action
on the non-zero vectors of $V$. It is obvious that the matrix of $g$ on the latter module is cyclic,
whence the claim (no matter what are dimensions of the \ir constituents of $M/N$).

Next, we assume that $G$ is unitary or symplectic.
Then $|g|=q^n-1$ and $g$ is almost
cyclic on $M$ (but not cyclic).
By Lemma \ref{ex7}(2),  there exists a one-dimensional
$g$-submodule $N$ of $M$ such that the matrix of $g$ on $M/N$ is  cyclic. It then follows
that
there is at most one $G$-composition factor of  $M$ on which $g$ is not
cyclic. (This is true for an arbitrary $G$-module $M$ admitting  a one-dimensional
$G$-submodule $N$
such that $g$ is cyclic on $M/N$.
For, let $M_1$ be a proper $G$-submodule of $M$. If
$N\subseteq M_1$, then $M/M_1$ is cyclic, and so are all composition
factors of $M/M_1$. As $M_1/N$ is cyclic, the claim follows by
induction on  $\dim M$. If $N\cap M_1=0$ then $M_1$ is isomorphic
to a submodule of  $M/N$, and hence $M_1$ is cyclic. So again the
claim follows by induction.)

Note that the dimensions of the composition factors  belong to the set
$\{1,|g|/|Z(G)|$, $(|g|/|Z(G)|)+1\}$. If  $\dim
\tau=(|g|/|Z(G)|)+1$, then $\tau(g)$ is not cyclic, but by the
above  $\deg\tau(g)=\dim\tau-1$, and
hence  $\deg\tau(g)=|g|/|Z(G)|$, as claimed. If $M$ has a composition
factor of degree $(|g|/|Z(G)|)+1$, then it is unique, and hence all the
other non-trivial factors are cyclic $g$-modules of degree
$|g|/|Z(G)|$. So if $\tau$ affords one of these factors then $
\deg\tau(g)=\dim\tau$, and the second claim of the lemma  follows in this case.

Finally, suppose that $M$ has no composition factor of degree
$|g|/|Z(G)|+1$. This is not the case for  $\ell=0$. So suppose $\ell>0$.
Then
 after realizing the generic Weil \rep by matrices over the
$\ell$-adic field integers, the module $A$, say, afforded by this
representation, has a submodule series $0=A_0\subset A_1\subset
\cdots \subset A_d$, where $d=|Z(G)|$ and the quotients
$A_{i+1}/A_{i}$ correspond to the irreducible Weil
representations. Therefore, $B:=A\pmod \ell$ has a submodule series
$0=B_0\subset B_1\subset \cdots \subset B_d$, where
$B_{i+1}/B_{i}$ is the reduction of $A_{i+1}/A_{i}$ modulo $\ell$.
Note that $A$ is a  Weil module in zero characteristic. So,
as mentioned above,  exactly one factor $A_{i+1}/A_{i}$ has dimension
$(|g|/|Z(G)|)+1$, whereas the others are of dimension $|g|/|Z(G)|$.
This is also true  for the factors $B_{i+1}/B_{i}$. The factors
$A_{i+1}/A_{i}$ of dimension $|g|/|Z(G)|$  remain \ir
modulo $\ell$, and the (single) one of dimension $(|g|/|Z(G)|)+1$
is reducible modulo $\ell$. Denote it by $D$, say. Clearly, the
matrix of $g$ on $D$ is not cyclic. The above reasoning ensures
that the matrix of $g$ on every factor $B_{i+1}/B_{i}$ other than
$D$ is cyclic. So we have to show that the matrix of $g$ on the
non-trivial composition factor $D'$, say, of $D$ is cyclic.
However, it is known that $D'$ lifts to characteristic 0, so $D'$
is isomorphic to a Weil \rep obtained by reduction modulo $\ell$
of an \irr of degree $|g|/|Z(G)|$. So the result follows.

\medskip
 
Our next aim is prove the converse of Corollary \ref{bc44}, by showing that, if $g$ is as in Lemma \ref{p33}, the condition
 $g\in \lan Z(G), h\ran $ is also necessary for $\tau (h)$ to be almost cyclic.
We shall do this below (see Lemmas \ref{sp55} and \ref {cw1}).

 As the \ir constituents of $\phi
(G)$ remain \ir under restriction to $G'$ (provided $n>2$ in the cases
of $GL(n,q)$ and $U(n,q)$, and $(n,q)\neq (2,2), (2,3)$ if
$G=Sp(2n,q)$), then this will imply the corresponding results for
$SL(n,q)$ and $SU(n,q)$.

\medskip
In order to make the proof of the subsequent lemma more transparent we explicitly
state the following:

\begin{lemma}\label{dd5} 
$(1)$
Let $S$ be a finite group and let $Z$ a cyclic central subgroup of $S$. Set $Z=X\times Y$, where $X=\lan x\ran $ is the Sylow $\ell$-subgroup of $Z$.
Then $\rho_S^{reg}=  \oplus _{\lam\in\Irr Y}\lam^S$, and for any fixed $\lam$ the quotient $FS$-modules
$(1-x)^i\lam^S/
(1-x)^{i+1}\lam^S$ for $i<|X|$ are isomorphic to each other and have dimension $|S:Z|$.
In addition, the above quotients can be identified with the induced modules $\overline{\lam}^S$,
where $\overline{\lam}\in \Irr Z$, $\overline{\lam}(X)=1$ and $\overline{\lam}(Y)=\lam$.

\medskip
$(2)$
Let $S$ be cyclic, $h\in S\setminus Z$ and $S_0=\lan Z, h\ran$. Then the matrix of $h$ on
$\overline{\lam}^S$ is almost cyclic \ii $S=S_0$. More precisely, if $f_c$ is the characteristic \po of $h$ on  $\overline{\lam}^S$  and $f_m$ is the minimum polynomial, then $f_c=f^d_m$, where $d=|S:S_0|$.
\end{lemma}

Proof. (1) By elementary properties of induced modules, as $X\cap Y=1$, the module
$\lam^S|_X$ is the direct sum of $|S:Z|$ copies of $\rho^{reg}_X$. The Jordan form of $x$
on $\lam^S$ shows that the quotient module $\lam^S/(1-x)\lam^S$ has dimension $|S:Z|$.
Clearly, this holds for any subsequent factor $(1-x)^{i}\lam^S/(1-x)^{i+1}\lam^S$.
On the other hand, the map which sends $v\in \lam^S$ to $(1-x)v\in (1-x)\lam^S$ induces an epimorphism
of $FS$-modules from $\lam^S/(1-x)\lam^S$ to  $(1-x)^{i}\lam^S/(1-x)^{i+1}\lam^S$. By dimension reasons, this is an isomorphism. The additional claim can be verified directly.

(2) Set $d=|S:S_0|$.  Consider the
restriction of  $\overline{\lam}^S$ to $S_0$. We claim that this restriction
is the direct sum of $d$ copies of
$\overline{\lam} ^{S_0 }$. Set $A=Z(G)$, 
$B=S_0$. By the general theory of induced modules, $(\overline{\lam}
|_A)^{S}|_B$ is the direct sum of $c$ copies of the modules $(\overline{\lam} |_{A\cap B})^B$
(as $S$ is abelian), where $c$ is the number of the double cosets
$A\backslash S/B$. In our case $c=|S|/|AB|=d$. Therefore, the
matrix of $h$ on $\overline{\lam} ^{S } $ is cyclic if $d=1$ (as $\overline{\lam}^S$ is a subquotient of $(1^{Z(G)})^S$), otherwise this
matrix is not even almost cyclic. In particular, this argument also proves the
last claim.

 \begin{lemma}\label{sp55}  Let $G\in \{Sp(2n,q)$, where $n>1$ 
  and $q$ is odd; $U(n,q)$, $n>2$, $(n,q)\neq (3,2)$; $GL(n,q)$, $n>2\}$.
 Let $g$ be as in Lemma $\ref{p33}  (1),(2)$;
 so in particular $g$ is of order $q^n\pm 1$.
Let   $\tau$ be an \ir 
 Weil $F$-representation of $G$, with $\dim \tau >1$. Suppose that $h\notin Z(G)$ and
 $h\in \lan g\ran$.
Then the matrix $\tau (h)$
is  almost cyclic only if $g\in \lan h,Z(G)\ran$. 
\end{lemma}

Proof.  Set $S=\lan g\ran$,  $\ep= \pm 1$ and  $|g|=q^n-\ep$.
Let $M$ be the $FG$-module 
afforded by a generic Weil $F$-\rep of $G$.
 By Lemma \ref{ex7} (in view of the construction of the Weil representations),
if $\ep=-1$ then $M|_S$ is isomorphic to a submodule of codimension 1 in $\rho_S^{reg}$,  the regular  $FS$-module; whereas, if $\ep=1$ then $\rho_S^{reg}$ is a submodule of
codimension 1 in $M|_S$. Observe that $Z(G)\subset S$, and  set $S_0=\lan h,Z(G)\ran$.
We want to prove that, if  $d:=|S:S_0|>1$, then the matrix $\tau (h)$ is not almost cyclic.

 So, assume $d>1$. Note that , for $G=U(n,q)$, this implies $(n,q)\neq (3,2)$. Indeed, if $G=U(3,2)$, then  $\ep=-1$ and $|g| = 9$. Thus $g^3\in Z(G)$, and $d>1$ forces $h\in\langle g^3\rangle =Z(G)$. 

We apply Lemma \ref{dd5}(2), choosing $Z=Z(G)=X\times Y$, 
where $X=\lan x\ran $ is the Sylow $\ell$-subgroup of $Z$, 
and  $Y=\lan y\ran$ (assuming $X=1$ for $\ell=0$). As $d>1$, 
the matrix of $h$ on  $\overline{\lam}^S$ is not almost cyclic.
For each ${\lam\in\Irr Y}$ and for each $i<|X|$, set 
$N_\lam^i=(1-x)^i\lam^{S}/ (1-x)^{i+1}\lam^{S}$.
Thus, $N_\lam^i$ affords the \rep $\overline{\lam}^S$ and
$\dim N_\lam^i=|S|/|Z|$ for every $\lam,i$ (see Lemma \ref{dd5}(1)).

According to Lemma \ref{dd5}(2),  
the action of $h$ on $N_\lam^i$ can be represented by a 
block-diagonal matrix
$\Delta=\diag(D\ld D)$, where the number of the blocks is equal to 
$d=|S:S_0|$, and  each block $D$ is a cyclic matrix of size 
$|S_0:Z|$. This implies that $\Delta$
is never almost cyclic. Moreover, denoting by $R$ the 
underlying space  $N_\lam^i$ of $\Delta$, and assuming that $R$ has a $\Delta$-stable 
subspace $R_1$ of dimension at least $\dim R-2$, we observe that $\Delta|_{R_1}$
is not almost cyclic, unless: either (i) $\dim R_1 = {\dim R} -1$, $d=2$ and $ |S_0:Z|=2$; (ii) or $\dim R_1 = {\dim R} -2$, $d=2$ and $ |S_0:Z| \leq3$. If case (i) holds, then
$g^4\in Z$. As $|g|=q^n\pm 1$ and $|Z|=2,q+1,q-1$
for $G=Sp(2n,q)$, $U(n,q)$, $GL(n,q)$, respectively, we must have that 
$q^n \pm 1$ divides  $8$, $q^n\pm1$ divides $4(q+1)$ and $q^n-1$ divides $6(q-1)$, respectively.
This implies $G= Sp(4,3)$. For this group, the statement follows from Lemma  \ref{423}. If case (ii) holds, then the above applies again if  $ |S_0:Z| =2$.  If  $ |S_0:Z| = 3$,  then  $g^6\in Z$. Arguing as before, since for $G=Sp(2n,q)$, $U(n,q)$, $GL(n,q)$, respectively, we must have that 
$q^n \pm 1$ divides  $12$, $q^n-1$ divides $6(q+1)$ and $q^n-1$ divides $4(q-1)$, respectively. This could only hold for $G=U(3,2)$, which is ruled out by our assumptions.

 Now, for each $\lam\in\Irr Y$ set 
$M_{\lam}=\{m\in M:ym=\lam (y)m\}$. 
Observe that $M=\oplus_{\lam\in\Irr Y} M_{\lam}$, 
where each $M_{\lam}$ is an $FG$-module, as $Y\subseteq Z(G)$. 
Moreover, every $M_{\lam}$ has a filtration
$M_{\lam}\supset (1-x)M_\lam\supset (1-x)^2M_\lam\supset\cdots $,
again because $X\subseteq Z(G)$. Clearly, $\lam^S$ is the 
$\lam(y)$-eigenspace of $y$
on $\rho_S^{reg}$, whereas $M_\lam$ is  the $\lam(y)$-eigenspace of $y$
on $M$. As mentioned above, $M|_S$ is isomorphic
to a submodule of codimension 1 in $\rho_S^{reg}$
if $\ep=-1$, whereas, if $\ep=1$, $\rho_S^{reg}$ is isomorphic a
submodule of codimension 1 in $M|_S$.

Set $M_\lam^i:=(1-x)^iM_\lam/(1-x)^{i+1}M_\lam$.
Then the \f holds: (i) if $\ep=-1$, either $M_\lam^i|_S\cong N_\lam^i$ or
$M_\lam^i|_S$ is isomorphic to a submodule of codimension 1 in 
$N_\lam^i$, so that $\dim M_\lam^i\leq (|S|/|Z|)$; (ii) if $\ep=1$, either $M_\lam^i|_S\cong N_\lam^i$ or
$N_\lam^i$ is isomorphic to a submodule of codimension 1 in 
$M_\lam^i|_S$, so that $\dim M_\lam^i\leq (|S|/|Z|)+1$. This gives 
us information on $\dim M_\lam^i$.

Let $T$ be the $FG$-module afforded by $\tau$. Clearly, 
we may identify $T$ with a composition factor of 
$M_{\lam}^i$ for some $i,\lam$, which we fix for the rest 
of our reasoning.
The core of our argument is to show that either the $FG$-module
 $M_{\lam}^i$ is
irreducible, or $M_{\lam}^i$ contains a composition factor of codimension 1, 
unless $G=GL(n,q)$,
in which case the codimension may be 2. As a consequence, either $T$ is isomorphic to $M_{\lam}^i$, 
or has codimension $1$ in $M_{\lam}^i$, or $G=GL(n,q)$ and $T$ has codimension $2$ in $M_{\lam}^i$. This, in view of the above 
formula $\diag(D\ld D)$ for the matrix of $h$ on $N_{\lam}^i$, will prove that $\tau(h)$ is not almost cyclic, unless possibly when $G=GL(n,q)$ and $M_{\lam}^i$ contains no composition factor of codimension $\leq 1$. In the latter case we shall adjust the matter (see below).
We finally observe that our strategy  depends on the comparison between the dimension 
of $M_{\lam}^i$ and the 
dimensions of the   irreducible constituents of the generic Weil 
\reps of $G$ in cross characteristic,
which have been described at the beginning of this section.

To avoid confusion, we prefer to argue case-by-case.  

(1) Suppose $G=Sp(2n,q)$, $n>1$, $q$ odd.

  First, let $\ep=1$. Then $\dim N_\lam^i=(q^n-1)/2$, and 
 $ (q^n-1)/2 \leq \dim T\leq \dim M_\lam^i\leq \dim N_\lam^i+1=(q^n+1)/2$. Whence the claim.

Let   $\ep=-1$. Then either $N_\lam^i\cong M_\lam^i|_S$  or
$M_\lam^i|_S$ is isomorphic to a submodule of  $N_\lam^i$
 of codimension 1. In the latter case  
 $M_\lam^i$ is irreducible, and $T=M_\lam^i$.
In the former case $(q^n-1)/2 \leq\dim T\leq \dim M_\lam^i\leq 
\dim N_\lam^i=(q^n+1)/2$, and the claim follows
again.

(2)    Suppose $G=U(n,q)$, $n>2$, $n$ even. 
Then 
$\dim T\geq (q^n-1)/(q+1)$ and
 $\dim N_\lam^i=|S/Z|=(q^n-1)/(q+1)$ (as only the case $\ep=1$ occurs).
Thus,  
either $M_\lam^i|_S\cong N_\lam^i$ and $M_\lam^i$  is irreducible, 
or $ M_\lam^i|_S$ contains a submodule of codimension 1 isomorphic 
to $N_\lam^i$. In both cases $\dim T\geq \dim M_\lam^i-1$, and we are done.

(3) Suppose $G=U(n,q)$, where $n\geq 3$ is odd. Here 
$\dim T\geq (q^n-q)/(q+1)$ and
  $\dim N_\lam^i=(q^n+1)/(q+1)$ (as only the case $\ep=-1$ occurs).

If
 $M_\lam^i|_S\cong N_\lam^i$, then either $\dim M_\lam^i =|S/Z|=
(q^n+1)/(q+1)$, and hence $\dim T\geq \dim M_\lam^i -1$.
If $M_\lam^i|_S<N_\lam^i$ then $\dim M_\lam^i\geq 
(q^n-q)/(q+1)$, and hence $T=M_\lam^i$.

(4) Suppose $G=GL(n,q)$, $n>2$. So $\ep=1$  and $|S/Z|=(q^n-1)/(q-1)$.

 Let $V$ be the underlying space for $GL(n,q)$, and let $\Pi$ be the permutation
$FG$-module associated with the natural action of $G$ on $V$. Recall
that $M=\Pi\otimes L$, where $L$ is some one-dimensional $FG$-module.
Therefore, $\tau$ is obtained from a constituent of $\Pi$ by tensoring with
a one-dimensional representation. Such tensoring does not affect almost cyclicity,
so we may assume that $M=\Pi$.
Let $P_0$ be the stabilizer in $G$ of a non-zero vector of $V$, and let  $P$ the 
stabilizer  of the line spanned by this  vector. Then the representation afforded by
$\Pi$ is $1_G\oplus 1_{P_0}^G$. Therefore $\tau$ is a constituent of $1_{P_0}^G$,
as $\dim\tau\neq 1$.
For a moment, denote by $M'$ the submodule of $M $
afforded by $1_{P_0}^G$. As $\Pi=1_G \oplus M'$, we can deal with $M'$ in place of $M$. However, to simplify notation, we
better rename as $M$ the module afforded by $1_{P_0}^G$.
As $S\cap P_0=1 $ and now $\dim M= |S| $, we have $M|_ S\cong \rho_ S^{reg}$.
This implies that $M_\lam^i|_S\cong N_\lam^i$. So $\dim M_\lam^i=\dim N_\lam^i=|S/Z|=
(q^n-1)/(q-1)$.
 
As $P=P_0\cdot Z(G)$, 
every one-dimensional \rep  $\lam$ of $Z(G)$ can be identified with a one-dimensional \rep of $P $ trivial on $P_0$. By the so-called 'Subgroup Theorem' for induced modules, $\lam^S=\lam^G|_S$, as $G=SP$ and $S\cap P=Z(G)$. Furthermore,
by the same theorem,  $\lam^G|_{G'}= \mu^{G'}$, where $\mu=\lam|_{P\cap G'}$.
By \cite[Theorem 9.1.4]{GPPS}, the  dimension of 
any non one-dimensional
\ir constituent of $ \mu^{G'}$ is at least  
    $((q^n-1)/(q-1))-e$, where $e=1$ if $\ell$ does not divide
$(q^n-q)/(q-1)$, and $e=2$ otherwise. Therefore, 
$\dim T\geq \frac{q^n-1}{q-1}-2$. 

It follows that $\dim M_\lam^i\leq \dim T+2$, 
as claimed. 

We conclude that, in all the cases examined,  the matrix of  $\tau(h)$
is not almost cyclic.

\medskip
{\bf Remark}. Recall that every \ir Weil $F$-\rep of   $GL(n,q), n>2,$ (respectively, $U(n,q)$, $n>2$) remains \ir under restriction to $SL(n,q)$ (respectively, $SU(n,q)$). (This follows by degree reasons from the lower bounds known for non-trivial \ir \reps of
$SU(n,q)$ and $SL(n,q)$, using Clifford's theorem.) 
 Therefore,  Lemma \ref{sp55} applies to the case $h\in SL(n,q), n>2,$ (respectively, $h\in SU(n,q)$, $n>2$).
Moreover, this allows us to limit ourselves to consider, with no loss of generality,
the groups $GL(n,q),U(n,q)$ instead of all the groups $G$ such that $SL(n,q)\subseteq G\subseteq  GL(n,q),SU(n,q)\subseteq  G\subseteq  U(n,q)$ in Lemmas \ref{hh3}, \ref{pr1} and \ref{s3s} below.

\medskip

Next, we examine in detail when the condition $g\in \lan h,Z(G)\ran$  holds, under the assumption that $h$ is a $p$-element. We distinguish two cases: (i) the case when $|h|=|g|$; (ii) the case when $|h|<|g|$.

\medskip

The case when $|h|=|g|$ is dealt with by the following:

\begin{lemma}\label{cw0} 

$(1)$ Let $G=Sp(2n,q)$,  where $n >1$ and $q$  is odd, or
  $G=U(n,q)$, where $n>1$ is odd. Let
  $g\in G$ be a Singer cycle for $G$.  Suppose that $h\in \lan g\ran$ and $|h|=|g|$ is a $p$-power.  Then $G=U(3,2)$ and $|g|=9$.

$(2)$ Let $G\in \{GL(n,q), n>2; U(n,q), n>2$;
$Sp(2n,q)$, $n>1$, $q$ odd$\}$ and let
  $g\in G$ be either a Singer cycle for $GL(n,q)$ or a Singer-type cycle for $G \neq GL(n,q)$.  Suppose that $h\in \lan g\ran$ and $|h|=|g|$ is a $p$-power. Then one of the following holds:
 
$(a)$ $G=Sp(4,3)$ and $|g|=8$;

$(b)$ $G=SL(n,2)$,  $n$ is an odd prime and $|g|$  is a Mersenne
prime.
 \end{lemma}

Proof. 
 Suppose that $|g|$ is a prime-power and (1) holds. Then, by
Lemma \ref{zgm}, one of the \f holds: (i) $p=2$ and $n=1$ (which
contradicts our assumptions), (ii) $q$ is even and $q^n+1=p$ is a
Fermat prime;  (iii) $q^n+1=9$, that is, $q^n=8$. If (ii) holds
then $G$ is not symplectic, as $q$ is even. As $q^n+1$ is a prime,
$n$ is even. So $G$ is not a unitary group. Finally, if (iii)
holds, then $|g|=9$, and $q^n=8$, that is, $G=U(3,2)$.

Next, suppose that $|g|$ is a prime-power and (2)  holds. Then, by
Lemma \ref{zgm}, either  $q$ is even and $|g|=q^n-1$ is an odd prime, or $|g|=q^n-1=8$, that is, $q^n=9$. 
In the latter case $G=Sp(4,3)$, in the former case $G=SL(n,2)$ and $n$ is a Mersenne prime.

\medskip

Recall that a prime $p$ is called a Zsigmondy prime for $q^n-1$ if $n$ is the least integer $i>0$ such that $p$ divides $q^i-1$. This can be expressed by saying that $n$ is the order of $p$ modulo $q$.
The classical Zsigmondy's theorem (\cite{Zg}) states that a Zsigmondy prime for $q^n-1$ exists for all
pairs of integers $n,q$ such that $n>2$, $q>1$ and $(n,q)\neq(6,2)$.

\medskip

The case when $|h|<|g|$ is dealt with by the following:

\begin{lemma}\label{cw1}  Let $G\in \{GL(n,q), n>2; U(n,q), n>2, (n,q)\neq (3,2)$;
$Sp(2n,q)$, $n>1$, $q$ odd$\}$ and let 
  $g\in G$ be either a Singer or a Singer-type cycle for $G$.
 Furthermore, suppose that $g\in \lan h,Z(G)\ran$, where 
$|h|$ is a $p$-power, $|h|<|g|$ and $h\in \lan g\ran$.
Then
$p>2, $  $(p,|Z(G)|)=1$, $\lan h\ran$ is a Sylow $p$-subgroup of $G$,
and one of the \f holds:

\medskip
$(1)$ $G= GL(n,q)$, $|h|=(q^n-1)/(q-1)$ and $n\neq p$ is an odd prime;

\medskip
$(2)$ $G= Sp(2n,q)$,  $|h|=(q^n+1)/2$ and $n$ is a $2$-power;

\medskip
$(3)$ $G=Sp(2n,3)$, $|h|=(3^n-1)/2$ and $n\neq p$ is an odd prime;

\medskip
$(4)$ $G= U(n,q)$, $|h|=(q^n+1)/(q+1)$ and $n\neq p$ is an odd prime; 

\medskip
$(5)$ $G= U(4,2)$  and $|h|=5$.
 \end{lemma}

Proof.  Obviously,   $h\notin Z(G)$. Under our assumptions, 
$|h|=p^k$ for some
integer $k>0$. Furthermore $|Z(G)> 1$, as $|h|<|g|$.

We first show that $(p,|Z(G)|)=1$. Suppose the contrary. Assume first that $G=Sp(2n,q)$. Then $p=2$ and $q^n \pm1$ is a $2$-power, which implies $q=3$ and $n=2$, that is,
$G=Sp(4,3)$ and $|g|=8$. But in this case   $\lan h\ran $ contains $Z(G)$, and hence
$  g\notin \lan h,Z(G)\ran$,  against our assumption. If $G=GL(n,q)$, $(n,q) \neq (6,2)$, we can apply Zsigmondy's theorem
to find a prime $s\neq p$ that divides $|g|$ and
does not divide $|Z(G)|$ (recall that $n>2$ in this case). This contradicts the assumption
$g\in \lan h,Z(G)\ran$. The case $(n,q) = (6,2)$ is trivial as $Z(G)=1$.
If $G=U(n,q)$, $n$ even, then $|g|=q^n-1$, and again by Zsigmondy's theorem, there is
a prime $t$, say, dividing $q^n-1$ but not $q^2-1$ (unless $(q,n)=(2,6)$,
but in this case   $g\notin \lan h,Z(G)$).
As $t\neq p$, we get a contradiction. Finally, if $G=U(n,q)$, $n$ odd, then $|g|=q^n+1$. By Zsigmondy's theorem, there is a prime $u$, say, dividing $q^{2n}-1$, but neither $q^n-1$ nor $q^2-1$.  Here again we reach a contradiction.

Thus, $p$ is coprime to $|Z(G)|$. This implies that $|g|=|h|\cdot |Z(G)|$ (by our assumption) and that $p>2$. The latter claim is obvious if $q$ is even, 
otherwise it  follows from the fact that  $|Z(G)|\in\{2,q\pm 1\}$.

Now, suppose that  $|g|=q^n-1$. First, observe that $(n,q) \neq
(6,2)$. (Otherwise $|Z(G)|=3$. But this implies that $|h|$ is not
a prime-power, against our assumptions.) Also, $p$ is the only
Zsigmondy prime for $q^n-1$. For, suppose that $t\neq p$, say, is
another Zsigmondy prime for $q^n-1$. Then, as $|g|=|h|\cdot
|Z(G)|$, $t$ divides $|Z(G)|$, whence $G$ is unitary and  $t |
(q^2-1)$. This in turn implies $n=2$, which is a contradiction as
for unitary groups we assume $n>2$.

Next, we claim that either (5) holds, or $n$ is an odd prime
different from $p$. For, suppose that $n= \nu t$ for some integers
$t,\nu$, where $1< \nu <n$. Then $p$ does not divide $q^ \nu -1$,
which in turn implies that $q^ \nu -1$ divides $|Z(G)|$. This
occurs if and only if (5) holds. So, assuming (5) does not holds,
$n$ is prime. Furthermore, $n$ must be odd. For, $n=2$ implies
that $G=Sp(4,q)$ and $4|(q^2-1)$. But then $q^2-1=|g|=|h|\cdot
|Z(G)|= p^k \cdot 2$, a contradiction as $p>2$. So $n$ is odd.
Now, suppose that $n=p$. Then $p$ is a Zsigmondy prime for
$q^p-1$. However, since the Galois group of $F_{q^p}$ over $F_q$
is  of order $p$, all the Galois group orbits on $F^{q^p}\setminus
F_q$ are of size $p$. So  $q^p-q=q(q^{p-1}-1)$ is divisible by
$p$. As $p$ is coprime to $q$, it follows that $p$ divides
$q^{p-1}-1$, a contradiction.

Additionally, we observe that 
if $G= Sp(2n,q)$ then $q=3$. Indeed,  we have $|g|=q^n-1=|Z(G)|\cdot |h|$ and $|h|<|g|$, so
  $|g|=2\cdot |h|=p^k\cdot 2$. As $p$ is a Zsigmondy prime for $q^n-1$,
 $(p,q-1)=1$. This forces $q-1=2,$ so (3) holds.

In conclusion: if $|g|=q^n-1$, then one of the cases (1), (3), (5) holds.

Next, let us consider the cases where $|g|=q^n+1$.  First, observe that  $p$ is the only Zsigmondy prime for $q^{2n}-1$. For, suppose that $t\neq p$, say, is another Zsigmondy prime for $q^{2n}-1$. Then, as $|g|=|h|\cdot |Z(G)|$, $t$ divides $|Z(G)|$, whence $t | (q^2-1)$, which is impossible, as $n>1$. (Notice that $(2n,q) \neq (6,2)$, since otherwise $G =U(3,2)$, which is excluded by our assumptions).

Suppose first that $G= U(n,q)$. Then $n>2$ is odd,  and $|g|=q^n+1=|h| \cdot (q+1)$, where $|h|=p^k$ for some $k>0$. We claim that $n$ is a prime different from $p$. For, suppose that $n= \nu s$, where $1< \nu <n$. By the above, $p$ is the unique Zsigmondy prime for $q^{2n}-1$. On the other hand, $q^n+1=(q^ \nu)^s+1=(q^ \nu +1)c=p^k(q+1)$, for some integer $c$. As $\nu>1$ is odd, this implies that $p$ must divide $q^\nu+1$, a contradiction. 
Assume that $n=p$. As above, by elementary Galois theory we obtain  $p$ divides $q^{2p}-q^2=q^2(q^{2p-2}-1)$. As $p$ is coprime to $q$ and does not divide $q^{2p-2}-1$, we get a contradiction. So we have case (4) of the statement.

Finally, suppose that $G= Sp(2n,q)$ and $|h|=(q^n+1)/2$. Then it is easily seen that $n$ must be a $2$-power. Indeed, suppose the contrary. Let $n=s\cdot d$, say, where $s$ is an odd prime. Then $q^n+1=(q^d)^s+1=(q^d+1)(q^{d(s-1)} -q^{d(s-2)}+ \cdots + 1)$, where both the factors in the last expression are greater than $2$. It follows that $p$ must divide $q^d+1$, and hence $q^{2d}-1$. A contradiction, as $p$ is a Zsigmondy prime for $q^{2n}-1$. So $n$ is a $2$-power, and we get case (2) of the statement.

 We are left to show that $\lan h \ran$ is a \syl of $G$. To this end, recall that $p$ is a Zsigmondy prime for $q^n-1$ if this is the order of $g$, and for $q^{2n}-1$ if $|g|=q^n+1$. Then the well-known formulas for the orders of classical groups (e.g. see \cite{KL}, p.19) show that $|G|_p=|q^n-1|_p$ in the first case, and  $|G|_p=|q^n+1|_p$ in the second case. It follows that, for each group $G$ under exam, the subgroup $\lan g \ran$ contains a Sylow $p$-subgroup of $G$ (which is therefore cyclic). Furthermore, as  $\lan g \ran= \lan h \ran \times Z(G)$ and $(p,|Z(G)|)=1$, we have $|h|_p=|g|_p$, and hence the Sylow $p$-subgroup of $G$ contained in  $\lan g \ran$
coincides with  $\lan h \ran$.
\smallskip

{\bf Remark}. The previous Lemma is clearly false when $G=U(3,2)$. This solvable group can be fully handled by direct computation, looking at the character table and the Brauer character tables of $G$. Note that $|G| = 2^3.3^4$, and we only need to examine the behaviour of non-scalar elements of order $3$ and elements of order $9$ (these are Singer cycles of $G$). Let $\tau$ be any irreducible $F$-representation of $G$. The following holds:
\smallskip

i)  Let $\ell=0$.    Then almost cyclicity for $g$ semisimple of prime-power order occurs if and only if:
$|g|=3$, $g$ belongs to any non-scalar class (in the GAP labelling: classes 3c,d,e,f,g,i), and $\dim \tau = 2, 3$ (that is, $\tau$ is Weil);
$|g|=9$, $g$ belongs to the classes 9a,9b (GAP labelling), and again $\dim \tau = 2, 3$. Here $g$ is in fact cyclic.

ii)  Let $\ell = 2$. The 2-modular irreducibile representations of $G$ have degrees 1,3,8. $\tau(g)$ is almost cyclic if and only if  $|g|=3$ or $9$ and $\dim \tau=3$. If $|g|=9$, $\tau(g)$ is cyclic.

iii) Let $\ell = 3$. There are just two 3-modular non-trivial  irreducibile representations, of degrees $2$ and $3$, namely the Weil representations. In both cases all the elements of $G$ of $3$-power order are obviously represented by almost cyclic matrices.

As for $G=SU(3,2)$, a group of order $2^3.3^3$ with Sylow $3$-subgroups of exponent $3$, the following holds. Let $g \in G$ be a non-scalar element of order $3$. Then:

i) if $\ell=0$, then $\tau(g)$ is almost cyclic if and only if $\dim \tau = 2, 3$. 

ii) if $\ell = 2$, then $G$ has no irreducible $2$-modular representations of degree $2$, and $\tau(g)$ is almost cyclic if and only if $\dim \tau = 3$.

iii) if $\ell = 3$, then $\tau(g)$ is almost cyclic if and only if $\dim \tau = 2, 3$.  ~

\subsection{Orthogonally decomposable elements}

In this subsection we  deal with the case when $g \in G$ is orthogonally decomposable.
We begin with an auxiliary Lemma:

\begin{lemma}\label{hh3} Let
$G=U(n,q)$, where $n>2$ is even and $(n,q)\neq (4,2)$, and let $g\in G$ be an element of $p$-power
order for some prime $p$, stabilizing a subspace $W$ of $V$ of dimension $n-1$ and
acting on $W$ irreducibly (so $(p,q)=1$). Let $\tau$ be an \ir Weil \rep of $G$.
Then $\tau (g)$ is not almost cyclic. \end{lemma}

Proof. Observe that $W^\perp$ is $g$-stable, and hence (as $n>2$)
$V=W \oplus W^\perp$ (so $W$ is non-degenerate). Thus, $g$ belongs
to a subgroup $H$ which can be identified with $U(W)\times
U(W^\perp)$ (where the latter group is cyclic of order $q+1$), and
hence $g$ is orthogonally decomposable. Let $g=g_1g_2$, where
$g_1\in U(W)$, $g_2\in U(W^\perp)$. Clearly both $g_1$ and $g_2$
are of $p$-power order.  Let $\tau_0$ be an \ir constituent of
$\tau |_H$  of dimension greater than 1. Then $\tau_0 (g)=\tau
_1(g_1)\otimes \tau_2(g_2)$, where $\tau_1$ is an
irreducible Weil representation of $U(W)$ of dimension greater than 1,
and $\tau_2$ is a $1$-dimensional  representation of $U(W^\perp)$ (see for instance
\cite[Lemma 4.2]{TZ2}).

By way of contradiction, suppose that $\tau (g)$ is almost cyclic.
Then $\tau _1(g_1)$ is almost cyclic (as $\tau_2$ is
1-dimensional). Since $g_1$ acts irreducibly on $W$, $g_1$ belongs
to a Singer subgroup of $U(W)$. By Lemmas \ref{p33} and
\ref{sp55}, $\lan g_1,Z(U(W))\ran $ is of order $q^{n-1}+1$.
 As $(n,q)\neq (4,2)$, the option $(n-1,q)=(3,2)$ recorded  in Lemma \ref{cw0},(1) is ruled out,
 and
therefore we may apply Lemma \ref{cw1} (where $G=U(n-1,q)$ and
$h=g_1$).
 We find that case (4) of Lemma  \ref{cw1}
must hold, and hence $|g_1|=\frac{q^{n-1}+1}{q+1}$, where $p\neq
n-1 $ is an odd prime. In addition, $P:=\lan g_1\ran$ is a \syl of
$U(W)$, $p$ is coprime to $q+1$. \itf $g_2=1$,
and hence $g=g_1$. Furthermore, $(p,q^{n}- 1)=1$.
 (Indeed, as $p$ divides
$q^{n-1}+1$, it does not divide $q^n-1=q(q^{n-1}+1)-(q+1)$, as  $(p,q+1)=1$.)
So, in fact $P$ is a \syl of $G$. Indeed, $|G|=q^a(q^n-1)\cdot |U(W)|$ for some natural number $a$,
and hence $p$ does not divide the index $|G:U(W)|$.

Recall that $\dim \tau\in \{\frac{q^{n}-1}{q+1},
\frac{q(q^{n-1}+1)}{q+1} \}$. Hence,
 $\dim \tau >|g|+1$, as $|g|=\frac{q^{n-1}+1}{q+1}$.
 If $\ell \neq p$, then we may assume  $\ell=0$, as the \ir
Weil \reps of $G$ lift to characteristic zero. Thus, we only need to consider the cases
$\ell =0$ and $\ell=p$. If $\dim \tau $ is divisible by $|g|$, then $\tau$ is of $p$-defect $0$, and hence $\tau (\lan g\ran)$ is a direct sum of regular  $F\lan g\ran$-modules.
As we are assuming that $\tau (g)$ is almost cyclic, this implies that
$\dim \tau =|g|$. But this is not the case.

So, suppose that $\dim \tau $ is
not divisible by $|g|$. Then $\dim\tau=\frac{q^{n}-1}{q+1}$.
It follows that  $\dim \tau +1=q\cdot |g|$.

 If $ \ell=0$, then
by \cite[Lemma 7.4]{RZ} (here $\tau$ is one of the $x'_i$ in \cite[Lemma 7.4]{RZ}),
 $\tau (\lan g\ran)$ contains
$q-1$ regular $F\lan g\ran$-modules, plus the quotient of the regular $F\lan g\ran$-module by a
one-dimensional submodule.  This gives a contradiction.

Next, suppose $\ell=p$.
In order to use Lemma \ref{ww1}, we show that the group
$N_G(P)/P$ is abelian.

Set $N:=N_G(P)$, and let $C_G(P)=P\cdot C$, where $C$ is a complement of $P$.
It is easy to see that $C$ is abelian (indeed, $C= Z(U(W))\times U(W^{\perp})$).
As $W^{\perp}$ is obviously the fixed-point subspace of $P$ on $V$, it follows that $W^{\perp}$,
and hence also
$W$, are $N$-stable. Thus,
$N\subseteq H= U(W)\times U(W^\perp)$. Then, obviously, $[N,C]=1$ and hence $C\subset Z(N)$.
Let $ T$ be a complement of $P$ in $N$. Then $C\subseteq T$ and $T$ acts on $P$ with kernel
$C$. Since $P$ is a cyclic $p$-group, where $p>2$,  
Aut$\, P$ is cyclic.
\itf $T$ is abelian, as $[T,C]=1$, and so is $N_G(P)/P$ (being a cyclic extension of a central subgroup).

As $P$ is a Sylow $p$-subgroup of $G$ and $p=\ell$,
by Lemma \ref{ww1} the restriction to $P$ of
the $FG$-module $M$ associated to $\tau$ decomposes (using the
notation of Lemma \ref{ww1}) as $M|_P=\frac{\dim M -\dim L}{|P|}
\rho^{reg}_P\oplus L$, where $L|_P$ is indecomposable and $\dim L<|P|$
(since $N_{G}(P)/P$ is abelian). By the above,  $\dim M=\dim \tau=q\cdot
 |P|-1$. This implies  $\dim L\equiv -1$ $( \mod |P|)$, whence $\dim L >1$.
As $\tau(g)$ is almost cyclic,
 this in turn forces $M|_P=L|_P$, which is not the case.

\medskip
The \f lemmas will be used for induction purposes:

\begin{lemma}\label{pr1} 
Let $G=U(n,q)$, $n>2$, $(n,q)\neq (3,2)$, and let $g\in G$ be a non-scalar
 semisimple element of prime-power order dividing $2(q\pm 1)$.
 Let $\tau$ be an \ir Weil representation of $G$, with $\dim \tau>1$. Then
  $\tau(g)$ is  almost cyclic \ii one of the \f holds:

 \medskip
 $(1)$ $G=U(3,3), $ $|g|=8$, and either
 $\dim\tau=6$, or   $\ell\neq 2$ and $\dim\tau=7;$

 \medskip
 $(2)$ $G=U(4,2)$, $|g|=3$, and either $\dim\tau=5$, or $\ell\neq 3$ and $\dim\tau=6$.

 \medskip
 In addition, $\tau(g)$ is cyclic \ii $G=U(3,3)$ and $\dim\tau=6$.

\end{lemma}

Proof. Let $G_1=\lan SU(n,q),g \ran$. Clearly $G_1$ is a normal subgroup of $G$.
Let $\tau_1$ be an \ir constituent of $\tau |_{G_1}$. Then, by Clifford's theorem,
$\dim\tau_1>1$. Indeed, otherwise, $SU(n,q)$ would lie in $\ker \tau_1$, and hence
in $\ker \tau$. As $U(n,q)/SU(n,q)$ is abelian, this would imply $\dim\tau=1$,
which is not the case.

So, $\dim\tau_1>1$. Suppose that $\tau(g)$ (and hence $\tau_1(g)$) is almost
 cyclic. Let $M_1$ be the module afforded by $\tau_1$ and suppose that
 neither $(n,q)=(3,3)$, nor $n=4$ and $|g|$ is a $2$-power. Then, by Propositions
 \ref{GS1} and \ref{GS2}, $n$ suitable $SU(n,q)$-conjugates of $g$ suffice to generate
 $G_1$. Thus, by Lemma \ref{nd5},
 $\dim M_1\leq n\cdot (|g|-1)\leq n\cdot (2q+1)$. As   $\dim M_1\geq
 (q^n-q)/(q+1)$ (see \cite{LS}; the exceptions for $SU(4,2),SU(4,3)$ recorded in \cite{LS}
occur only for projective representations),
it follows that $n\cdot (2q+1) \geq
 (q^n-q)/(q+1)$, or equivalently $n(q+1)(2q+1)\geq q(q^{n-1}-1)$. But this only
 holds if either  $n=3$ and $q\leq 7$, or $n=4$ and $q\leq 3$, or $n=5,6$ and $q=2$.

 Direct computations using the GAP package show that, if $n=3$ and $4 \leq q
 \leq7$, $\tau(g)$ is not almost cyclic, whereas if
 $(n,q)=(3,3)$ the exceptional case listed in $(1)$ arises.

 If $n=4$ and $|g|$ is not a $2$-power, then, by the above, only the case
 $G=U(4,2)$ needs to be examined. This case is dealt with by Lemma \ref{442}, yielding the exceptional item listed in (2).

 So, suppose that $n=4$ and $|g|$ is a $2$-power. Then, by Proposition
 \ref{GS2}(2),  as $g$
is semisimple, four suitable $SU(4,q)$-conjugates of $g$ suffice to generate
 $G_1$. In this case, the condition $(q^4-q)/(q+1) \leq \dim M_1\leq 4\cdot
 (2q+1)$ must hold, and this only happens if $q \leq 3$. Computations using GAP rule out the case
 $q=3$.
The case $q=2$ may be ignored, as $g$ is semisimple.

Finally, suppose that $n=5,6$ and $q=2$. The case $n=6$ is easily ruled out, since on one hand $\dim M_1\leq n(|g|-1)=12$, but 
on the other hand $(q^n-q/q+1)\geq 20$. So, assume that $G=U(5,2)$ and let $V$ be the natural module for $G$.  Since $|g|=3$, $g$ is diagonalizable on $V$, and hence has an eigenspace of dimension at least $2$ on $V$, by dimension reasons. Therefore, $g$ stabilizes an isotropic 1-dimensional subspace, say, $W$. Then $g\in P$, where $P$ is the stabilizer of $W$ in $G$. Let $U$ be the unipotent radical of $P$. We know that $g$ acts faithfully by conjugation on $U/Z(U)$. Set $X=\langle g,U\rangle$, and let $\phi$ be
an irreducible constituent of $\tau|_X$, non-trivial on $Z(U)$. Set $E=\phi(U)$. Then $E$ is a group of symplectic type, and $E/Z(E) \cong U/Z(U)$ has order $2^6$ (see \cite{DMZ0}). However, as $\phi(g)$ is assumed to be almost cyclic, Lemma \ref{a55} implies that $|g|=2^3\pm 1$, a contradiction, as $|g|=3$.

\begin{lemma}\label{rd6}  Let $G\in\{ U(n,q)$,
$n>2$, $Sp(2n,q)$, $n>1$ and $q$ odd\},
 and let $H=G_1\times G_2$, where $H$ is the stabilizer in $G$  of
a non-degenerate $m$-dimensional subspace
$W$ of the natural module for $G$ (so
$G_1=C_H(W^\perp),G_2=C_H(W)$). Suppose that $G_1$ is non-solvable.
Let $\tau$ be an \ir Weil $F$-\rep of $G$. Then either $\tau|_H$ contains
an \ir  constituent
$\phi$ such that 
$\phi=\phi_1\otimes \phi_2$, where $\phi_1, \phi_2$ are \ir Weil $F$-\reps of $G_1,G_2$,
 both of dimension greater than $1$, or  one of the \f holds:

\med
$(1)$   $G=U(n,q)$ and  $n-m=1;$

\med
$(2)$ $G= Sp(2n,3)$,  $G_2\cong Sp(2,3)$ and $\ell=2;$ in this case
 the restriction of $\tau$ to the derived subgroup of $H$ contains at least $2$
 isomorphic composition factors of dimension greater than $1;$

\med
$(3)$ $G=U(n,2)$, $\ell=3$ and $G_2\cong U(2,2);$ in this case
the restriction of $\tau$ to the derived subgroup of $H$  contains at least $2$
 isomorphic composition factors  of dimension greater than $1$.
\end{lemma}

Proof. Suppose that (1) does not hold.

Case (i): $G_2$ is non-solvable.

Suppose first that $\ell=0$. Let $\om$ be a generic Weil \rep   of $G$.
Then   the restriction of $\om$ to $G_i$
($i=1,2$) is the sum of  generic Weil \reps of $G_i$ by Lemma \ref{wr3}. \itf $\om|_H$ is the sum of
\ir \reps of shape $\phi_1\otimes \phi_2$, where $\phi_1, \phi_2$ are \ir Weil \reps of $G_1,G_2$,
respectively.
Therefore, this is also true for $\tau|_H$. As we assume that $G_1,G_2$ are both non-solvable,
neither of them has an \ir Weil \rep of dimension 1. So we are done in the case $\ell=0$.

Now, suppose that $\ell>0$. Recall that $\tau$ lifts to characteristic 0. Let
$\overline{\tau}$ be the lift of $\tau$, and let $\overline{\phi}$ be a
composition factor of $\overline{\tau}|_H$. Then every composition
factor of $\overline{\phi}\pmod \ell$ is a composition factor of
$\tau|_H$. Let $ \overline{\phi}=\overline{\phi}_1\otimes
\overline{\phi}_2$, where $\overline{\phi}_1,\overline{\phi}_2$
are \ir \reps of $G_1,G_2$, respectively (both of dimension greater than 1).
Then $ \overline{\phi}({\rm mod}\,
\ell)$ contains all the composition factors of $\overline{\phi}_1({\rm mod}\,
\ell)\otimes \overline{\phi}_2({\rm mod}\, \ell)$. Clearly, each
$\overline{\phi}_i({\rm mod}\, \ell)$ ($i=1,2$) 
contains a composition factor of dimension greater than 1, otherwise
$\overline{\phi}_i({\rm mod}\, \ell)$ would be solvable. Setting $\phi=
\overline{\phi}\pmod \ell$, the
result follows.

\medskip
Case (ii) $G_2$ is solvable.

Then either $G_2\cong Sp(2,3)$ or $G_2\cong U(2,2), U(2,3), U(3,2)$. In each case
$G_2$ is non-abelian.
Again, let us start with the case $\ell=0$. As $\tau$ is faithful on $G_2$,  $\tau|_{G_2}$
has an \ir constituent $\phi_2$, say, of dimension at least 2. \itf
$\tau|_{H}$
has an \ir constituent $\phi=\phi_1\otimes \phi_2$, where $\phi_1$ is an \ir Weil \rep of $G_1$, as required.
Now, let  $\ell>0$, and let $\phi$ be chosen as above.
If $\ell$ is coprime to the order of $G_2$, that is, $\ell\notin\{2,3\}$,
  $\phi({\rm mod}\,
\ell)$ behaves like $\phi$, and we are done.
So we may assume that
 $\ell=2$ if $q=3$, and  $\ell=3$ if $q=2$. If $G_2\cong Sp(2,3)$, then $G_2/O_2(G_2)$
is of order 3. So the reduction modulo 2 of the 2-dimensional Weil \rep of $G_2$
is a completely reducible non-trivial \rep of dimension 2. \itf $\phi ({\rm mod}\,2)$
is the direct sum of two \reps of $H$, which are isomorphic under restriction
to $H'=G_1\times O_2(G_2)$.
This gives us
case $(2)$ of the statement. If $G_2\cong U(2,2)$, then $G_2/O_3(G_2)$ has order 2,
so the reduction mod 3  of the 2-dimensional Weil \rep of $G_2$
is a completely reducible non-trivial \rep of dimension 2.
\itf $\phi({\rm mod}\,3)$
is the direct sum of two \reps of $H$, which are isomorphic under restriction
to $H'=G_1\times O_3(G_2)$.
Finally, in both  $U(2,3)$ and $U(3,2)$, the reduction mod $\ell$ of an ordinary \ir Weil \rep
contains a composition factor
of dimension greater than 1. This gives case (3) of the statement.

\begin{lemma}\label{s3s}  Suppose that $G=U(n,q)$,
 where $n>2$ and $(n,q)\neq (5,2),(4,2),(3,3),(3,2)$. Let $g\in G$ be
 a non-scalar semisimple element of
$p$-power order for some prime $p$, and let $\tau$ be an
irreducible Weil $F$-representation of $G$. Suppose that $\tau (g)$ is
almost cyclic. Then $n\neq p$ is an odd prime, $g$ is irreducible of order  $(q^n+1)/(q+1)$,
and $\lan g \ran$ is a Sylow $p$-subgroup of $G$.\end{lemma}~

Proof.
Suppose  that $g\in G$ is orthogonally indecomposable. Then
(taking into account  Lemma \ref{cw0}) our $g$ satisfies the
assumptions on  $h$ in Lemma  \ref{cw1}, and therefore the result follows from that lemma.

So, assume  that $g$ is orthogonally decomposable. We aim to show that this cannot occur. 

 Let $V$ be the natural module for $G$, let  $W$ be a non-degenerate $g$-stable subspace of $V$ such that  $g|_W$ is orthogonally indecomposable, and choose $W$ such that $|g| = |g|_W|$.  Set $m=\dim W$.
 By Lemma \ref{pr1}, we can assume that $m>2$. (Indeed, otherwise, $|g|$ divides $q^2-1$; as $|g|$
is a prime power, $|g|$ divides either $2(q+1)$ or $2(q-1)$. This is impossible by Lemma \ref{pr1},
as we exclude the cases $(n,q)=(4,2),(3,3)$.)

So $m>2$, and $g$ belongs to a
subgroup $H=G_1\times G_2$,  where $G_1\cong U(m,q)$ and $G_2\cong
U(n-m,q)$.
Let $g=g_1g_2$, where $g_1\in G_1$, $g_2\in G_2$. Then $|g|=|g_1|$; moreover, either $m$ is 
odd and $|g|$ divides $q^{m}+1$, or $m$ is even and  $|g|$
divides $q^{m}-1$.

We first rule out the second possibility. Indeed, let $|g|$ divide $q^{m}-1$ and set
$G_3=\lan G_1,g\ran=\lan G_1,g_2\ran$. Then $G_3= G_1\cdot Z(G_3)$. Let
$\phi$ be an \ir constituent of  $\tau|_{ G_3}$ of dimension greater than 1.
By Schur's lemma, $\phi(g)$ is a
scalar multiple of $\phi(g_1)$. As $\phi|_{ G_1}$ is a Weil \rep
of $ G_1$ by Lemma \ref{wr3},    then, by Lemmas
\ref{cw0} and \ref{cw1},  $G_1\cong U(4,2)$ and $|g_1|=5$. In this case
$\dim\phi=5$ or 6.  By Corollary \ref{bc44},
$\deg\phi(g_1)=5$. If $g_2=\Id$, then $G_3=G_1$ is contained in a subgroup
$H_1=SU(5,2)$ (the pointwise stabilizer of a non-degenerate subspace
of $W^\perp$ of dimension $n-5$). So, we may assume that $g\in H_1$.
Suppose first that $\ell = 5$. As a Sylow $5$-subgroup $S$  of $H_1$
is cyclic and $N_{H_1}(S)$ is abelian, it follows from Lemma \ref{ww1} and Corollary \ref{ddw} that every \ir constituent of $\tau|_{H_1}$
is of degree at most 6. However, inspection of the Brauer character tables for
all $\ell\neq 2$ in  \cite{MAtl} shows that the minimum dimension of a
non-trivial \ir $F$-representation of
$H_1$ is 10, a contradiction. Next, suppose that $\ell \neq 5$. By Proposition \ref{GS1}, $H_1$ can be generated by five conjugates of $g_1$, and hence, by Lemma \ref{nd5}, we only need to examine $F$-representations of $SU(5,2)$ of degree at most 20. However, inspection of the character table and the Brauer character tables of  $H_1$ shows that such representations can only have degree 10 or 11, with character or Brauer character having respectively value equal to 0 or 1 on elements of order 5. This obviously implies that $\tau (g)$ is not almost cyclic, a contradiction.
So, suppose that $g_2\neq\Id$, and hence $|g_2| = 5$. Then $\dim W^\perp\geq 4$ (as $U(n-m,2)$ does not have elements of order $5$ for $n-m \leq3$); in
particular, both $G_1$ and $G_2$ are not solvable. By Lemma \ref{rd6}, we
can choose $\tau$ so that $\tau|_{H}$ contains an irreducible constituent  of shape $\tau_1\otimes \tau_2$,
where $\dim\tau_i>1$. Then, by Lemmas \ref{tp1} and \ref{p99}, the
matrix $\tau_1(g_1)\otimes \tau_2(g_2)$ is not almost cyclic, again a contradiction.

Thus, $|g|$ divides $q^{m}+1$. Observe that $n-m>1$. Otherwise $m=n-1$,
but this option
is ruled out by Lemma \ref{hh3} (note that $g$ acts irreducibly on $W$). In fact, this lemma also rules out the case
$n-m=2$. Indeed,  if $g_2\in U(2,q)$, then $|g_2|$ divides $2(q\pm 1)$
and $q^m+1$. As $m$ is odd, $\frac{q^m+1}{q+1}$ is odd,
so $(q^m+1, 2(q+1))=q+1$, and either $(q^m+1, 2(q-1))=
2$, or $4 | (q+1)$ and hence $(q^m+1, 2(q-1))=
4$. In both cases $|g_2|$ divides $q+1$. \itf $g_2$ stabilizes a non-degenerate subspace of dimension 1 on
$W^\perp$. In this case $g$ belongs to a subgroup $X$, say,
isomorphic to $X_1\times X_2$, where $X_1=U(n-1,q)$, $X_2= U(1,q)$,
and the restriction $\tau|_X$
contains an \ir constituent $\phi$, say, non-trivial on the commutator subgroup
 of $X_1$. Thus,
we may apply Lemma \ref{hh3} to $\phi$, getting that $n-m>2$, unless $G=U(5,2)$. But this case is ruled out by our assumptions.

Now, suppose first that $G_1$ is not
solvable. Then, by the above, $2<m<n-2$. Again by Lemma \ref{rd6}, the restriction of $\tau$ to
$G_1\times G_2$ contains a composition factor
 $\lambda$ of shape $\lambda_1\otimes \lambda_2$, where $\lambda_i$ is an
\ir Weil \rep of $G_i$ and $\dim \lambda_i>1$ for $i=1,2$. 

Let $\ell \neq p$. Then the
matrix $\lambda_1(g_1)$ is  cyclic (see Lemma \ref{p74}), and hence, by Lemmas
 \ref{cw0} and   \ref{cw1} applied to $\lambda_1(g_1)$, 
  case (4) of Lemma \ref{cw1} holds
 for $\lambda_1(g_1)$. In particular, $p$ is coprime to $q+1$, whence, by
Corollary \ref{bc44}, $\deg\lambda_1(g_1)\geq |g_1|-1=|g|-1$.

 Clearly, neither  $\lambda _1(g_1)$ nor $\lambda _2(g_2)$ are scalar, otherwise  the matrix of $\lambda
(g)$ is not almost cyclic. In particular, 
$\deg\lambda_2(g_2)\geq 2$.
 Set $k=\deg\lambda_1(g_1)$, $l=\deg\lambda_2(g_2)$. 
Then $k\geq l$. (Indeed, otherwise $l=|g|$, as $k\geq |g|-1$. By  Lemma  \ref{p74}, $|g|\geq kl-\min\{k,l\}+1\geq ( |g|-1)|g|-( |g|-1)+1= |g|^2-2 |g|+2$, whence $|g|\leq 2$, a contradiction.) So $\min\{k,l\} = l$, and hence, by Lemma  \ref{p74},
$|g|\geq kl-l+1=(k-1)l+1\geq 2( |g|-2) +1$ (as $l\geq 2$), whence $|g|=3$. Since $\lambda_1(g_1)$ is cyclic, it follows that  $\dim\lam_1\leq 3$, which is not the case, as $G_1$ is not solvable.

So, let $\ell =p$.  By Lemma \ref{p99}, $\deg\lambda_1(g_1)\leq 2$. Since $\deg\lambda_1(g_1)\geq |g|-1$, this implies $|g| \leq 3$. As above, this leads to a contradiction with the assumption  $2<m<n-2$.

We are left with the case where $G_1$  is solvable.  Thus $q=2$, and 
$m=3$.  Suppose first that  $G_2$ is  also solvable. Then $n-m=3$, and hence $G=U(6,2)$ and $|g|=9$, by Lemma \ref{pr1}. But this case is ruled out by Lemma \ref{u62}. 
Next, suppose that $G_2$ is not solvable. Then $n \geq 7$. Furthermore, by Lemma \ref{pr1}, we may assume that $|g| = 9$. Let $W'$ be a minimal non-zero  non-degenerate $g$-stable subspace of $V$ lying in $W^\perp$. It is easy to see that either $ \dim{W'}=1$ or $ \dim{W'} = 3$. Indeed, suppose that $ \dim{W'}=2$. Then $g|_{W'}$ has order 3 and acts irreducibly on $W'$; hence its minimum polynomial has degree 3 and splits over $F_4$. But this implies that $g$ is reducible on $W'$, a contradiction. Next, suppose that $ \dim{W'}>2$. Observe that $g^3$ is diagonalizable over $F_4$, and hence acts scalarly on $W'$. This implies that $ \dim{W'}<4$, and we conclude that $ \dim{W'}=3$.  Since $n\geq 7$,  it follows that there exists 
 a $6$-dimensional non-degenerate $g$-stable subspace $W_1$  of $V$  containing $W$. Set $Y_1=U(W_1)$, $Y_2=U(W_1^\perp)$, and let $Y$ be the stabilizer of $W_1$ in $G$. Clearly $Y\cong Y_1\times Y_2$. Let us write $g=y_1y_2$, where $y_1\in Y_1$ and $y_2\in Y_2$. Note that $|y_1|=9$.
Let $\sigma$ be an irreducible constituent of $\tau|_Y$  non-trivial on $Y_1'$.
Then $\sigma=\sigma_1\otimes \sigma_2$ and $\sigma(g)=\sigma_1(y_1)
\otimes \sigma_2(y_2)$. As $\tau(g)$ is almost cyclic, so are $\sigma(g)$ and
$\sigma_1(y_1)$. This contradicts Lemma \ref{u62}.

\begin{propo}\label{sy4}  Let $G=Sp(2n,q)$, where $n>1$, $q$ is odd and $G\neq Sp(4,3)$.
Let $g\in G$ be a non-scalar
semisimple orthogonally decomposable element of $p$-power order,
and let $\tau$ be an \ir Weil $F$-\rep
of $G$ of dimension greater
than $1$. Then $\tau(g)$ is not almost cyclic.
\end{propo}

Proof. (A) We first rule out the case where $q=3$, $\ell=2$ and $V$ contains a $g$-stable
$2$-dimensional non-degenerate subspace $X$, say. In this case,  $g$ is contained
in a subgroup $K=K_1\times K_2$, where $K_1= Sp(X)$ and $K_2= Sp(X^\perp) $.
In particular,  $ K_2$
is not solvable (otherwise $G=Sp(4,3)$).
Let $g=g_1g_2$, where $g_1\in K_1,g_2\in K_2$. As $g_1$ is semisimple
and $K_1/K_1'$ is of order 3, we have
 $g_1\in K_1$. As
 $K_2$ is perfect,   $g\in K'$.
By Lemma \ref{rd6}(2),
$\tau|_{K'}$ contains at least two isomorphic composition factors of dimension
greater than 1.
Therefore, $\tau(g)$ is not almost cyclic, and the claim follows.

\medskip
(B) Next, let us show that $g$ has no  $1$- or $-1$-eigenspace on $V$. Indeed,
suppose the contrary. Observe that any such
eigenspace is non-degenerate. So, clearly, in any of them we can choose  a non-degenerate  $g$-stable subspace
 $X$, say, of dimension 2. Thus, $g$ is contained
in a subgroup $H=G_1\times G_2$, where $G_1= Sp(X^\perp) $ and $G_2= Sp(X)$.
Suppose that $G_1$ is solvable.  Then $G=Sp(4,3)$, which is against our assumption.
So we may assume that $G_1$ is not solvable.
It follows, by Lemma \ref{rd6}, that there is an \ir constituent $\phi$ of $\tau|_H$ such that
$\phi=\phi_1\otimes \phi_2$, where $\phi_1\in\Irr G_1,$ $\phi_2\in\Irr G_2$ and
 either both $\phi_1,\phi_2$ are of dimension greater than 1, or
$q=3,\ell=2$. The latter case is ruled out in (A). 
In the former case $\phi(g)$ is obviously not almost cyclic.

\medskip
From now on,
we choose $W$ to be a non-degenerate $g$-stable  subspace of $V$ such that  $g|_W$ is
orthogonally indecomposable, and  $W$ is of maximal
dimension with this property.  Set $2m=\dim W$. So $g\in H=G_1\times G_2$, where
 $G_1\cong Sp(2m,q)$ and $G_2 \cong  Sp(2n-2m,q)$. Set $g=g_1g_2$,
where $g_1\in G_1,g_2\in G_2$.

\medskip
(C) Suppose  that  $G_1$ is solvable. Then $G_1=Sp(2,3)$,
and $g$ stabilizes a direct sum of non-degenerate two-dimensional subspaces of $V$
(by our choice of $W$). So $|g|\leq 4$.
We can assume $\ell\neq 2$ by (A).
If $g^2$ is not scalar, $g$ must have
a  $1$- or a $-1$-eigenspace. But this case has been ruled out in (B).
So we are left with the case where
  $g^2=\pm \Id$. As above, since $\ell\neq 2$, by Lemma \ref{rd6} there is an
\ir constituent $\phi$ of $\tau|_H$ such that
$\phi=\phi_1\otimes \phi_2$, where $\phi_1\in\Irr G_1,$ $\phi_2\in\Irr G_2$ and
 both $\phi_1,\phi_2$ are of dimension greater than 1. Therefore,
$\phi_i(g_i)^2=\pm \Id$
for $i=1,2$. However, the tensor product of any two matrices over $F$
of size greater than 1 whose squares are scalar cannot be almost cyclic (by Lemma \ref{p74}).

\med In view of the above, from now on we  may
 assume that $G_1$ is not solvable.

\med
$(D)$ Suppose  that $G_2\cong Sp(2n-2m,q)$ is solvable.
So $G_2=Sp(2,3)$.
By (B),  $g_2\neq \Id$, so $|g_2|$, and hence also $|g_1|$, is a non-trivial 2-power.
Moreover, by our choice of $W$, $g_1$ is orthogonally indecomposable.
By Lemma \ref{rd6}, there is an \ir constituent $\phi$ of $\tau|_{H}$ such that  $\phi=\phi_1\otimes \phi_2$, and $\dim\phi_i>1$ for $i=1,2$. Assume that $\tau(g)$ is almost cyclic. Then $\phi(g)$ is also almost cyclic, and therefore $\phi_i(g_i)$
is  cyclic for $i=1,2$, in view of Lemma \ref{p74}. As $g_1$ is orthogonally indecomposable,
$g_1$ is either a power of  a Singer cycle, or a power of a  Singer-type cycle
of $Sp(2r,3)$. It follows from Lemmas \ref{sp55} and \ref{cw1}
applied to $\phi_1(g_1)$, that
 $ g_1$ itself is either a Singer cycle or a Singer-type cycle.
By Lemma \ref{zgm},
 either $m=2$ and  $|g_1|=3^2-1=8$,  or $m=1$ and $G_1=Sp(2,3)$.
The  latter case is ruled out, as $G_1$ is non-solvable.
In the former case $G=Sp(6,3)$ and $G_1=Sp(4,3)$.  By (A) applied to $Sp(6,3)$,
we can assume $\ell\neq 2$.
As $g_1^4\in Z(G_1)$,  $\phi_1(g_1^4)$ is scalar.
   By Lemma \ref{423},
$\dim\phi_1(g_1)=4$ (as $\phi_1(g_1)$ is  cyclic) and the spectrum  of
$\phi_1(g_1)$ consists of four distinct $4$-roots of $-1$.
Denote this set by $S$, say. In turn, as
$G_2\cong Sp(2,3)$, we have $g_2^4=1$. Note that $S\cdot \al=S$
for every 4-root $\al$ of $1$. Therefore, $\phi(g)$ consists of
4-roots of $-1$, each of \mult equal to $\dim \phi_2$, a contradiction.

  \med
$(E)$ Suppose that $|g|$ divides $q+1$ or $q-1$ (it is convenient to consider this case separately).
 By Propositions \ref{GS1} and \ref{GS2}, $G$ can be
generated by at most $2n$ conjugates of $g$
unless $n=2$ and $g^2\in Z(G)$, in which case  $G$ can be generated by at most $5$ conjugates of  $g$.
Suppose that $
\tau(g)$ is almost cyclic. In the exceptional case, it follows from Lemma \ref{nd5} that $\dim\tau\leq 5$; as $\dim\tau\geq (q^2-1)/2$,
this implies $q=3$, which is ruled out by our assumptions. Otherwise, again by Lemma \ref{nd5},
$\dim \tau\leq 2n(|g|-1)\leq 2nq$. As $\dim\tau\geq (q^n-1)/2$, this implies
$q^n\leq 1+4nq$, whence either $n=2,q\leq 7$ or $n=3,q=3$.

Suppose that $n=2,q\leq 7$. Then either $p=2$ or $|g|=3$. In the latter case, the above bound  reduces to $\dim \tau\leq 8 $, whence $q^2\leq 17$. This implies $q=3$, which contradicts our assumptions. Now,  we are left with the cases
$p=2$ and $(n,q)\in\{ (2,5), (2,7), (3,3)\}$. Clearly, in view of $(B)$, $|g| \neq 2$. Suppose that $|g| = 4$. As $g^2\not\in Z(G)$, $g^2\neq -\Id$,
and hence $g$ acts as an element of order $2$ on the 1-eigenspace of $g^2$. But this contradicts $(B)$. Therefore, the case $(n,q) = (3,3)\}$ is ruled out, and are left with the case $|g|=8$, $G=Sp(4,7)$.

Thus, $G_1\cong G_2\cong SL(2,7)$.
The restriction of
$\tau$ to $H=G_1\times G_2$ is easy to describe for $G=Sp(4,q)$ and $G_1\cong G_2\cong Sp(2,q)$ (e.g. see \cite{Za85}, Theorem 2). Namely. let $\lam_i,\mu_i$ $(i=1,2)$ be the \ir Weil \reps of
$G_i$ of degree $(q-1)/2$ and $(q+1)/2$, respectively, over the complex numbers.
Then $\tau|_{H}= (\lam_1\otimes \mu_2) \oplus
(\lam_2\otimes \mu_1)$ if $\dim\tau=(q^2-1)/2$, whereas $\tau|_{H}= (\lam_1\otimes \mu_1)\oplus
(\lam_2\otimes \mu_2)$ if $\dim\tau=(q^2+1)/2$.  Recall that the representations
$\lam_i,\mu_i$ remain \ir  modulo any $\ell \neq2$, so these formulae are valid for any $\ell\neq 2$. Furthermore, as $q=7$,  we have $\dim\lam_i=3$ and $\dim\mu_i=4$. On the other hand, the representations $\mu_i$ under reduction mod 2 contain a composition factor isomorphic to $\lam_i$ mod 2. Thus, by Lemma \ref{p99}, we may assume that $\ell \neq 2$.

By Lemma \ref{p74}, $\lam_1(g_1)\otimes \mu_2(g_2)$ is cyclic if $\tau(g)$ is so. This implies $\deg \tau(g)\geq 12-3+1=10$, which is a contradiction, as
$|g|=8$. Similarly, we get a contradiction considering $\lam_1(g_1)\otimes \mu_1(g_2)$. This completes
the argument for $p=2$.

\medskip
(F) Finally, suppose that  both $G_1$ and $G_2$ are not solvable.
Then, again by Lemma \ref{rd6}, there is an \ir constituent $\phi$, say, of $\tau|_{H}$
such that $\phi=\phi_1\otimes \phi_2$, where $\phi_1,\phi_2$ are \ir Weil \reps of $G_1,G_2$, respectively,  both  of dimension at least 2.  Thus, $\phi(g)=\phi_1(g_1)\otimes \phi_2(g_2)$.
As in (D), assuming that $\tau(g) $, and hence $\phi(g)$, is almost cyclic, it follows from Lemma \ref{p74} that $\phi_1(g_1)$ and $\phi_2(g_2)$
are cyclic. In particular, $g_1$ and $g_2 $ are not scalar.

\med
(i) Assume first that $\ell=p$.  By Lemma \ref{p99}, $\phi(g)$ is not almost cyclic unless
$\ell\neq 2$ and $\dim\phi_i=2$ for $i=1,2$. As $\dim \phi_2\geq  (q^{n-m}-1)/2$ and
$n-m\geq 1$, the equality  $\dim\phi_2=2$ implies $n-m=1$ and $q=3$ or 5. As $G_2$ is not solvable, we are left to examine the case where $G_2 = Sp(2,5)$.
Similarly, $\dim\phi_1=2$ implies $G_1= Sp(2,5)$. In addition, $\ell=p=3$, as $p=\ell \neq 2,5$. Therefore, $G=Sp(4,5)$ and $|g|=3$.   Note that, by Lemma \ref{wr3},
$\tau|_{G_i}$ is a sum  of \ir Weil \reps of $G_i$. As $\ell=3$, none of them is one-dimensional.
(Indeed, the \ir Weil \reps of $Sp(2,5)$ are of dimension 2 and 3 in characteristic 0. Both of them remain \ir modulo 3, the former
by dimension reasons, and the latter by the fact that it is of $3$-defect 0.)  So $\tau(g)$ is not almost cyclic, unless $\dim \tau=\dim\phi=4$. However, $\dim\tau>4$.

\med
(ii) Now, assume that $\ell\neq p$.
 The case where $m=1$ is ruled out by (E).  Indeed, if $m=1$, then every orthogonally indecomposable $g$-stable subspace of $V$ is of dimension 2, and hence  we may choose $W$ so that $|g|=|g_1|\geq |g_2|$.
So, we may assume $m>1$. As $\phi_1(g_1)$ is cyclic,  $g_1=g|_W$ is orthogonally indecomposable on $W$, and $|g_1|$  is a $p$-power, it follows from Lemma \ref{sp55} 
that $\lan g_1,Z(G_1)\ran$ is of order $q^m\pm 1$.

Suppose first that
 $p>2$. Then $|g_1|=(q^m\pm 1)/2$. As $W$ is chosen of maximum dimension, we again get
$|g|= |g_1|\geq |g_2|$.  Recall that $\phi_1$ is an \ir Weil $F$-\rep of $G_1$ (Lemma \ref{wr3}),
and hence has dimension $(q^m-1)/2$ or $(q^m+1)/2$. As the matrix of $\phi_1(g_1)$ is  cyclic,
it has size $k=|g_1|$ or $|g_1|-1$, and is similar to $\diag(\ep_1\ld \ep_k)$,
where the $\ep_i$'s are pairwise distinct $|g_1|$-roots of unity.
On the other hand, $\phi(g)=\phi_1(g_1)\otimes \phi_2(g_2)$ has order at most $|g|\leq k+1$.
As $\phi(g)$ is almost  cyclic, this contradicts Lemma \ref{p74} (unless $m=1$, which is not the case here).

Next, let $p=2$ (and hence $\ell\neq 2$ by (i)).
Then $|g_1|=q^m\pm 1$ is a $2$-power.
 As  $m>1$, in view of Lemma \ref{zgm}  this implies that
 $q^{m}=3^2$, $|g_1|=8$, $G_1=Sp(4,3)$ and $G=Sp(2n,3)$. Let $t=g_1^{4}$ and $h=g^{4}$.
Then $t=-\Id$ 
and $h=\diag(t,t')$, where $t'=g_2^{4}$. Note that,
 as $\phi_1(g_1)$ is cyclic, we have $\dim\phi_1=4$ by Lemma \ref{423}; in turn,
this implies that $\phi_1(t)=-\Id$.

 Suppose first that  $n>3$. Then $\dim\phi_2\geq (3^2-1)/2=4$. By Lemma \ref{p74} (as $\phi_1(g_1),\phi_2(g_2)$ are cyclic), we have $\deg\phi(g)\geq 4^2-4+1=13$, which is false as $|g|=|g_1|=8$. Let $n=3$. Then $G_2=Sp(2,3)$ is solvable, which is false. This completes
the proof of the Proposition.

\begin{propo}\label{35} 
Let   $SL(n,q)\subseteq G\subseteq GL(n,q)$, where $n>2$ and $(n,q)\neq (4,2), (3,3)$ or $(4,3)$.
Let  $\phi $ be an \ir Weil $F$-\rep of  $G$, with
$\dim\phi>1$, and let $g \in G$ be a non-scalar semisimple element of prime-power order $p^a$ for some prime $p$. Then  $\phi (g)$ is almost cyclic if and only if  $g$ is irreducible and $|\lan g,Z(GL(n,q))\ran |=q^n-1$. 
\end{propo}

Proof. Observe first that the 'if' part of the statement follows from Lemma  \ref{ss3} and Corollary  \ref{bc44} (recall that, as $n>2$, the irreducible Weil $F$-representations of $G$ extend to $GL(n,q)$). So, from now on, we assume that $\phi (g)$ is almost cyclic.

 Let $V$ be the natural $G$-module, and let $W\subseteq V$ be a $g$-stable subspace of $V$ on which $g$ acts irreducibly,
and  such that $|g|$
coincides with the order of  $g|_{W}$ (observe that this choice is possible as $|g|$ is a prime-power).

If $V=W$,  then the statement  follows from
Lemma \ref{sp55}. Otherwise,  by
Lemma \ref{sz33}, $g$ stabilizes no one-dimensional subspace (and hence also no subspace of codimension 1, by Maschke's theorem).
Therefore, setting $\dim W = d$, we have
 $n-1> d >1$ (so $n>3$). Also, $|g|$ does not divide $q-1$ (otherwise $g|_{W}$ would be scalar). Furthermore, we
may assume that  $(d,q)\neq (2,2),(2,3)$. (Indeed, if $(d,q)= (2,2)$ then,
since $(n,q)\neq (4,2)$, we have $n\geq5$. Then $G$ is generated by at most $n$ conjugates of $g$ (by Proposition  \ref{GS1}).
As $\phi(g)$ is almost cyclic, and $|g| = |g|_W| = 3$, we have $\dim\phi\leq 2n$ (by Lemma \ref{nd5}). However, the lower bound for $\dim \phi$ is $2^n-n-1>2n$ for $n\geq5$ (see \cite{SZ}),
which is a contradiction. If $(d,q)=(2,3)$, then $|g| = |g|_W| \leq 8$, and hence $V$ is the direct sum of $2$-dimensional $g$-stable subspaces. It follows that $g$ stabilizes a direct sum of subspaces   $W'\oplus W''$,
where $\dim W'=4$,  $g|_{W'}\neq \Id$ and $\dim W''=n-4$. Then the claim follows from Lemma \ref{s43} for $(n,q)=(4,3)$.)  

Now, we can write $V=W\oplus V'$, where $V'$ is a $g$-stable subspace of $V$. Set $X=
\{x\in SL(V):x|_{V'}=\Id
\}$ and $Y=\lan X,g\ran$. Clearly, $X\cong SL(W)$. Let $g_1=\diag(g|_{W}, \Id)$ and $g_2=gg_1\up$
(note that $g_1,g_2$ may not belong to $G$, but $|g_1|=|g|$). Set $Y_1=\lan X,g_1\ran$. 
Then $Y\subset \lan Y,g_2\ran =Y_1\times \lan g_2\ran$.

Let $\tau$ be any \ir constituent of $\phi|_{Y}$. Then, by our assumption on $\phi(g)$,  $\tau(g)$ is almost cyclic. Since $g_2$ centralizes $Y$,
$\tau $ extends to a \rep $\tau'$, say, of $\lan Y,g_2\ran $. As   $\tau'(g_2)$ is scalar and $g_1=gg_2\up$,
the matrix of $\tau'(g_1)$ is almost cyclic. As $Y_1=\lan X,g_1\ran$, one observes that
$\tau'(Y_1)$ contains an almost cyclic matrix $\tau'(g_1)$. Set $X_1=Y_1|_{W} \cong Y_1$. Then we can view  $\tau'|_{Y_1}$ as a \rep of $X_1$. Note that $SL(W)\subset  X_1\subset GL(W)$. By Lemma \ref{wr3},
$\tau'|_X=\tau|_X$ is a Weil \rep of $X$, and hence so is $\tau'|_{X_1}$.

Therefore, we can apply results obtained earlier to $\tau'(X_1)$. Namely, by  Lemma \ref{p33}(2), Lemma \ref{cw1},
Lemma \ref{sp55} and the remark following it, it follows that either
$d=2$, or $d$ is an odd prime and $g_1$ is a multiple of $(q^d-1)/(q-1)$, which must be a $p$- power.
Moreover, if $d=2$, then the order of $g_1$ must divide $2(q+1)$, as $g$ stabilizes no line on $V$.

Suppose first that $d>2$, or $d=2$ and $p$ is odd.
In this case, $p$ is coprime to $q-1$ (otherwise, by Zsigmondy's theorem,  $(q^d-1)$ would be divisible by a prime different from $p$).
It follows that $\det g_1=1=\det g_2$, and hence $g$ is contained in $H:=X\times X_2$, where $X \cong SL(W)$ and $X_2\cong SL(V')$.
 As $(d,q) \neq (2,2),(2,3)$, by \cite[Corollary 3.8]{TZ08}, $\phi|_{H}$ contains an \ir constituent $\tau$
 that is non-trivial on both $X$ and $
X_2$. Let  $\tau=\tau_1\otimes \tau_2$, where $\tau_1 \in \Irr(X)$  and $\tau_2 \in \Irr(X_2)$.
Note that $\dim \tau_2>1$ unless  $q=2,3$ and $n-d=2$.  

 We need to examine the following cases:

(a)  $d>2$, $(n-d,q)\neq (2,2),(2,3)$. 
 In this case we may apply Lemma \ref{sp55} to $g_1$. Namely,   by Lemma \ref{sp55}
and Corollary  \ref{bc44}, $\deg\tau_1(g_1)\geq |g_1|-1$. Recall that $g_2$ is not scalar, as $g$ stabilizes no one-dimensional subspace, and hence, as  $SL(V')$ is quasi-simple, $\tau_2(g_2)$ is not scalar.
Therefore, $\deg\tau(g)=\deg (\tau_1(g_1)\otimes \tau_2(g_2))=|g|$. Note that the mappings $g\ra \tau_i(g_i)$, $i=1,2$, yield representations
of the group $\lan g \ran$. If $p\neq \ell$, then $\tau(g)$ is not almost cyclic by Lemma \ref{tp1}. So, let $p=\ell$. Then, by Lemma \ref{p99},
$\tau(g)$ is not almost cyclic unless $\ell\neq 2$ and $\dim\tau_1=\dim\tau_2=2$. However, this implies $d=2$, which  is not the case.

Next suppose that  $d>2$, $(n-d,q)\in\{ (2,2),(2,3)\}$. As $p$ is coprime to $q-1$ and $|GL(2,3)|=48$, the case $q=3$ is ruled out. So, let $q=2$.
 Note that $g_2$ is irreducible in $GL(n-d,q)=GL(2,q)$, as
$g$ stabilizes no line of $V$. As $\tau_1(g_1)$ is almost cyclic and $\tau_1$ is a Weil \rep of $X$,
we have $|g_1|=2^{n-2}-1$ by Lemmas \ref{sp55} and \ref{cw1}. 
As $g_1$ is irreducible on $W$, the eigenvalues of $g_1$ in $GL(n-2, \overline{F}_2)$ are pairwise distinct primitive $|g_1|$-roots of unity (by Galois theory), whereas
the eigenvalues of $g_2$ in $GL(2, \overline{F}_2)$ are distinct primitive $|g_2|$-roots of unity. Therefore, all the eigenvalues of  $g$
in $GL(n, \overline{F}_2)$ are pairwise distinct,  unless $|g_2|=|g_1|$. In the latter case $d=n-2=2$, which is not the case.

Thus, the eigenvalues of  $g$
in $GL(n, \overline{F}_2)$ are distinct, that is, $g$ a is regular semisimple element and
$C_G(g)$ has no unipotent element. By Lemma \ref{g12},
$G$ is generated by three conjugates of $g$. As  $\phi(g)$ is almost cyclic and $\deg\phi(g)\leq |g|$, by
Lemma \ref{nd5},
$\dim\phi\leq 3|g|-3$. As $|g|=|g_1|=2^{n-2}-1$, we have $\dim\phi\leq 6(2^{n-3}-1)$. On the other hand, the  dimension of an \ir
$F$-\rep of $G$ is at least 
$2^{n-1}-n-1$. This gives us a contradiction.

(b) $d=2$, $p>2$, $n\geq 4$. As $p$ is odd here, $|g|\leq q+1$, and $G$ is generated by at most $n$ conjugates of $g$ (see Propositions \ref{GS1} and \ref{GS2}).
Then $\dim \phi \leq nq$, whereas $\dim \phi \geq (q^n-1)/(q-1)-2$ as $\phi$ is a Weil representation. 
This is a contradiction.

(c) $d=p=2$, $n\geq 4$. Then $|g|\leq 2(q+1)$ and hence $\dim \phi \leq n(2q+1)$,
whereas $\dim \phi \geq (q^n-1)/(q-1) - 2$. This implies  $q=2$, and hence $p>2$, a contradiction.

 \medskip

\section{Proof of Theorem  \ref{mt1}}

At this stage we are in a position to prove Theorem \ref{mt1}. 


\medskip

{\bf Proof of Theorem} \ref {mt1}: Suppose first that $G=Sp(2n,q)$, $G \neq Sp(4,3)$. Then $\tau(g)$ cannot be almost cyclic unless $g$ is orthogonally indecomposable, by Proposition \ref{sy4}. So, assume that $g$ is orthogonally indecomposable. Assume first that $g$ is either a Singer or a Singer-type cycle. These cases are ruled out in Lemma  \ref{cw0}, items (1) and (2) respectively, applying Lemma  \ref{zgm}.  So, assume that $|g|$ is a proper divisor of $q^n +1$, $q^n -1$, respectively. Then the claims  (a) and (b) in item (1) of the statement follow from Lemma  \ref{cw1}. Next, suppose that $G=Sp(4,3)$. This group is dealt with in Lemma  \ref{423}, yielding the cases under (c) in item (1) of the statement.

Now, suppose that  $SU(n,q)\subseteq G\subseteq U(n,q)$,
 where  $n>2$ and $(n,q)\neq (5,2)$, $(4,2)$, $(3,3)$, $(3,2)$. Remember that $|g|$ is assumed to be a prime-power, and recall the Remark following Lemma \ref{sp55}. Take into account Lemma \ref{p33},  Corollary \ref{bc44} and Lemma \ref{cw1}.  Then claim (2),(a) of the statement follows from Lemma \ref{s3s}. Now, let us consider the 'exceptional cases' $(n,q) = (5,2),(4,2),(3,3),(3,2)$. The case $(n,q) = (5,2)$ is dealt with in Lemma \ref{452}, yielding item (2),(b) of the statement.  Next, suppose that $(n,q)=(4,2)$. Observe that we may assume that $G=SU(4,2)$. Then the claims in item (2),(c) of the statement follow from Lemma \ref{442}. Now, suppose that $(n,q)=(3,3)$. If $g$ is orthogonally indecomposable, then Lemma \ref{p33} and Lemma \ref{cw1} yield $|g|=7$. Otherwise $g$ is orthogonally decomposable, and Lemma \ref{pr1} applies. So the claims in item (2),(d) of the statement hold. Finally, suppose that $(n,q)=(3,2)$. This case has been handled in detail by direct computation (see the Remark following Lemma \ref{cw1}), yielding item (2),(e) of the statement.
 
 We are left with the case where $SL(n,q)\subseteq G\subseteq GL(n,q)$, with $n>2$. Suppose first that $(n,q)\notin \{(3,3),(4,3),(4,2)\}$. It then follows from Proposition \ref{35} that $\tau(g)$ is almost cyclic if and only if $g$ is irreducible and $|\lan g,Z(GL(n,q))\ran |=q^n-1$. This yields items (3),(a) and (3),(b) of the statement, according to Lemma \ref{p33} and Lemma \ref{cw1}. Now, suppose that $(n,q) = (3,3)$. This case is dealt with in Lemma \ref{333}, which gives $|g|=13$, that is an instance of item (3),(b) of the statement.
 Next, suppose that $(n,q) = (4,3)$. Then $\tau(g)$ is not almost cyclic, by Lemma \ref{s43}. Finally, observe that the case $(n,q)=(4,2)$ is ruled out by Lemma \ref{ss3}, (2), as we are assuming that $\tau$ is Weil, and hence has degree $> 7$. 

\bigskip
ACKNOWLEDGEMENT.  We wish to thank Marco Antonio Pellegrini for having helped us with his skills in the use of the MAGMA and GAP packages for dealing with character tables and Brauer character tables, and for having provided efficient ad hoc routines for testing cyclicity and almost cyclicity of matrices.

\

\bigskip
Authors' addresses: Dipartimento di Matematica e Applicazioni, Universita degli Studi di Milano-Bicocca,
Via R. Cozzi 53, Milano, 20125, Italy

e-mail:
lino.dimartino@unimib.it,
alexandre.zalesskii@gmail.com

\end{document}